\documentclass[12pt,leqno]{article}
\usepackage[mathscr]{eucal}
\usepackage{amsmath,amssymb,latexsym,theorem,bbm,amsfonts}
\usepackage{verbatim}
\usepackage{color}
\usepackage{appendix}
\usepackage{array,epsfig,theorem,bbm,graphics,booktabs,url}

\newcommand{\CC}{\mathbb{C}}

\newcommand{\NN}{\mathbb{N}}

\newcommand{\RR}{\mathbb{R}}

\newcommand{\ZZ}{\mathbb{Z}}

\newcommand{\bone}{{\boldsymbol{1}}}

\newcommand{\cA}{{\mathcal A}}

\newcommand{\cD}{{\mathcal D}}

\newcommand{\dd}{\mathrm{d}}
\newcommand{\ee}{\mathrm{e}}
\newcommand{\ff}{\mathrm{f}}
\newcommand{\ii}{\mathrm{i}}

\newcommand{\EE}{\operatorname{\mathbb{E}}}

\newcommand{\PP}{\operatorname{\mathbb{P}}}

\newcommand{\OO}{\operatorname{O}}

\renewcommand{\Re}{\operatorname{Re}}

\newcommand{\tS}{\widetilde{S}}

\newcommand{\vare}{\varepsilon}

\renewcommand{\mid}{\,|\,}

\renewcommand{\leq}{\leqslant}
\renewcommand{\geq}{\geqslant}

\newcommand{\distr}{\stackrel{\cD}{\longrightarrow}}
\newcommand{\distrf}{\stackrel{\cD_\ff}{\longrightarrow}}
\newcommand{\distre}{\stackrel{\cD}{=}}

\newcommand{\proofend}{\hfill\mbox{$\Box$}}

\numberwithin{equation}{section}

\theoremstyle{change} \theorembodyfont{\em}
\newtheorem{Lem}{Lemma.}[section]
\newtheorem{Thm}[Lem]{Theorem.}
\newtheorem{Pro}[Lem]{Proposition.}

\theorembodyfont{\rm}
\newtheorem{Rem}[Lem]{Remark.}

\begin{document}

\begin{center}
 {\bfseries\Large
   On simultaneous limits for aggregation of stationary
   randomized INAR(1) processes with Poisson innovations}

\vskip0.5cm

 {\sc\large
  M\'aty\'as $\text{Barczy}^{*,\diamond}$,
  Fanni K. $\text{Ned\'enyi}^{*}$,
  Gyula $\text{Pap}^{**}$}

\end{center}

\vskip0.1cm

\centerline{\sl \large To the memory of Gyula}

\vskip0.2cm

\noindent
 * MTA-SZTE Analysis and Stochastics Research Group,
   Bolyai Institute, University of Szeged,
   Aradi v\'ertan\'uk tere 1, H--6720 Szeged, Hungary.

\noindent
 ** Bolyai Institute, University of Szeged,
     Aradi v\'ertan\'uk tere 1, H--6720 Szeged, Hungary.

\noindent e--mails: barczy@math.u-szeged.hu (M. Barczy),
                    nfanni@math.u-szeged.hu (F. K. Ned\'enyi).

\noindent $\diamond$ Corresponding author.

\vskip0.2cm


\renewcommand{\thefootnote}{}
\vspace*{0.2cm}
\footnote{M\'aty\'as Barczy is supported by the J\'anos Bolyai Research Scholarship of the Hungarian Academy of Sciences,
 and by the \'UNKP-19-4 New National Excellence Program of the Ministry for Innovation and Technology.
Fanni K. Ned\'enyi is supported by the \'UNKP-19-3 New National Excellence Program of the Ministry for Innovation and Technology.
Gyula Pap was supported by the Ministry for Innovation and Technology, Hungary grant TUDFO/47138-1/2019-ITM.}

\vspace*{-10mm}

\begin{abstract}
We investigate joint temporal and contemporaneous aggregation of \ $N$ \ independent
 copies of strictly stationary INteger-valued AutoRegressive processes of
 order 1 (INAR(1)) with random coefficient \ $\alpha \in (0, 1)$ \ and with idiosyncratic Poisson
 innovations.
Assuming that \ $\alpha$ \ has a density function of the form
 \ $\psi(x) (1 - x)^\beta$, \ $x \in (0, 1)$, \ with \ $\beta\in(-1,\infty)$ \ and
 \ $\lim_{x\uparrow 1} \psi(x) = \psi_1 \in (0, \infty)$, \ different limits of
 appropriately centered and scaled aggregated partial sums are shown to exist for
 \ $\beta\in(-1,0]$ \ in the so-called simultaneous case, i.e., when both \ $N$ \ and the time scale \ $n$ \ increase to infinity at a given rate.
The case \ $\beta\in(0,\infty)$ \ remains open.
We also give a new explicit formula for the joint characteristic functions of finite dimensional distributions of the appropriately
 centered aggregated process in question.
\end{abstract}

\section{Introduction and main results}\label{Sec_Int_Main}

Studying temporal and contemporaneous (also called cross-sectional) aggregations of independent stationary stochastic processes
 goes back to Granger \cite{Gra}.
He started to investigate contemporaneous aggregation of
 random-co\-ef\-fi\-cient autoregressive processes of order 1 in order to obtain the long memory phenomenon in aggregated time series.
Random-coefficient autoregressive processes of order 1 were introduced by Robinson \cite{Rob}, and some of its statistical properties were studied as well.
The field of aggregation of stochastic processes became an important area of statistics,
 for surveys on aggregation of different kinds of stochastic processes, see, e.g., Pilipauskait{\.e} and Surgailis \cite{PilSur},
 Jirak \cite[page 521]{Jir} or the arXiv version \cite{BarNedPap_Arxiv} of Barczy et al.\ \cite{BarNedPap}.
For historical fidelity, we note that Theil \cite{The} already considered contemporaneous aggregations of linear regression models
 with non-random coefficients, and later Zellner \cite{Zel} investigated the case of random coefficients.

Recently, Puplinskait{\.e} and Surgailis \cite{PupSur1, PupSur2} have studied iterated aggregation of random coefficient autoregressive processes of order 1
 with common innovations and with so-called idiosyncratic innovations, respectively, belonging to the domain of attraction of an $\alpha$-stable law.
They described the weak limits of finite dimensional distributions of appropriately centered and scaled aggregated partial
 sum processes when first the number of copies \ $N \to \infty$ \ and
 then the time scale \ $n \to \infty$.
\ Very recently, Pilipauskait{\.e} et al.\ \cite{PilSkoSur} have extended the results of Puplinskait{\.e} and Surgailis \cite{PupSur2} (idiosyncratic case)
 deriving the weak limits of finite dimensional distributions of appropriately centered and scaled aggregated partial sum processes
 when first the time scale \ $n \to \infty$ \ and then the number of copies \ $N \to \infty$, \ and when \ $n \to \infty$ \ and
 \ $N \to \infty$ \ simultaneously with possibly different rates.
We note that similar kinds of results were derived for the total accumulated work process of
 the aggregation (also called superposition) of independent and identically distributed stationary ON/OFF processes,
 see, e.g., Taqqu et al.\ \cite{TaqWilShe}, Mikosch et al.\ \cite{MikResRooSte} and Dombry and Kaj \cite{DomKaj}.
Namely, there are two kinds of results, iterated ones and simultaneous ones for the total accumulated work process in question:
 first the number of aggregated processes \ $M$ \ tends to infinity and then the time-scaling parameter \ $t$ \ converges to infinity,
 and in reversed order (iterated cases), and the simultaneous cases in which both \ $M$ \ and \ $t$ \ go to infinity at the same time
 possibly at different rates.
In the simultaneous cases it turned out that there are three subcases, where so-called fast, slow and intermediate growth conditions hold, respectively,
 see, e.g., Dombry and Kaj \cite[page 35]{DomKaj}.
In Pilipauskait{\.e} and Surgailis \cite[page 1013]{PilSur}, one can find a comparison of their results on aggregation
 of random coefficient autoregressive processes of order 1 and the above mentioned results on the total accumulated work process for ON/OFF processes.
For some random coefficient autoregressive processes of order 1,
 Leipus et al.\ \cite{LeiPhiPilSur} have also described the asymptotic behaviour of sample covariances in \ $N \times n$ \ panel data
 (see formula (1.5) in \cite{LeiPhiPilSur}) when both \ $N$ \ and \ $n$ \ tend to \ $\infty$, \ possibly at different rate.

The above mentioned references are mainly about aggregation schemes for randomized autoregressive processes and ON/OFF processes.
In the present paper we study aggregation procedures for randomized INteger-valued Autoregressive Processes of order 1 (INAR(1) processes)
 with Poisson innovations in the so-called simultaneous case, and this work can be considered as a continuation of the papers
 Barczy et al.\ \cite{BarNedPap} and Ned\'enyi and Pap \cite{NedPap}, where the iterated cases have been studied.
According to our knowledge, simultaneous limits have not been derived for aggregations of randomized INAR(1) processes
 (or more generally for those of randomized branching processes with immigration), our results are the first ones
 in this direction.
In our forthcoming Theorems \ref{simultaneous_aggregation_random_2} and \ref{simultaneous_aggregation_random_0}
 the number of aggregated copies of a stationary randomized INAR(1) process with Poisson innovations
 and the time scale both tend to \ $\infty$ \ simultaneously at a rate which could be considered
 analogous to the fast growth condition for ON/OFF processes mentioned above.

Let \ $\ZZ_+$, \ $\NN$, \ $\RR$, \ $\RR_+$, \ and \ $\CC$ \ denote the set of
 non-negative integers, positive integers, real numbers, non-negative real numbers,
 and complex numbers, respectively.
For \ $x,y\in\RR$, \ let \ $x \vee y:=\max(x,y)$.
\ We will use \ $\distrf$ \ for the weak convergence of the finite dimensional distributions of stochastic processes
 with sample paths in \ $D(\RR_+,\RR)$, \ where \ $D(\RR_+,\RR)$ \ denotes the space of real-valued
 c\`adl\`ag functions defined on \ $\RR_+$.
\ Equality in distribution will be denoted by \ $\distre$.

An INAR(1) time series model was first introduced by McKenzie \cite{McK} and Al-Osh and Alzaid \cite{AloAlz}, and it
 is a stochastic process \ $(Y_k)_{k\in\ZZ_+}$
 \ satisfying the recursive equation
 \begin{equation}\label{INAR1}
  Y_k = \sum_{j=1}^{Y_{k-1}} \xi_{k,j} + \vare_k, \qquad k \in\NN ,
 \end{equation}
 where \ $(\vare_k)_{k\in\NN}$ \ are independent and identically distributed (i.i.d.) non-negative integer-valued
 random variables, \ $(\xi_{k,j})_{k,j\in\NN}$ \ are i.i.d.\ Bernoulli random
 variables with mean \ $\alpha \in (0, 1)$, \ and \ $Y_0$ \ is a non-negative
 integer-valued random variable such that \ $Y_0$,
 \ $(\xi_{k,j})_{k,j\in\NN}$ \ and \ $(\vare_k)_{k\in\NN}$ \ are
 independent, and we define \ $\sum_{j=1}^0:=0$.
\ With the binomial thinning operator \ $\alpha\,\circ$ \ due to Steutel and van Harn
 \cite{SteHar},
the INAR(1) model in \eqref{INAR1} can be written as
 \begin{equation*}
  Y_k = \alpha \circ Y_{k-1} + \vare_k , \qquad k \in\NN ,
 \end{equation*}
which is very similar to an autoregressive model of order 1 (where \ $\circ$ \ is replaced by the usual multiplication).
An INAR(1) process can also be considered as a special branching process
 with immigration having Bernoulli offspring distribution.
We point out the fact that the theory and application of integer-valued time series models (such as INAR(1) processes)
 are rapidly developing and important fields (see, e.g., the survey paper of Wei{\ss} \cite{Wei} and Chapter 5
 in the book of Turkman et al.\ \cite{TurScoBer}).

As in Barczy et al.\ \cite{BarNedPap}, we will consider a certain randomized INAR(1) process \ $(X_k)_{k\in\ZZ_+}$ \ with randomized thinning parameter
 \ $\alpha$, \ given formally by the recursive equation
 \begin{equation*}
  X_k = \alpha \circ X_{k-1} + \vare_k , \qquad k \in\NN ,
 \end{equation*}
 where \ $\alpha$ \ is a random variable with values in \ $(0, 1)$ \  and \ $X_0$ \ is
 some appropriate random variable.
We will construct a process \ $(X_k)_{k\in\ZZ_+}$ \ such that, conditionally on \ $\alpha$, \ it is a strictly stationary INAR(1)
 process with thinning parameter \ $\alpha$ \ and with Poisson immigrations.
Conditionally on \ $\alpha$, \ the i.i.d.\ innovations
 \ $(\vare_k)_{k\in\NN}$ \ have a Poisson distribution with
 parameter \ $\lambda \in (0, \infty)$, \ and the conditional distribution of the initial
 value \ $X_0$ \ given \ $\alpha$ \ is the unique stationary distribution,
 namely, a Poisson distribution with parameter \ $\lambda/(1-\alpha)$.
\ More precisely, let \ $\lambda \in (0, \infty)$, \ and let \ $\PP_\alpha$ \ be a probability
 measure on \ $(0, 1)$.
\ Then there exist a probability space \ $(\Omega, \cA, \PP)$, \ a random variable
 \ $\alpha$ \ with distribution \ $\PP_\alpha$ \ and random variables
 \ $\{X_0, \, \xi_{k,j}, \, \vare_k : k, j \in \NN\}$, \ conditionally independent
 given \ $\alpha$ \ on \ $(\Omega, \cA, \PP)$ \ such that
 \begin{align*}
  &\PP(\xi_{k,j} = 1 \mid \alpha)
  = \alpha = 1 - \PP(\xi_{k,j} = 0 \mid \alpha) , \qquad
  k, j \in \NN , \\
  &\PP(\vare_k = \ell \mid \alpha)
  = \frac{\lambda^\ell}{\ell!} \ee^{-\lambda} , \qquad \ell \in \ZZ_+ ,
  \qquad k \in \NN ,  \\
  &\PP(X_0 = \ell \mid \alpha)
  = \frac{\lambda^\ell}{\ell!(1-\alpha)^\ell} \ee^{-(1-\alpha)^{-1}\lambda} ,
  \qquad \ell \in \ZZ_+ ,
 \end{align*}
 for details see Barczy et al.\ \cite[Section 4]{BarNedPap_Arxiv}.
Note that the conditional distribution of \ $\vare_k$, \ $ k\in\NN$, \ does not depend on \ $\alpha$.
\ Define a process \ $(X_k)_{k\in\ZZ_+}$ \ by
 \[
  X_k = \sum_{j=1}^{X_{k-1}} \xi_{k,j} + \vare_k , \qquad k \in \NN .
 \]
Then, conditionally on \ $\alpha$, \ the process \ $(X_k)_{k\in\ZZ_+}$ \ is a strictly stationary INAR(1) process
 with thinning parameter \ $\alpha$ \ and with Poisson immigrations having parameter \ $\lambda$,
 \ see, e.g., Turkman et al.\ \cite[Section 5.2.1]{TurScoBer}.
\ The process \ $(X_k)_{k\in\ZZ_+}$ \ can be called a randomized INAR(1)
 process with Poisson immigrations, and the distribution of \ $\alpha$ \ is the
 so-called mixing distribution of the model.
We note that \ $(X_k)_{k\in\ZZ_+}$ \ is a strictly stationary sequence, but it
 is not even a Markov chain (so it is not an INAR(1) process) if \ $\alpha$ \ is not
 degenerate, see Section 2 and Appendix A in Barczy et al.\ \cite{BarNedPap_Arxiv}.
Further, a strong law of large numbers does not hold for \ $(X_k)_{k\in\ZZ_+}$ \ in the sense that
 \ $\frac{1}{n}\sum_{k=0}^n X_k$ \ does not converge to a constant as \ $n\to\infty$ \ with probability one.

The conditional generator function of \ $X_0$ \ given \ $\alpha \in (0, 1)$ \ takes the form
 \[
  F_0( z_0\mid \alpha)
  := \EE( z_0^{X_0} \mid \alpha)
  = \ee^{(1-\alpha)^{-1}\lambda({z_0}-1)} ,
  \qquad z_0 \in D,
 \]
 where \ $D:=\{z\in\CC : \vert z\vert\leq 1 \}$, \ i.e., conditionally on \ $\alpha$, \ $X_0$ \ has a Poisson distribution with parameter \ $(1-\alpha)^{-1}\lambda$, \
 and consequently the conditional expectation of \ $X_0$ \ given \ $\alpha$ \ is
 \ $\EE(X_0 \mid \alpha) = (1-\alpha)^{-1}\lambda$.
\ Here and hereinafter the conditional expectation \ $\EE(X_0 \mid \alpha)$ \ is meant in the
 generalized sense, see, e.g., in Stroock \cite[\S\,5.1.1]{Str}.
Then, as the negative part of \ $X_0$ \ is \ $0$, \ which is integrable, the conditional expectation in question
 does exist in this generalized sense.
The joint conditional generator function of \ $X_0, X_1, \ldots, X_k$ \ given
 \ $\alpha$ \ will be denoted by \ $F_{0,\ldots,k}(z_0, \ldots, z_k\mid\alpha)$,
 \ $z_0, \ldots, z_k \in D$.
\ Let us remark that the choice of Poisson-distributed innovations serves a technical purpose.
It allows us to calculate explicitly the stationary distribution of the model and also the joint characteristic function
of finite dimensional distributions of the randomized process itself (see Proposition \ref{Pro_egyuttes_char_aggregalt}).

Following the setup of our former paper Barczy et al.\ \cite{BarNedPap},
 we assume that the distribution of the random variable \ $\alpha$, \ i.e., the mixing distribution, has a
 probability density of the form
 \begin{equation}\label{alpha}
  \psi(x) (1 - x)^\beta , \qquad x \in (0, 1) ,
 \end{equation}
 where \ $\psi$ \ is a function on \ $(0, 1)$ \ having a limit
 \ $\lim_{x\uparrow 1} \psi(x) = \psi_1 \in (0, \infty)$.
\ This is the same mixing distribution as the one in Pilipauskait{\.e} and Surgailis \cite[equation (1.5)]{PilSur}
 used for randomized autoregressive processes of order 1.
Note that necessarily \ $\beta \in (-1, \infty)$ \ (otherwise
 \ $\int_0^1 \psi(x) (1 - x)^\beta \, \dd x = \infty$) \ and
 the function \ $(0,1)\ni x\mapsto \psi(x)$ \ is integrable on \ $(0,1)$.
\ Further, in case of
 \ $\psi(x) = \frac{\Gamma(a+\beta+2)}{\Gamma(a+1)\Gamma(\beta+1)} x^a$,
 \ $x \in (0, 1)$,  \ with some \ $a \in (-1, \infty)$, \ the random variable
 \ $\alpha$ \ is Beta distributed with parameters \ $a + 1$ \ and \ $\beta + 1$. \
This is an important special case from the historical point of view,
 since the Nobel prize winner Clive W. J. Granger used the square root of a
 Beta distribution as a mixing distribution for independent random coefficient AR(1) processes,
 and considered their contemporaneous aggregations, see Granger \cite{Gra}.
Note also that certain \ $\circ$ \ operators, where the summands are random parameter Bernoulli distributions
 with a parameter having Beta distribution, appear in catastrophe models.
One can check that, under \eqref{alpha}, for each \ $\ell\in\RR$, \ the expectation \ $\EE\big(\frac{1}{(1-\alpha)^\ell}\big)$ \ is finite
 if and only if \ $\beta > \ell-1$ \ (see, e.g., Barczy et al.\ \cite[Remark 4.5]{BarNedPap}).

Let \ $\alpha^{(j)}$, \ $j \in \NN$, \ be a sequence of independent copies of the
 random variable \ $\alpha$ \ having density function given in \eqref{alpha},
 and let \ $(X^{(j)}_k)_{k\in\ZZ_+}$, \ $j \in \NN$, \ be a sequence of independent copies of the process \ $(X_k)_{k\in\ZZ_+}$
 \ with idiosyncratic innovations (i.e., the innovations
 \ $(\vare^{(j)}_k)_{k\in\ZZ_+}$, $j\in\NN$, \ belonging to
 \ $(X^{(j)}_k)_{k\in\ZZ_+}$, \ $j \in \NN$, \ are independent) such that
 \ $(X^{(j)}_k)_{k\in\ZZ_+}$ \ conditionally on \ $\alpha^{(j)}$ \ is a strictly
 stationary INAR(1) process with thinning parameter \ $\alpha^{(j)}$ \ and
 with Poisson innovations having parameter \ $\lambda$ \ for all \ $j \in \NN$.

For each \ $N, n \in \NN$, \ consider the stochastic process
 \ $S^{(N,n)} = (S_t^{(N,n)})_{t\in\RR_+}$ \ given by
 \[
   S_t^{(N,n)}
   := \sum_{j=1}^N \sum_{k=1}^{\lfloor nt \rfloor}
       (X^{(j)}_k - \EE(X^{(j)}_k \mid \alpha^{(j)}))
   =  \sum_{j=1}^N \sum_{k=1}^{\lfloor nt \rfloor}
            \left( X^{(j)}_k - \frac{\lambda}{1-\alpha^{(j)}}\right),
   \qquad t\in\RR_+.
 \]
We remark that if \ $\beta\in(-1,0]$, \ then the first moment of \ $\frac{1}{1-\alpha}$ \ is infinite,
 so the centralization \ $\EE(X^{(j)}_k \mid \alpha^{(j)})$ \ in \ $ S^{(N,n)}$ \ could not be replaced
 by \ $\EE(X^{(j)}_k)$ \ in case of \ $\beta\in(-1,0]$.
\ From a statistical point of view, it is also reasonable to consider a process similar to \ $S^{(N,n)}$ \ given by
 \ $\widehat S^{(N,n)}_t := \sum_{j=1}^N \sum_{k=1}^{\lfloor nt \rfloor}\big(X^{(j)}_k -   \frac{\sum_{\ell=1}^n X^{(j)}_\ell}{n}\big)$, \ $t\in\RR_+$,
 \ which does not require the conditional expectations of the processes \ $X^{(j)}$, \ $j\in\NN$.

An INAR(1) process may be used to model migration,
 which is an important task nowadays all over the world.
More precisely, given a camp, for all \ $k\in\ZZ_+$, \ the random variable \ $X_k$ \ can be
 interpreted as the number of migrants present in the camp at time $k$, and every migrant
 stays in the camp with probability \ $\alpha\in (0,1)$ \ independently of each other
 (i.e., with probability \ $1 - \alpha$ \ each migrant leaves the camp) and,
 at any time \ $k\in\NN$, \ new migrants may come to the camp.
Given several camps in a country, we may suppose that the corresponding INAR(1) processes are independent,
 and each one can have independent parameters \ $\alpha$ \ coming from a certain distribution
 (in our case having a density function given in \eqref{alpha}).
So, the temporal and contemporaneous aggregates of these INAR(1) processes is the total usage of the camps
 in terms of the number of migrants in the given country in a given time period, and this
 quantity may be worth studying.

In Barczy et al.\ \cite{BarNedPap} and Ned\'enyi and Pap \cite{NedPap} limit theorems for appropriately scaled versions of
 \ $S^{(N,n)}$ \ have been derived in the so-called iterated cases, i.e., first taking the limit \ $N \to \infty$ \
 and then \ $n \to \infty$ \ or vica versa for all possible \ $\beta\in(-1,\infty)$.
\ (We note that in \cite{BarNedPap} and \cite{NedPap}, \ $S_t^{(N,n)}$ \ was denoted by \ $\tS_t^{(N,n)}$.)
As the main result of the paper, in case of \ $\beta\in(-1,0]$, \ we derive limit theorems
for appropriately scaled versions of \ $S^{(N,n)}$ \ in the so-called simultaneous case, i.e.,
 when both \ $N$ \ and \ $n$ \ increase to infinity at a given rate.
The case \ $\beta\in(0,\infty)$ \ is of the greatest interest, but it remains open,
 since our present technique is not suitable for this case (for more details, see Remark \ref{Rem_discussion2}).

\begin{Thm}\label{simultaneous_aggregation_random_2}
If \ $\beta \in (-1, 0)$, \ then
 \[
   n^{-1} N_n^{-\frac{1}{2(1+\beta)}} \,
   S^{(N_n,n)}
   \distrf (V_{2(1+\beta)} t)_{t\in\RR_+} \qquad
   \text{as \ $n \to \infty$ \ and \ $N_n^{\frac{-\beta}{1+\beta}} n^{-1} \to \infty$,}
 \]
 where \ $V_{2(1+\beta)}$ \ is a symmetric \ $2(1+\beta)$-stable random variable
 (not depending on \ $t$) \ with characteristic function
 \[
   \EE(\ee^{\ii \theta V_{2(1+\beta)}})
   =\ee^{-K_\beta|\theta|^{2(1+\beta)}}, \qquad \theta \in \RR ,
 \]
 where \ $ K_\beta := \psi_1(\frac{\lambda}{2})^{1+\beta} \frac{\Gamma(-\beta)}{1+\beta}$.
\end{Thm}

We note that Theorem \ref{simultaneous_aggregation_random_2} can be considered as a counterpart of
  Theorem 4.8 in Barczy et al.\ \cite{BarNedPap}, which is about the iterated aggregation case first taking the limit \ $N\to\infty$ \
  and then \ $n\to\infty$ \ in case of \ $\beta\in(-1,0)$.
The scaling factors and the limit processes coincide in these two theorems.
Heuristically, one might think that it is a consequence of the condition \ $N_n^{\frac{-\beta}{1+\beta}} n^{-1} \to \infty$ \
 as \ $n\to\infty$ \ in Theorem \ref{simultaneous_aggregation_random_2}, which in case of \ $\beta\in(-\frac{1}{2},0)$ \
 can be interpreted in a way that \ $N_n$ \ tends to \ $\infty$ \ much faster than \ $n$
 \ (indeed, \ $N_nn^{-1} = N_n^{\frac{1+2\beta}{1+\beta}} N_n^{-\frac{\beta}{1+\beta}} n^{-1}\to\infty\cdot\infty = \infty$ \ as \ $n\to\infty$ \
 in case of \ $\beta\in(-\frac{1}{2},0)$).
\ So this simultaneous case is more or less the above mentioned iterated case.
We mention that the same phenomenon occurs for randomized autoregressive processes of order (1), see Pilipauskait{\.e} et al.\ \cite[(2.15) and (2.20)]{PilSkoSur}.

\begin{Thm}\label{simultaneous_aggregation_random_0}
If \ $\beta=0$, \ then
 \[
   n^{-1} (N_n\log N_n)^{-\frac{1}{2}} \,
   S^{(N_n,n)}
   \distrf (W_{\lambda \psi_1} t)_{t\in\RR_+} \qquad
   \text{as \ $n \to \infty$ \ and \ $(\log N_n)^2 n^{-1} \to \infty$,}
 \]
 where \ $W_{\lambda \psi_1}$ \ has a normal distribution with mean \ $0$ \ and variance \ $\lambda\psi_1$.
\end{Thm}

We note that Theorem \ref{simultaneous_aggregation_random_0} can be considered as a counterpart of
  Theorem 4.9 in Barczy et al.\ \cite{BarNedPap}, which is about the iterated aggregation case first taking the limit \ $N\to\infty$ \
  and then \ $n\to\infty$ \ in case of \ $\beta=0$.
The scaling factors and the limit processes coincide in these two theorems.
For this fact one might give a similar heuristic explanation as we did in case of Theorem \ref{simultaneous_aggregation_random_2}
 (indeed, \ $N_nn^{-1} = (\log N_n)^{-2} N_n (\log N_n)^2 n^{-1} \to\infty\cdot\infty = \infty$ \ as \ $n\to\infty$),
 \ and note also that the same phenomenon occurs for randomized autoregressive processes of order (1), see Pilipauskait{\.e}
 et al.\ \cite[(2.16) and (2.21)]{PilSkoSur}.

In both Theorems \ref{simultaneous_aggregation_random_2} and \ref{simultaneous_aggregation_random_0} the limit processes are lines with random slopes.
So, similarly as it was explained at the end of Section 4 in Barczy et al.\ \cite{BarNedPap}, under
 the assumptions of Theorems \ref{simultaneous_aggregation_random_2} and \ref{simultaneous_aggregation_random_0}
 we have \ $n^{-1} N_n^{-\frac{1}{2(1+\beta)}} \, \widehat S^{(N_n,n)} \distrf 0$ \ as \ $n\to\infty$ \ and
 \ $n^{-1}(N_n\log N_n)^{-\frac{1}{2}} \, \widehat S^{(N_n,n)} \distrf 0$ \ as \ $n\to\infty$, \ respectively.
In principle, by applying some smaller scaling factors, one could try to achieve a non-degenerate weak limit of
 \ $\widehat S^{(N_n,n)}$ \ as \ $n\to\infty$ \ in these cases.

The proofs of Theorems \ref{simultaneous_aggregation_random_2} and \ref{simultaneous_aggregation_random_0} are based on an explicit formula for
 the joint characteristic function of \ $(S^{(1,n)}_{t_1},\ldots,S^{(1,n)}_{t_m})$, \ where
 \ $n,m\in\NN$ \  and \ $0=:t_0< t_1<t_2<\cdots<t_m$.
\ In fact, we derive two formulae for the characteristic function in question in the next proposition.

\begin{Pro}\label{Pro_egyuttes_char_aggregalt}
Let \ $n,m\in\NN$ \  and \ $0=t_0< t_1<t_2<\cdots<t_m$.
\ Then the joint characteristic function of \ $(S^{(1,n)}_{t_1},\ldots,S^{(1,n)}_{t_m})$ \ takes the form
 \begin{align}\label{egyuttes_char_aggregalt1}
  \EE\left(\exp\left\{ \ii \sum_{\ell=1}^m \theta_\ell S^{(1,n)}_{t_\ell} \right\}\right)
    = \int_0^1 \ee^{\frac{\lambda}{1-a} K_n(a)} \psi(a) (1-a)^\beta\,\dd a
 \end{align}
 for all \ $\theta_i\in\RR$, \ $i=1,\ldots,m$, \ where for all \ $a\in[0,1]$,
 \begin{align*}
   &K_n(a) :=  \sum_{\ell=1}^m
                \Bigl(\ee^{\ii \theta_{\ell,m}} - 1 - \ii \theta_{\ell,m}\Bigr)
                (\lfloor nt_\ell\rfloor - \lfloor nt_{\ell-1}\rfloor) \\
   &+\sum_{1\leq\ell_1<\ell_2\leq m}
      \sum_{k_1=\lfloor nt_{\ell_1-1}\rfloor+1}^{\lfloor nt_{\ell_1}\rfloor}
       \sum_{k_2=\lfloor nt_{\ell_2-1}\rfloor+1}^{\lfloor nt_{\ell_2}\rfloor}
        a^{k_2-k_1}
        \bigl(\ee^{\ii \theta_{\ell_1 , m}} - 1\bigr)
        \bigl(\ee^{\ii \theta_{\ell_2 , m}} - 1\bigr) \\
   &\phantom{+\sum_{1\leq\ell_1<\ell_2\leq m}
               \sum_{k_1=\lfloor nt_{\ell_1-1}\rfloor+1}^{\lfloor nt_{\ell_1}\rfloor}
                \sum_{k_2=\lfloor nt_{\ell_2-1}\rfloor+1}^{\lfloor nt_{\ell_2}\rfloor}}
        \times
        \ee^{\ii
             \bigl((\lfloor nt_{\ell_1}\rfloor-k_1)\theta_{\ell_1,m}
                   +\sum_{\ell=\ell_1+1}^{\ell_2-1}
                     \theta_{\ell,m} (\lfloor nt_\ell\rfloor - \lfloor nt_{\ell-1}\rfloor)
                   +(k_2 -1 - \lfloor nt_{\ell_2-1}\rfloor)\theta_{\ell_2, m}\bigr)} \\
   &+\sum_{\ell=1}^m
      \sum_{\lfloor nt_{\ell-1}\rfloor+1\leq k_1<k_2\leq\lfloor nt_\ell\rfloor}
       a^{k_2-k_1}
       \bigl(\ee^{\ii \theta_{\ell,m}} - 1\bigr)^2
        \ee^{\ii (k_2-k_1-1)\theta_{\ell,m}}
 \end{align*}
 with the notation \ $\theta_{j,m}:=\theta_j+ \cdots + \theta_m$, \ $j=1,\ldots,m$.

Further, we also have
 \begin{align}\label{egyuttes_char_aggregalt2}
  \EE\left(\exp\left\{ \ii \sum_{\ell=1}^m \theta_\ell S^{(1,n)}_{t_\ell} \right\}\right)
    = \int_0^1 \ee^{\lambda \widetilde K_n(a)} \psi(a) (1-a)^\beta\,\dd a,
 \end{align}
 where for all \ $a\in[0,1]$,
 \begin{align*} \widetilde K_n(a)
      &:= \sum_{\ell=-\infty}^{\lfloor nt_m \rfloor}
         \Bigg[(1-a)\sum_{k=\ell\vee 1}^{\lfloor nt_m \rfloor-1}
               a^{k-\ell} (\ee^{\ii \widetilde \theta_{\ell\vee 1,k}} - 1 - \ii \widetilde \theta_{\ell\vee 1,k})
               + a^{\lfloor nt_m \rfloor -\ell} (\ee^{\ii \widetilde \theta_{\ell\vee 1,\lfloor nt_m \rfloor}} - 1 - \ii \widetilde \theta_{\ell\vee 1,\lfloor nt_m \rfloor})
         \Bigg]\\
       & \; =  \sum_{k=1}^{\lfloor nt_m \rfloor-1} a^k (\ee^{\ii \widetilde \theta_{1,k}} - 1 - \ii \widetilde \theta_{1,k})
          + \frac{a^{\lfloor nt_m \rfloor}}{1-a}
            (\ee^{\ii \widetilde \theta_{1,\lfloor nt_m \rfloor}} - 1 - \ii \widetilde \theta_{1,\lfloor nt_m \rfloor})
\\
      &\phantom{\;=}+\sum_{\ell=1}^{\lfloor nt_m \rfloor}
         \Bigg[(1-a)\sum_{k=\ell}^{\lfloor nt_m \rfloor-1}
               a^{k-\ell} (\ee^{\ii \widetilde \theta_{\ell,k}} - 1 - \ii \widetilde \theta_{\ell,k})
               + a^{\lfloor nt_m \rfloor -\ell} (\ee^{\ii \widetilde \theta_{\ell,\lfloor nt_m \rfloor}} - 1 - \ii \widetilde \theta_{\ell,\lfloor nt_m \rfloor})
         \Bigg]
 \end{align*}
 with the notation \ $\widetilde\theta_{\ell,k}:=\widetilde\theta_\ell+ \cdots + \widetilde\theta_k$, \ $1\leq \ell\leq k \leq \lfloor nt_m\rfloor$,
 \ where \ $\widetilde\theta_j:= \sum_{i=1}^m \theta_i \bone_{\{ j\leq \lfloor nt_i \rfloor\}}$, \ $j=1,\ldots,\lfloor nt_m\rfloor$,
 \  and we define \ $\sum_{j={\lfloor nt_m \rfloor}}^{{\lfloor nt_m \rfloor}-1}:=0$.
\end{Pro}

Formulae \eqref{egyuttes_char_aggregalt1} and \eqref{egyuttes_char_aggregalt2} in Proposition \ref{Pro_egyuttes_char_aggregalt}
 have quite a different structure, and it seems to be difficult to check their equality
 not using any ingredients of the proof of Proposition \ref{Pro_egyuttes_char_aggregalt}.
However, in Section \ref{Proofs}, we present such a proof in case of \ $m=1$.
\ In the proofs we will use only \eqref{egyuttes_char_aggregalt1}, but we present \eqref{egyuttes_char_aggregalt2} as well, since it
 is interesting on its own right and it can be useful later on for handling the case \ $\beta\in(0,\infty)$ \ as well.
We note that formula \eqref{egyuttes_char_aggregalt2} is based on an infinite series representation of strictly stationary INAR(1) processes,
 recalled in Appendix \ref{App_Stac_series}.

The proofs of Theorems \ref{simultaneous_aggregation_random_2} and \ref{simultaneous_aggregation_random_0} are based on
 the explicit formula of the characteristic function of \ $(S^{(N_n,n)}_{t_1},\ldots,S^{(N_n,n)}_{t_m})$ \
 given in \eqref{egyuttes_char_aggregalt1}, where \ $0<t_1<t_2<\ldots < t_m$, \ $m\in\NN$, \ and an auxiliary Lemma \ref{ordo},
 which gives a set of sufficient conditions for the convergence of the integral
 \ $N_n\int_0^1 \left( 1- \ee^{\frac{\lambda}{1-a} z_n(a)}\right) \psi(a) (1-a)^\beta\,\dd a$ \ as \ $n\to\infty$, \
 where \ $(z_n(a))_{n\in\NN}$ \ is a sequence of complex numbers.
This proof technique is not suitable for handling other possible cases, e.g., the case \ $\beta\in(0,\infty)$, \ these remain for future work
 (for more details, see Remark \ref{Rem_discussion2}).

In the next remark we compare our assumptions in Theorems \ref{simultaneous_aggregation_random_2} and \ref{simultaneous_aggregation_random_0}
 with the corresponding assumptions in Pilipauskait{\.e} et al.\ \cite{PilSkoSur} for analogous results about simultaneous aggregation of
 random coefficient autoregressive processes of order 1.

\begin{Rem}\label{Rem_discussion1}
In Theorem \ref{simultaneous_aggregation_random_2} (where \ $\beta\in(-1,0)$), \ the condition \ $N_n^{\frac{-\beta}{1+\beta}} n^{-1} \to \infty$ \ as
 \ $n\to\infty$ \ yields that \ $N_n\to\infty$ \ as \ $n\to\infty$ \ and
 \[
   N_n^{\frac{1}{1+\beta}} n^{-1} = N_n N_n^{\frac{-\beta}{1+\beta}} n^{-1}\to\infty
   \qquad \text{as \ $n\to\infty$,}
 \]
 which is the form of the condition in Pilipauskait{\.e} et al.\ \cite{PilSkoSur} for their convergence (2.22)
 for simultaneous aggregation of random coefficient autoregressive processes of order 1 with the same mixing distribution given in \eqref{alpha}.
However, in case of \ $\beta\in(-1,0)$, \ the condition \ $N_n^{\frac{1}{1+\beta}} n^{-1}\to\infty$ \ as \ $n\to\infty$ \ does not imply that
 \ $N_n^{\frac{-\beta}{1+\beta}} n^{-1} \to \infty$ \ as \ $n\to\infty$ \ in general.
Indeed, for example, if \ $N_n:=\lfloor n^\gamma \ln n\rfloor$ \ with some \ $\gamma\in(1+\beta, -1- \frac{1}{\beta})$, \
 then \ $N_n^{\frac{1}{1+\beta}} n^{-1} \sim  n^{-1+\frac{\gamma}{1+\beta}} (\ln n)^{\frac{1}{1+\beta}}\to\infty$ \
 as \ $n\to\infty$, \ since \ $-1+\frac{\gamma}{1+\beta}>0$, \
 but \ $N_n^{\frac{-\beta}{1+\beta}} n^{-1} \sim n^{\frac{-1-\beta-\gamma\beta}{1+\beta}} (\ln n)^{\frac{-\beta}{1+\beta}}\to 0$ \
 as \ $n\to\infty$, \ since \ $\frac{-1-\beta-\gamma\beta}{1+\beta}<0$.
\ We note that the condition \ $N_n^{\frac{-\beta}{1+\beta}} n^{-1} \to \infty$ \ as \ $n\to\infty$ \ in Theorem \ref{simultaneous_aggregation_random_2}
 might be replaced by \ $N_n^{\frac{1}{1+\beta}} n^{-1}$ \ as \ $n\to\infty$.
\ However, a new proof technique would be needed, since our present one uses effectively that \ $N_n^{\frac{-\beta}{1+\beta}} n^{-1} \to \infty$ \ as \ $n\to\infty$,
 \ for example, in the proof of Theorem \ref{simultaneous_aggregation_random_2} we argued that for large enough \ $n$ \ and for any
 \ $z\in (N_n^{-1},1]$, \ we have \ $z^{-1} n N_n^{\frac{-1}{1+\beta}} \vert O(1)\vert \leq \vert O(1)\vert$ \ (see \eqref{help_cond_needed1}).

In Theorem \ref{simultaneous_aggregation_random_0} (where \ $\beta = 0$), \ the condition \ $(\log N_n)^2 n^{-1} \to \infty$ \ as \ $n\to\infty$ \ yields that
 \ $N_n\to\infty$ \ as \ $n\to\infty$ \ and \ $N_nn^{-1} = n^{-1} (\log N_n)^2\frac{N_n}{(\log N_n)^2}\to\infty\cdot \infty=\infty$ \ as \ $n\to\infty$, \
 which is the form of the condition in Pilipauskait{\.e} et al. \cite{PilSkoSur} for their convergence (2.21).
However, the condition \ $N_nn^{-1}\to\infty$ \ as \ $n\to\infty$ \ does not imply that
  \ $(\log N_n)^2 n^{-1} \to \infty$ \ as \ $n\to\infty$ \ in general.
Indeed, for example, if \ $N_n:=n^2$, \ then \ $N_n n^{-1} = n\to\infty$ \ as \ $n\to\infty$, \ but
 \ $(\log N_n)^2 n^{-1} = 4 n^{-1} (\ln n)^2\to0$ \ as \ $n\to\infty$.
\ Further, one can check that
 \[
 \lim_{n\to\infty}
     \frac{n^{-1} (N_n \log N_n )^{-\frac{1}{2}}}{ n^{-1} (N_n \log(N_n/n))^{-\frac{1}{2}} } = 1,
 \]
 where \ $n^{-1} (N_n \log(N_n/n))^{-\frac{1}{2}}$ \ is the scaling factor in (2.21) in Pilipauskait{\.e} et al. \cite{PilSkoSur}.
Indeed,
 \begin{align*}
  \lim_{n\to\infty}
     \frac{n^{-1} (N_n \log N_n)^{-\frac{1}{2}}}{ n^{-1} (N_n \log(N_n/n))^{-\frac{1}{2}} }
   = \lim_{n\to\infty}
     \left( \frac{N_n\log N_n}{N_n \log N_n  - N_n \log n} \right)^{-\frac{1}{2}}
   = \lim_{n\to\infty}
     \left( \frac{1}{1 - \log n/\log N_n} \right)^{-\frac{1}{2}}
   =1,
 \end{align*}
 since \ $(\log n /\log N_n )^2 = [ (\log N_n )^{-2} n] n^{-1}(\log n )^2\to 0\cdot 0 = 0$ \ as  \ $n\to\infty$
 \ under the condition \ $(\log N_n)^2 n^{-1} \to \infty$ \ as \ $n\to\infty$.
We note that the condition \ $(\log N_n)^2 n^{-1} \to \infty$ \ as \ $n\to\infty$ \ in Theorem \ref{simultaneous_aggregation_random_0}
 might be replaced by \ $N_nn^{-1}$ \ as \ $n\to\infty$.
\ However, a new proof technique would be needed, since our present one uses effectively that \ $(\log N_n)^2 n^{-1} \to \infty$ \ as \ $n\to\infty$,
 \ for example, in the proof of Theorem \ref{simultaneous_aggregation_random_0} we argued that
 \ $\frac{n}{N_n \log N_n} \int_{(\log N_n)^{-1}}^{N_n(\log N_n)^{-1}} 1 \, \dd z = \frac{n}{(\log N_n)^2} \left(1 - \frac{1}{N_n}\right) \to 0$
 \ as \ $n\to\infty$ \ (see \eqref{help_int2}).
\proofend
\end{Rem}

In the next remark we shed some light on why our proof technique is not suitable for other cases, e.g., the case \ $\beta\in(0,\infty)$ \
 being of the greatest interest, which remain for future work.

\begin{Rem}\label{Rem_discussion2}
Motivated by the fact that Theorem \ref{simultaneous_aggregation_random_2} can be considered as a counterpart of
 Theorem 4.8 in Barczy et al.\ \cite{BarNedPap}, using Theorem 4.10 in Barczy et al.\ \cite{BarNedPap} (which is about the iterated aggregation case
 first taking the limit \ $n\to\infty$ \ and then \ $N\to\infty$ \ in case of \ $\beta\in(-1,1)$), \ one might suspect that
 in case of \ $\beta\in(-1,1)$ \ in the simultaneous aggregation case an appropriate limit theorem might hold with the same scaling
 \ $n^{-\frac{1}{2}} N_n^{-\frac{1}{1+\beta}}$ \ and with the same limit distribution as in Theorem 4.10 in Barczy et al.\ \cite{BarNedPap}
 when both \ $n$ \ and \ $N_n$ \ increase to infinity at some appropriate rate.
Unfortunately, our proof technique used in the present paper, i.e. an application of Lemma \ref{ordo}, is not suitable for handling this case,
 since in order to fulfill condition \eqref{help_ordo1} of Lemma \ref{ordo} one is forced to choose a sequence \ $(\vare_n)_{n\in\NN}$ \ in \ $(0,1)$ \
 satisfying \ $\vare_n \leq L nN_n^{\frac{-1+\beta}{1+\beta}}$ \ with some constant \ $L>0$ \ and \ $\vare_n\to0$ \ as \ $n\to\infty$.
\ However, with such a possible choice of \ $(\vare_n)_{n\in\NN}$ \ we were not able to check the other two conditions of Lemma \ref{ordo}.

Very recently, for randomized autoregressive processes of order 1,
 Pilipauskait{\.e} et al.\ \cite{PilSkoSur} have found a somewhat new approach for studying
 simultaneous limits.
Namely, they used an infinite series representation of the stationary distribution of their model for
 calculating the characteristic function of finite dimensional distributions in question.
In our case, i.e., in case of randomized INAR(1) processes, we also derived such a formula
 given in \eqref{egyuttes_char_aggregalt2}, and it is much more complicated.
As a future work, using it, we plan to handle the remaining case \ $\beta\in(0,\infty)$, \
 where a completely different limit behaviour is expected.
\proofend
\end{Rem}

The paper is structured as follows.
Section  \ref{Proofs} contains the proofs, first the proof of Proposition \ref{Pro_egyuttes_char_aggregalt}
 and a direct proof of the equality of the formulae \eqref{egyuttes_char_aggregalt1} and \eqref{egyuttes_char_aggregalt2}
 in case of \ $m=1$, \ then the proofs of Theorems \ref{simultaneous_aggregation_random_2} and \ref{simultaneous_aggregation_random_0}.
We close the paper with three appendices.
In Appendix \ref{prel} we recall the generator function of finite-dimensional distributions of stationary INAR(1) processes with
 Poisson immigrations.
Appendix \ref{App_Stac_series} is devoted to an infinite series representation of strictly stationary INAR(1) processes in question.
Finally, Appendix \ref{App1} contains some approximations of the exponential function and some of its integral
 extensively used in the proofs of Theorems \ref{simultaneous_aggregation_random_2} and \ref{simultaneous_aggregation_random_0}.

\section{Proofs}
\label{Proofs}

For a non-negative integer-valued random variable \ $\zeta$, \ its characteristic function and generating function will be
  denoted by \ $\varphi_\zeta$ \ and \ $G_\zeta$, \ respectively.
For a non-negative integer-valued random variable \ $L$ \ and \ $p\in[0,1]$, \ $\textrm{Bin}(L,p)$ \ will denote a random variable having a binomial distribution
 with parameters \ $L$ \ and \ $p$ \ (meaning that the conditional distribution of \ $\textrm{Bin}(L,p)$ \ given \ $L=\ell$, \ $\ell\in\ZZ_+$, \ is a binomial
 distribution with parameters \ $\ell$ \ and \ $p$).
The notations \ $\OO(1)$ \ and \ $|\!\OO(1)|$ \ stand for a possibly complex and real sequence \ $(a_k)_{k\in\NN}$, \ respectively, that is bounded and can only depend
 on the parameters \ $\lambda$, \ $\psi_1$, \ $\beta$, \ and on some fixed \ $m \in \NN$
 \ and \ $\theta_1,\dots, \theta_m \in \RR$.
\ Further, we call the attention that the multiple \ $\OO(1)$  and
 \ $|\!\OO(1)|$ \ notations in the same formula do not necessarily mean the same bounded
 sequence.

\noindent{\bf Proof of Proposition \ref{Pro_egyuttes_char_aggregalt}.}
First, we prove \eqref{egyuttes_char_aggregalt1}.
Since
 \begin{align*}
   \sum_{\ell=1}^m \theta_\ell S^{(1,n)}_{t_\ell}
     =  \begin{bmatrix}
          \theta_1 \\
          \vdots \\
          \theta_m \\
        \end{bmatrix}^\top
        \begin{bmatrix}
          S^{(1,n)}_{t_1} \\
          \vdots \\
          S^{(1,n)}_{t_m} \\
        \end{bmatrix}
     =  \begin{bmatrix}
          \theta_1+\cdots+\theta_m \\
          \theta_2+\cdots+\theta_m \\
          \vdots \\
          \theta_{m-1}+\theta_m\\
          \theta_m \\
        \end{bmatrix}^\top
        \begin{bmatrix}
          S^{(1,n)}_{t_1} \\[1mm]
          S^{(1,n)}_{t_2} - S^{(1,n)}_{t_1}\\
          \vdots \\
          S^{(1,n)}_{t_{m-1}} - S^{(1,n)}_{t_{m-2}}\\[1mm]
          S^{(1,n)}_{t_m} - S^{(1,n)}_{t_{m-1}} \\
        \end{bmatrix},
 \end{align*}
 by the law of total expectation, we have
 \begin{align*}
   &\EE\left(\exp\left\{ \ii \sum_{\ell=1}^m \theta_\ell S^{(1,n)}_{t_\ell} \right\}\right)
     = \EE\left(\exp\left\{ \ii \sum_{\ell=1}^m \theta_{\ell,m} ( S^{(1,n)}_{t_\ell} - S^{(1,n)}_{t_{\ell-1}} ) \right\}\right)\\
   & = \int_0^1 \EE\left(\exp\left\{ \ii \sum_{\ell=1}^m \theta_{\ell,m} \sum_{k=\lfloor nt_{\ell-1}\rfloor+1}^{\lfloor nt_\ell\rfloor}
        \left(X_k - \frac{\lambda}{1-\alpha}\right) \right\} \Bigg\vert\; \alpha=a\right)
       \psi(a) (1-a)^\beta\,\dd a,
 \end{align*}
 where for all \ $a\in[0,1]$ \ we have
 \begin{align*}
  &\EE\left(\exp\left\{ \ii \sum_{\ell=1}^m \theta_{\ell,m} \sum_{k=\lfloor nt_{\ell-1}\rfloor+1}^{\lfloor nt_\ell\rfloor}
        \left(X_k - \frac{\lambda}{1-\alpha}\right) \right\} \Bigg\vert\; \alpha=a\right)\\
  & = \EE\Bigl[\ee^{-\ii \frac{\lambda}{1-a}
            \sum_{\ell=1}^m \theta_{\ell,m} (\lfloor{nt_{\ell}}\rfloor - \lfloor nt_{\ell-1}\rfloor)} \\
  & \phantom{= \EE\Bigl[}
     \times F_{0,\dots,\lfloor nt_m\rfloor-1}
      \Bigl( \underbrace{ \ee^{\ii \theta_{1,m}},\ldots, \ee^{\ii \theta_{1,m}} }_{\lfloor nt_1\rfloor \ \text{items}},
             \underbrace{ \ee^{\ii \theta_{2,m}},\ldots, \ee^{\ii \theta_{2,m}} }_{\lfloor nt_2\rfloor - \lfloor nt_1\rfloor \ \text{items}},
             \dots,
             \underbrace{ \ee^{\ii \theta_{m,m}}, \ldots, \ee^{\ii \theta_{m,m}} }_{\lfloor nt_m\rfloor - \lfloor nt_{m-1}\rfloor \ \text{items}}  \,\Big|\, \alpha=a\Bigr) \Bigr]
 \end{align*}
 where recall that \ $F_{0,\dots,\lfloor nt_m\rfloor-1}(z_0,\ldots,z_{\lfloor nt_m\rfloor-1}\mid \alpha=a)$ \ denotes the conditional
 joint generating function of \ $(X_0,X_1,\ldots, X_{\lfloor nt_m\rfloor-1})$ \ given \ $\alpha=a$ \ at \ $(z_0,\ldots,z_{\lfloor nt_m\rfloor-1})
  \in\CC^{\lfloor nt_m\rfloor}$.
\ Then an application of \eqref{help1_alter} yields \eqref{egyuttes_char_aggregalt1}.

Now we turn to prove \eqref{egyuttes_char_aggregalt2}.
Using again the law of total expectation we have
 \begin{align*}
   &\EE\left(\exp\left\{ \ii \sum_{i=1}^m \theta_i S^{(1,n)}_{t_i} \right\}\right)\\
   & = \int_0^1 \EE\left(\exp\left\{ \ii \sum_{i=1}^m \theta_i \sum_{k=1}^{\lfloor nt_i\rfloor} \left(X_k - \frac{\lambda}{1-\alpha}\right) \right\} \Bigg\vert\; \alpha=a\right)
       \psi(a) (1-a)^\beta\,\dd a,
 \end{align*}
 where, by \eqref{repr} and the fact that \ $\EE(X_k \mid \alpha) = \frac{\lambda}{1-\alpha} = \lambda \sum_{\ell=0}^\infty \alpha^\ell$, \ $k\in\NN$,
 \ we get Lebesgue a.e. \ $a\in[0,1]$,
 \begin{align*}
   &\EE\left(\exp\left\{ \ii \sum_{i=1}^m \theta_i \sum_{k=1}^{\lfloor nt_i\rfloor} \left(X_k - \frac{\lambda}{1-\alpha}\right) \right\} \Bigg\vert\; \alpha=a\right) \\
   & = \EE\left(\exp\left\{ \ii \sum_{i=1}^m \theta_i \sum_{k=1}^{\lfloor nt_i\rfloor} \sum_{\ell=0}^\infty
          \Big( a_k^{(k-\ell)} \circ \cdots \circ a_{k-\ell+1}^{(k-\ell)} \circ \vare_{k-\ell} - \lambda a^\ell \Big) \right\} \right)\\
   & = \EE\left(\exp\left\{ \ii \sum_{i=1}^m \theta_i \sum_{k=1}^{\lfloor nt_i\rfloor} \sum_{\ell=-\infty}^k
          \Big( a_k^{(\ell)} \circ \cdots \circ a_{\ell+1}^{(\ell)} \circ \vare_\ell - \lambda a^{k-\ell} \Big) \right\} \right)
 \end{align*}
 \begin{align*}
   & = \EE\Bigg( \exp\Bigg\{ \ii \sum_{i=1}^m  \theta_i \sum_{\ell=-\infty}^0 \sum_{k=1}^{\lfloor nt_i\rfloor}
            \Big( a_k^{(\ell)} \circ \cdots \circ a_{\ell+1}^{(\ell)} \circ \vare_\ell - \lambda a^{k-\ell} \Big) \\
   &\phantom{= \EE\Bigg( \exp\Bigg\{ }
           + \ii \sum_{i=1}^m  \theta_i \sum_{\ell=1}^{\lfloor nt_i\rfloor}  \sum_{k=\ell}^{\lfloor nt_i\rfloor}
            \Big( a_k^{(\ell)} \circ \cdots \circ a_{\ell+1}^{(\ell)} \circ \vare_\ell - \lambda a^{k-\ell} \Big) \Bigg\}
           \Bigg),
 \end{align*}
 where for \ $\ell=0$ \ and \ $k\in\NN$, \ $a_k^{(k-\ell)} \circ \cdots \circ a_{k-\ell+1}^{(k-\ell)} \circ \vare_{k-\ell}$ \ is defined to be \ $\vare_k$.
Here
 \begin{align*}
   &\sum_{i=1}^m  \theta_i \sum_{\ell=-\infty}^0 \sum_{k=1}^{\lfloor nt_i\rfloor}
            \Big( a_k^{(\ell)} \circ \cdots \circ a_{\ell+1}^{(\ell)} \circ \vare_\ell - \lambda a^{k-\ell} \Big) \\
   &\qquad = \sum_{\ell=-\infty}^0 \sum_{k=1}^{\lfloor nt_m\rfloor} \left( \sum_{i=1}^m \theta_i \bone_{\{ k\leq \lfloor nt_i\rfloor \}}\right)
       \Big( a_k^{(\ell)} \circ \cdots \circ a_{\ell+1}^{(\ell)} \circ \vare_\ell - \lambda a^{k-\ell} \Big),
 \end{align*}
  and
 \begin{align*}
  &\sum_{i=1}^m  \theta_i \sum_{\ell=1}^{\lfloor nt_i\rfloor}  \sum_{k=\ell}^{\lfloor nt_i\rfloor}
            \Big( a_k^{(\ell)} \circ \cdots \circ a_{\ell+1}^{(\ell)} \circ \vare_\ell - \lambda a^{k-\ell} \Big) \\
  &\qquad = \sum_{\ell=1}^{\lfloor nt_m\rfloor} \sum_{k=\ell}^{\lfloor nt_m\rfloor} \left( \sum_{i=1}^m \theta_i \bone_{\{ k\leq \lfloor nt_i\rfloor \}} \right)
      \Big( a_k^{(\ell)} \circ \cdots \circ a_{\ell+1}^{(\ell)} \circ \vare_\ell - \lambda a^{k-\ell} \Big).
 \end{align*}
Consequently, for Lebesgue a.e. \ $a\in[0,1]$,
 \begin{align*}
  &\EE\left(\exp\left\{ \ii \sum_{i=1}^m \theta_i \sum_{k=1}^{\lfloor nt_i\rfloor} \left(X_k - \frac{\lambda}{1-\alpha}\right) \right\} \Bigg\vert\; \alpha=a\right) \\
  &\qquad  = \EE\Bigg( \exp\Bigg\{ \ii \sum_{\ell=-\infty}^{\lfloor nt_m\rfloor} \sum_{k=\ell\vee 1}^{\lfloor nt_m\rfloor}
            \widetilde \theta_k \Big( a_k^{(\ell)} \circ \cdots \circ a_{\ell+1}^{(\ell)} \circ \vare_\ell - \lambda a^{k-\ell} \Big) \Bigg\} \Bigg).
 \end{align*}
Next we show that for all \ $t>0$, \ $n\in\NN$, \ and \ $\ell\leq \lfloor nt\rfloor$, \ $\ell\in\ZZ$, \ we have
 \begin{align}\label{help_joint_characteristic2}
  \begin{split}
 &\EE\Bigg( \exp\Bigg\{ \ii  \sum_{k=\ell\vee 1}^{\lfloor nt\rfloor}
            \widetilde \theta_k \Big( a_k^{(\ell)} \circ \cdots \circ a_{\ell+1}^{(\ell)} \circ \vare_\ell - \lambda a^{k-\ell} \Big) \Bigg\} \Bigg) \\
 &= \exp\Bigg\{\lambda
   \Bigg[(1-a)\sum_{k=\ell\vee 1}^{\lfloor nt \rfloor-1}
               a^{k-\ell} (\ee^{\ii \widetilde \theta_{\ell\vee 1,k}} - 1 - \ii \widetilde \theta_{\ell\vee 1,k})
               + a^{\lfloor nt \rfloor -\ell} (\ee^{\ii \widetilde \theta_{\ell\vee 1,\lfloor nt \rfloor}} - 1 - \ii \widetilde \theta_{\ell\vee 1,\lfloor nt \rfloor})
         \Bigg]\Bigg\}
  \end{split}
 \end{align}
 for all \ $a\in[0,1]$, \ which together with  the independence of \ $\{\vare_k : k\in\ZZ\}$ \ and \ $a^{(\ell)}_k$, \ $\ell,k\in\ZZ$, \ yield \eqref{egyuttes_char_aggregalt2}.
First we prove \eqref{help_joint_characteristic2} in case of \ $1\leq \ell\leq \lfloor nt\rfloor$, \ $\ell\in\ZZ$, \ yielding that \ $\ell\vee 1=\ell$.
\ By the tower rule we have for any \ $a\in[0,1]$,
 \begin{align*}
  &\EE\Bigg( \exp\Bigg\{ \ii  \sum_{k=\ell\vee 1}^{\lfloor nt\rfloor}
            \widetilde \theta_k \Big( a_k^{(\ell)} \circ \cdots \circ a_{\ell+1}^{(\ell)} \circ \vare_\ell - \lambda a^{k-\ell} \Big) \Bigg\} \Bigg)
=
\EE\Bigg( \exp\Bigg\{ \ii  \sum_{k=\ell}^{\lfloor nt\rfloor}
            \widetilde \theta_k \Big( a_k^{(\ell)} \circ \cdots \circ a_{\ell+1}^{(\ell)} \circ \vare_\ell - \lambda a^{k-\ell} \Big) \Bigg\} \Bigg)\\
  & =\EE\Bigg(\exp\Bigg\{\ii\Big(\widetilde\theta_\ell (\vare_\ell - \lambda) + \widetilde\theta_{\ell+1} (a^{(\ell)}_{\ell+1}\circ\vare_\ell - \lambda a) + \cdots
     +\widetilde\theta_{\lfloor nt\rfloor -1}( a_{\lfloor nt\rfloor-1}^{(\ell)} \circ \cdots \circ a_{\ell+1}^{(\ell)} \circ \vare_\ell - \lambda a^{\lfloor nt\rfloor -1-\ell} )
        \Big) \Bigg\}\\
  &\quad \times\EE\Bigg( \exp\Bigg\{\ii
      \widetilde\theta_{\lfloor nt\rfloor}( a_{\lfloor nt\rfloor}^{(\ell)} \circ \cdots \circ a_{\ell+1}^{(\ell)} \circ \vare_\ell - \lambda a^{\lfloor nt\rfloor-\ell} )
     \Bigg\} \Bigg\vert \; \vare_\ell, a_{\ell+1}^{(\ell)}\circ\vare_\ell, \ldots, a_{\lfloor nt\rfloor-1}^{(\ell)} \circ \cdots \circ a_{\ell+1}^{(\ell)} \circ \vare_\ell \Bigg)
      \Bigg)\\
  &= \EE\Bigg(\exp\Bigg\{\ii\Big(\widetilde\theta_\ell (\vare_\ell - \lambda) + \widetilde\theta_{\ell+1} (a^{(\ell)}_{\ell+1}\circ\vare_\ell - \lambda a) +\cdots
     +\widetilde\theta_{\lfloor nt\rfloor -1}( a_{\lfloor nt\rfloor-1}^{(\ell)} \circ \cdots \circ a_{\ell+1}^{(\ell)} \circ \vare_\ell - \lambda a^{\lfloor nt\rfloor -1-\ell} )
        \Big) \Bigg\}\\
  & \phantom{= \EE\Bigg(}
      \times \ee^{-\ii \widetilde\theta_{\lfloor nt\rfloor} \lambda a^{\lfloor nt\rfloor-\ell}}
            \varphi_{\textrm{Bin}\big(a_{\lfloor nt\rfloor-1}^{(\ell)} \circ \cdots \circ a_{\ell+1}^{(\ell)} \circ \vare_\ell,a\big)}(\widetilde\theta_{\lfloor nt\rfloor})
   \Bigg)\\
  &=  \ee^{-\ii \lambda \widetilde\theta_{\lfloor nt\rfloor}  a^{\lfloor nt\rfloor-\ell}}\\
  &\phantom{=}\times\EE\Bigg(\exp\Bigg\{\ii\Big(\widetilde\theta_\ell (\vare_\ell - \lambda) + \widetilde\theta_{\ell+1} (a^{(\ell)}_{\ell+1}\circ\vare_\ell - \lambda a) +\cdots
     +\widetilde\theta_{\lfloor nt\rfloor -1}( a_{\lfloor nt\rfloor-1}^{(\ell)} \circ \cdots \circ a_{\ell+1}^{(\ell)} \circ \vare_\ell - \lambda a^{\lfloor nt\rfloor -1-\ell} )
        \Big) \Bigg\}\\
  & \phantom{= \times\EE\Bigg(}
      \times (1-a + a\ee^{\ii\widetilde\theta_{\lfloor nt\rfloor}})^{a_{\lfloor nt\rfloor-1}^{(\ell)} \circ \cdots \circ a_{\ell+1}^{(\ell)} \circ \vare_\ell}
      \Bigg)\\
  &=  \ee^{-\ii \lambda \big( \widetilde\theta_{\lfloor nt\rfloor - 1}  a^{\lfloor nt\rfloor-1-\ell} + \widetilde\theta_{\lfloor nt\rfloor}  a^{\lfloor nt\rfloor-\ell}\big) }\\
  &\phantom{=}\times\EE\Bigg(\exp\Bigg\{\ii\Big(\widetilde\theta_\ell (\vare_\ell - \lambda) + \widetilde\theta_{\ell+1} (a^{(\ell)}_{\ell+1}\circ\vare_\ell - \lambda a) +\cdots
     +\widetilde\theta_{\lfloor nt\rfloor -2}( a_{\lfloor nt\rfloor-2}^{(\ell)} \circ \cdots \circ a_{\ell+1}^{(\ell)} \circ \vare_\ell - \lambda a^{\lfloor nt\rfloor -2-\ell} )
        \Big) \Bigg\}\\
  & \phantom{= \times\EE\Bigg(}
        \times \big( (1-a)\ee^{\ii\widetilde\theta_{\lfloor nt\rfloor-1}}  + a \ee^{\ii(\widetilde\theta_{\lfloor nt\rfloor-1} + \widetilde\theta_{\lfloor nt\rfloor})}\big)^{a_{\lfloor nt\rfloor-1}^{(\ell)} \circ \cdots \circ a_{\ell+1}^{(\ell)} \circ \vare_\ell}
      \Bigg)=
  \end{align*}
 \begin{align*}
  &=  \ee^{-\ii \lambda \big( \widetilde\theta_{\lfloor nt\rfloor - 1}  a^{\lfloor nt\rfloor-1-\ell} + \widetilde\theta_{\lfloor nt\rfloor}  a^{\lfloor nt\rfloor-\ell}\big) }\\
  &\phantom{=}\times\EE\Bigg(\exp\Bigg\{\ii\Big(\widetilde\theta_\ell (\vare_\ell - \lambda) + \widetilde\theta_{\ell+1} (a^{(\ell)}_{\ell+1}\circ\vare_\ell - \lambda a) +\cdots
     +\widetilde\theta_{\lfloor nt\rfloor -2}( a_{\lfloor nt\rfloor-2}^{(\ell)} \circ \cdots \circ a_{\ell+1}^{(\ell)} \circ \vare_\ell - \lambda a^{\lfloor nt\rfloor -2-\ell} )
        \Big) \Bigg\}\\
  & \phantom{= \times\EE\Bigg(}
      \times G_{\textrm{Bin}\big(a_{\lfloor nt\rfloor-2}^{(\ell)} \circ \cdots \circ a_{\ell+1}^{(\ell)} \circ \vare_\ell,a\big)}
         \big( (1-a)\ee^{\ii\widetilde\theta_{\lfloor nt\rfloor-1}}  + a \ee^{\ii(\widetilde\theta_{\lfloor nt\rfloor-1} + \widetilde\theta_{\lfloor nt\rfloor})}\big)
      \Bigg)\\
  &=  \ee^{-\ii \lambda \big( \widetilde\theta_{\lfloor nt\rfloor - 1}  a^{\lfloor nt\rfloor-1-\ell} + \widetilde\theta_{\lfloor nt\rfloor}  a^{\lfloor nt\rfloor-\ell}\big) }\\
  &\phantom{=}\times\EE\Bigg(\exp\Bigg\{\ii\Big(\widetilde\theta_\ell (\vare_\ell - \lambda) + \widetilde\theta_{\ell+1} (a^{(\ell)}_{\ell+1}\circ\vare_\ell - \lambda a) +\cdots
     +\widetilde\theta_{\lfloor nt\rfloor -2}( a_{\lfloor nt\rfloor-2}^{(\ell)} \circ \cdots \circ a_{\ell+1}^{(\ell)} \circ \vare_\ell - \lambda a^{\lfloor nt\rfloor -2-\ell} )
        \Big) \Bigg\}\\
  & \phantom{= \times\EE\Bigg(}
      \times
      \bigg( 1-a + a(1-a)\ee^{\ii\widetilde\theta_{\lfloor nt\rfloor-1}}
         + a^2 \ee^{\ii(\widetilde\theta_{\lfloor nt\rfloor-1} + \widetilde\theta_{\lfloor nt\rfloor})}\bigg)^{a_{\lfloor nt\rfloor-2}^{(\ell)} \circ \cdots \circ a_{\ell+1}^{(\ell)} \circ \vare_\ell}
      \Bigg)\\
  & = \cdots =\\
  &= \ee^{-\ii\lambda \sum_{k=\ell}^{\lfloor nt\rfloor} \widetilde\theta_k a^{k-\ell}}
          G_{\vare_\ell}\Big( (1-a) \sum_{k=\ell}^{\lfloor nt\rfloor-1}  a^{k-\ell} \ee^{\ii\widetilde\theta_{\ell,k}}
                       + a^{\lfloor nt\rfloor-\ell}\ee^{\ii\widetilde\theta_{\ell,\lfloor nt\rfloor}} \Big)\\
  &= \ee^{-\ii\lambda \sum_{k=\ell}^{\lfloor nt\rfloor} \widetilde\theta_k a^{k-\ell}}
     \exp\Big\{ -\lambda + \lambda \Big( (1-a) \sum_{k=\ell}^{\lfloor nt\rfloor-1}  a^{k-\ell} \ee^{\ii\widetilde\theta_{\ell,k}}
                        +a^{\lfloor nt\rfloor-\ell}\ee^{\ii\widetilde\theta_{\ell,\lfloor nt\rfloor}} \Big) \Big\},
 \end{align*}
 which coincides with
 \begin{align}\label{help_joint_characteristic1}
    \exp\Bigg\{\lambda
   \Bigg[(1-a)\sum_{k=\ell}^{\lfloor nt \rfloor-1}
               a^{k-\ell} (\ee^{\ii \widetilde \theta_{\ell,k}} - 1 - \ii \widetilde \theta_{\ell,k})
               + a^{\lfloor nt \rfloor -\ell} (\ee^{\ii \widetilde \theta_{\ell,\lfloor nt \rfloor}} - 1 - \ii \widetilde \theta_{\ell,\lfloor nt \rfloor})
         \Bigg]\Bigg\},
 \end{align}
 as desired.
Indeed, the coefficient of the constant term (not depending on \ $\widetilde \theta_j$, \ $j=\ell,\ldots,\lfloor nt \rfloor$) \ in the exponential
 of \eqref{help_joint_characteristic1} is
 \begin{align*}
   -\lambda(1-a)\sum_{k=\ell}^{\lfloor nt \rfloor - 1} a^{k-\ell} - \lambda a^{\lfloor nt \rfloor - \ell}
      = -\lambda(1-a)\frac{a^{\lfloor nt \rfloor - \ell} - 1}{a-1} - \lambda a^{\lfloor nt \rfloor - \ell}
      = -\lambda,
 \end{align*}
 the coefficient of \ $\widetilde \theta_j$, \ $j\in\{\ell,\ldots,\lfloor nt \rfloor-1\}$, \ in the exponential of \eqref{help_joint_characteristic1} is
 \[
  -\ii \lambda (1-a) \sum_{k=j}^{\lfloor nt \rfloor - 1} a^{k-\ell} - \ii \lambda a^{\lfloor nt \rfloor-\ell}
    = -\ii\lambda (1-a) a^{j-\ell} \frac{a^{\lfloor nt \rfloor - j} - 1}{a-1} - \ii \lambda a^{\lfloor nt \rfloor-\ell}
    = -\ii \lambda a^{j-\ell},
 \]
the coefficient of \ $\widetilde \theta_{\lfloor nt \rfloor}$ \ in the exponential of \eqref{help_joint_characteristic1} is
 \ $-\ii\lambda a^{\lfloor nt \rfloor - \ell}$,
 \ the coefficient of \ $\ee^{\ii\widetilde \theta_{\ell,\lfloor nt \rfloor}}$ \ in the exponential of \eqref{help_joint_characteristic1} is
 \ $\lambda a^{\lfloor nt \rfloor - \ell}$, \ and the remaining term \ $\lambda(1-a)\sum_{k=\ell}^{\lfloor nt \rfloor-1}
 a^{k-\ell} \ee^{\ii \widetilde \theta_{\ell,k}}$ \ coincide as well.

Finally, we prove \eqref{help_joint_characteristic2} in case of \ $\ell\leq 0$, \ $\ell\in\ZZ$, \ yielding that \ $\ell\vee 1 = 1$.
\ By the tower rule we have for any \ $a\in[0,1]$,
 \begin{align*}
  &\EE\Bigg( \exp\Bigg\{ \ii  \sum_{k=\ell\vee 1}^{\lfloor nt\rfloor}
            \widetilde \theta_k \Big( a_k^{(\ell)} \circ \cdots \circ a_{\ell+1}^{(\ell)} \circ \vare_\ell - \lambda a^{k-\ell} \Big) \Bigg\} \Bigg)
    = \EE\Bigg( \exp\Bigg\{ \ii  \sum_{k=1}^{\lfloor nt\rfloor}
            \widetilde \theta_k \Big( a_k^{(\ell)} \circ \cdots \circ a_{\ell+1}^{(\ell)} \circ \vare_\ell - \lambda a^{k-\ell} \Big) \Bigg\} \Bigg)  \\
  & =\EE\Bigg(\exp\Bigg\{\ii\Big(\widetilde\theta_1 (a_1^{(\ell)}\circ \cdots \circ a^{(\ell)}_{\ell+1}\circ\vare_\ell - \lambda a^{1-\ell}) + \cdots
     +\widetilde\theta_{\lfloor nt\rfloor -1}( a_{\lfloor nt\rfloor-1}^{(\ell)} \circ \cdots \circ a_{\ell+1}^{(\ell)} \circ \vare_\ell - \lambda a^{\lfloor nt\rfloor -1-\ell} )
        \Big) \Bigg\}\\
  & \times\EE\Bigg( \exp\Bigg\{\ii
      \widetilde\theta_{\lfloor nt\rfloor}( a_{\lfloor nt\rfloor}^{(\ell)} \circ \cdots \circ a_{\ell+1}^{(\ell)} \circ \vare_\ell - \lambda a^{\lfloor nt\rfloor-\ell} )
     \Bigg\} \Bigg\vert \; a_1^{(\ell)}\circ \cdots \circ a_{\ell+1}^{(\ell)}\circ\vare_\ell, \ldots, a_{\lfloor nt\rfloor-1}^{(\ell)} \circ \cdots \circ a_{\ell+1}^{(\ell)} \circ \vare_\ell \Bigg)
      \Bigg)\\
  &= \EE\Bigg(\exp\Bigg\{\ii\Big(\widetilde\theta_1 (a_1^{(\ell)}\circ \cdots \circ a^{(\ell)}_{\ell+1}\circ\vare_\ell - \lambda a^{1-\ell}) +\cdots
     +\widetilde\theta_{\lfloor nt\rfloor -1}( a_{\lfloor nt\rfloor-1}^{(\ell)} \circ \cdots \circ a_{\ell+1}^{(\ell)} \circ \vare_\ell - \lambda a^{\lfloor nt\rfloor -1-\ell} )
        \Big) \Bigg\}\\
  & \phantom{= \EE\Bigg(}
      \times \ee^{-\ii \widetilde\theta_{\lfloor nt\rfloor} \lambda a^{\lfloor nt\rfloor-\ell}}
            \varphi_{\textrm{Bin}\big(a_{\lfloor nt\rfloor-1}^{(\ell)} \circ \cdots \circ a_{\ell+1}^{(\ell)} \circ \vare_\ell,a\big)}(\widetilde\theta_{\lfloor nt\rfloor})
   \Bigg)\\
  &=  \ee^{-\ii \lambda \widetilde\theta_{\lfloor nt\rfloor}  a^{\lfloor nt\rfloor-\ell}}\\
  &\phantom{=}\times\EE\Bigg(\exp\Bigg\{\ii\Big( \widetilde\theta_1 (a_1^{(\ell)}\circ \cdots \circ a^{(\ell)}_{\ell+1}\circ\vare_\ell - \lambda a^{1-\ell}) +\cdots
     +\widetilde\theta_{\lfloor nt\rfloor -1}( a_{\lfloor nt\rfloor-1}^{(\ell)} \circ \cdots \circ a_{\ell+1}^{(\ell)} \circ \vare_\ell - \lambda a^{\lfloor nt\rfloor -1-\ell} )
        \Big) \Bigg\}\\
  & \phantom{= \times\EE\Bigg(}
      \times (1-a + a\ee^{\ii\widetilde\theta_{\lfloor nt\rfloor}})^{a_{\lfloor nt\rfloor-1}^{(\ell)} \circ \cdots \circ a_{\ell+1}^{(\ell)} \circ \vare_\ell}
      \Bigg)\\
  &=  \ee^{-\ii \lambda \big( \widetilde\theta_{\lfloor nt\rfloor - 1}  a^{\lfloor nt\rfloor-1-\ell} + \widetilde\theta_{\lfloor nt\rfloor}  a^{\lfloor nt\rfloor-\ell}\big) }\\
  &\phantom{=}\times\EE\Bigg(\exp\Bigg\{\ii\Big( \widetilde\theta_1 (a_1^{(\ell)}\circ \cdots \circ a^{(\ell)}_{\ell+1}\circ\vare_\ell - \lambda a^{1-\ell}) +\cdots
     +\widetilde\theta_{\lfloor nt\rfloor -2}( a_{\lfloor nt\rfloor-2}^{(\ell)} \circ \cdots \circ a_{\ell+1}^{(\ell)} \circ \vare_\ell - \lambda a^{\lfloor nt\rfloor -2-\ell} )
        \Big) \Bigg\}\\
  & \phantom{= \times\EE\Bigg(}
        \times \big( (1-a)\ee^{\ii\widetilde\theta_{\lfloor nt\rfloor-1}}  + a \ee^{\ii(\widetilde\theta_{\lfloor nt\rfloor-1} + \widetilde\theta_{\lfloor nt\rfloor})}\big)^{a_{\lfloor nt\rfloor-1}^{(\ell)} \circ \cdots \circ a_{\ell+1}^{(\ell)} \circ \vare_\ell}
      \Bigg)\\
  &=  \ee^{-\ii \lambda \big( \widetilde\theta_{\lfloor nt\rfloor - 1}  a^{\lfloor nt\rfloor-1-\ell} + \widetilde\theta_{\lfloor nt\rfloor}  a^{\lfloor nt\rfloor-\ell}\big) }\\
  &\phantom{=}\times\EE\Bigg(\exp\Bigg\{\ii\Big( \widetilde\theta_1 (a_1^{(\ell)}\circ \cdots \circ a^{(\ell)}_{\ell+1}\circ\vare_\ell - \lambda a^{1-\ell}) +\cdots
     +\widetilde\theta_{\lfloor nt\rfloor -2}( a_{\lfloor nt\rfloor-2}^{(\ell)} \circ \cdots \circ a_{\ell+1}^{(\ell)} \circ \vare_\ell - \lambda a^{\lfloor nt\rfloor -2-\ell} )
        \Big) \Bigg\}\\
  & \phantom{= \times\EE\Bigg(}
      \times G_{\textrm{Bin}\big(a_{\lfloor nt\rfloor-2}^{(\ell)} \circ \cdots \circ a_{\ell+1}^{(\ell)} \circ \vare_\ell,a\big)}
         \big( (1-a)\ee^{\ii\widetilde\theta_{\lfloor nt\rfloor-1}}  + a \ee^{\ii(\widetilde\theta_{\lfloor nt\rfloor-1} + \widetilde\theta_{\lfloor nt\rfloor})}\big)
      \Bigg)
 \end{align*}
 \begin{align*}
  &=  \ee^{-\ii \lambda \big( \widetilde\theta_{\lfloor nt\rfloor - 1}  a^{\lfloor nt\rfloor-1-\ell} + \widetilde\theta_{\lfloor nt\rfloor}  a^{\lfloor nt\rfloor-\ell}\big) }\\
  &\phantom{=}\times\EE\Bigg(\exp\Bigg\{\ii\Big(\widetilde\theta_1 (a_1^{(\ell)}\circ \cdots \circ a^{(\ell)}_{\ell+1}\circ\vare_\ell - \lambda a^{1-\ell}) +\cdots
     +\widetilde\theta_{\lfloor nt\rfloor -2}( a_{\lfloor nt\rfloor-2}^{(\ell)} \circ \cdots \circ a_{\ell+1}^{(\ell)} \circ \vare_\ell - \lambda a^{\lfloor nt\rfloor -2-\ell} )
        \Big) \Bigg\}\\
  & \phantom{= \times\EE\Bigg(}
      \times
      \bigg( 1-a + a(1-a)\ee^{\ii\widetilde\theta_{\lfloor nt\rfloor-1}}
         + a^2 \ee^{\ii(\widetilde\theta_{\lfloor nt\rfloor-1} + \widetilde\theta_{\lfloor nt\rfloor})}\bigg)^{a_{\lfloor nt\rfloor-2}^{(\ell)} \circ \cdots \circ a_{\ell+1}^{(\ell)} \circ \vare_\ell}
      \Bigg)\\
  & = \cdots =\\
  &= \ee^{-\ii\lambda \sum_{k=1}^{\lfloor nt\rfloor} \widetilde\theta_k a^{k-\ell}}
          G_{a_1^{(\ell)}\circ \cdots \circ a^{(\ell)}_{\ell+1} \circ \vare_\ell}\Big( (1-a) \sum_{k=1}^{\lfloor nt\rfloor-1}  a^{k-1} \ee^{\ii\widetilde\theta_{1,k}}
                       + a^{\lfloor nt\rfloor-1}\ee^{\ii\widetilde\theta_{1,\lfloor nt\rfloor}} \Big).
 \end{align*}
Since for all \ $\ell\leq 0$, \ $\ell\in\ZZ$, \ and \ $a\in[0,1]$, \
 \[
    a_1^{(\ell)}\circ \cdots \circ a^{(\ell)}_{\ell+1} \circ \vare_\ell \distre \textrm{Bin}(\vare_\ell, a^{1-\ell}),
 \]
 which can be checked by calculating the characteristic function of both sides
 (see also Turkman et al. \cite[Lemma 5.1.1]{TurScoBer}),
 using again the tower rule, we have
 \begin{align*}
  &\EE\Bigg( \exp\Bigg\{ \ii  \sum_{k=\ell\vee 1}^{\lfloor nt\rfloor}
            \widetilde \theta_k \Big( a_k^{(\ell)} \circ \cdots \circ a_{\ell+1}^{(\ell)} \circ \vare_\ell - \lambda a^{k-\ell} \Big) \Bigg\} \Bigg) \\
  &= \ee^{-\ii\lambda \sum_{k=1}^{\lfloor nt\rfloor} \widetilde\theta_k a^{k-\ell}}
     G_{\vare_\ell} \Big( 1 - a^{1-\ell} + a^{1-\ell}\Big( (1-a) \sum_{k=1}^{\lfloor nt\rfloor-1}  a^{k-1} \ee^{\ii\widetilde\theta_{1,k}}
                       + a^{\lfloor nt\rfloor-1}\ee^{\ii\widetilde\theta_{1,\lfloor nt\rfloor}} \Big)\Big)\\
  &= \ee^{-\ii\lambda \sum_{k=1}^{\lfloor nt\rfloor} \widetilde\theta_k a^{k-\ell}}
     \exp\Big\{ -\lambda + \lambda \Big( 1 - a^{1-\ell} + a^{1-\ell}\Big( (1-a) \sum_{k=1}^{\lfloor nt\rfloor-1}  a^{k-1} \ee^{\ii\widetilde\theta_{1,k}}
                       + a^{\lfloor nt\rfloor-1}\ee^{\ii\widetilde\theta_{1,\lfloor nt\rfloor}} \Big)  \Big) \Big\}\\
  & = \ee^{-\ii\lambda \sum_{k=1}^{\lfloor nt\rfloor} \widetilde\theta_k a^{k-\ell}}
      \exp\Big\{ - \lambda a^{1-\ell} + \lambda  (1-a) \sum_{k=1}^{\lfloor nt\rfloor-1}  a^{k-\ell} \ee^{\ii\widetilde\theta_{1,k}}
                       + \lambda a^{\lfloor nt\rfloor-\ell}\ee^{\ii\widetilde\theta_{1,\lfloor nt\rfloor}}  \Big\},
 \end{align*}
 which coincides with
 \begin{align}\label{help_joint_characteristic3}
    \exp\Bigg\{\lambda
   \Bigg[(1-a)\sum_{k=1}^{\lfloor nt \rfloor-1}
               a^{k-\ell} (\ee^{\ii \widetilde \theta_{1,k}} - 1 - \ii \widetilde \theta_{1,k})
               + a^{\lfloor nt \rfloor -\ell} (\ee^{\ii \widetilde \theta_{1,\lfloor nt \rfloor}} - 1 - \ii \widetilde \theta_{1,\lfloor nt \rfloor})
         \Bigg]\Bigg\},
 \end{align}
 as desired.
Indeed, the coefficient of the constant term (not depending on \ $\widetilde \theta_j$, \ $j=1,\ldots,\lfloor nt \rfloor$) \ in the exponential
 of \eqref{help_joint_characteristic3} is
 \begin{align*}
   -\lambda(1-a)\sum_{k=1}^{\lfloor nt \rfloor - 1} a^{k-\ell} - \lambda a^{\lfloor nt \rfloor - \ell}
      = -\lambda(1-a) a^{1-\ell}\frac{a^{\lfloor nt \rfloor - 1} - 1}{a-1} - \lambda a^{\lfloor nt \rfloor - \ell}
      = -\lambda a^{1-\ell},
 \end{align*}
 the coefficient of \ $\widetilde \theta_j$, \ $j\in\{1,\ldots,\lfloor nt \rfloor-1\}$, \ in the exponential of \eqref{help_joint_characteristic3} is
 \[
  -\ii \lambda (1-a) \sum_{k=j}^{\lfloor nt \rfloor - 1} a^{k-\ell} - \ii \lambda a^{\lfloor nt \rfloor-\ell}
    = -\ii\lambda (1-a) a^{j-\ell} \frac{a^{\lfloor nt \rfloor - j} - 1}{a-1} - \ii \lambda a^{\lfloor nt \rfloor-\ell}
    = -\ii \lambda a^{j-\ell},
 \]
the coefficient of \ $\widetilde \theta_{\lfloor nt \rfloor}$ \ in the exponential of \eqref{help_joint_characteristic3} is
 \ $-\ii\lambda a^{\lfloor nt \rfloor - \ell}$,
 \ the coefficient of \ $\ee^{\ii \widetilde \theta_{1, \lfloor nt \rfloor}}$ \ in the exponential of \eqref{help_joint_characteristic3} is
 \ $\lambda a^{\lfloor nt \rfloor - \ell}$, \ and the remaining term \ $\lambda(1-a)\sum_{k=1}^{\lfloor nt \rfloor-1}
 a^{k-\ell} \ee^{\ii \widetilde \theta_{1,k}}$ \ coincide as well.
\proofend

\medskip

\noindent{\bf Direct proof of the equality of formulae \eqref{egyuttes_char_aggregalt1} and \eqref{egyuttes_char_aggregalt2}
 in Proposition \ref{Pro_egyuttes_char_aggregalt} in case of \ $\mathbf{m=1}$.}
\ In case of \ $m=1$, for all \ $a\in(0,1)$ \ and \ $n\in\NN$, \ we have
 \begin{align*}
  \frac{1}{1-a}K_n(a)
   & = \frac{1}{1-a}
      \Bigg[
      (\ee^{\ii \theta_{1,1}} - 1-\ii \theta_{1,1})
               (\lfloor nt_1\rfloor - 0)
       + \sum_{1\leq k_1<k_2\leq\lfloor nt_1\rfloor}
          a^{k_2-k_1}
          \bigl(\ee^{\ii \theta_{1,1}} - 1\bigr)^2
          \ee^{\ii (k_2-k_1-1)\theta_{1,1}}
        \Bigg] \\
  & = \frac{1}{1-a} \Bigg[  (\ee^{\ii \theta_1} - 1 - \ii \theta_1)\lfloor nt_1\rfloor
                            + \sum_{1\leq k_1<k_2\leq\lfloor nt_1\rfloor}
                              a^{k_2-k_1}
                            \bigl(\ee^{\ii \theta_1} - 1\bigr)^2
                             \ee^{\ii (k_2-k_1-1)\theta_1}
                    \Bigg],
 \end{align*}
and
 \begin{align*}
  \widetilde K_n(a)
       &=  \sum_{k=1}^{\lfloor nt_1 \rfloor-1} a^k (\ee^{\ii \widetilde \theta_{1,k}} - 1 - \ii \widetilde \theta_{1,k})
           + \frac{a^{\lfloor nt_1 \rfloor}}{1-a}
            (\ee^{\ii \widetilde \theta_{1,\lfloor nt_1 \rfloor}} - 1 - \ii \widetilde \theta_{1,\lfloor nt_1 \rfloor}) \\
       &\quad   +\sum_{\ell=1}^{\lfloor nt_1 \rfloor}
         \Bigg[(1-a)\sum_{k=\ell}^{\lfloor nt_1 \rfloor-1}
               a^{k-\ell} (\ee^{\ii \widetilde \theta_{\ell,k}} - 1 - \ii \widetilde \theta_{\ell,k})
               + a^{\lfloor nt_1 \rfloor -\ell} (\ee^{\ii \widetilde \theta_{\ell,\lfloor nt_1 \rfloor}} - 1 - \ii \widetilde \theta_{\ell,\lfloor nt_1 \rfloor})
         \Bigg] ,
 \end{align*}
 where
 \[
   \widetilde\theta_{\ell,k} = \widetilde\theta_\ell + \cdots + \widetilde\theta_k = (k-\ell +1)\theta_1,
     \qquad 1\leq \ell\leq k\leq \lfloor nt_1 \rfloor,
 \]
 since \ $\widetilde \theta_j=\theta_1 \bone_{\{ j\leq \lfloor nt_1 \rfloor \}} = \theta_1$ \ for every \ $j=1,\ldots,\lfloor nt_1 \rfloor$.
\ Next, we derive a simpler form of \ $\widetilde K_n(a)$.
\ Namely,
 \begin{align*}
  \widetilde K_n(a)
      &=  \sum_{k=1}^{\lfloor nt_1 \rfloor-1} a^k (\ee^{\ii k \theta_1 } - 1 - \ii k\theta_1)
           + \frac{a^{\lfloor nt_1 \rfloor}}{1-a}
            \big(\ee^{\ii \lfloor nt_1 \rfloor \theta_1 } - 1 - \ii \lfloor nt_1 \rfloor \theta_1 \big) \\
      &\phantom{=\;}  + (1-a)\sum_{\ell=1}^{\lfloor nt_1 \rfloor}
                 \sum_{k=\ell}^{\lfloor nt_1 \rfloor-1}
                 a^{k-\ell} \big(\ee^{\ii (k-\ell +1)\theta_1} - 1 - \ii (k-\ell +1)\theta_1 \big)\\
      &\phantom{=\;} + \sum_{\ell=1}^{\lfloor nt_1 \rfloor} a^{\lfloor nt_1 \rfloor -\ell}
                  \big(\ee^{\ii (\lfloor nt_1 \rfloor - \ell + 1)\theta_1 } - 1 - \ii (\lfloor nt_1 \rfloor - \ell + 1)\theta_1 \big) \\
      &=  \sum_{k=1}^{\lfloor nt_1 \rfloor-1} a^k (\ee^{\ii k \theta_1 } - 1 - \ii k\theta_1)
           + \frac{a^{\lfloor nt_1 \rfloor}}{1-a}
            \big(\ee^{\ii \lfloor nt_1 \rfloor \theta_1 } - 1 - \ii \lfloor nt_1 \rfloor \theta_1 \big) \\
      &\phantom{=\;}  + (1-a)\sum_{k=1}^{\lfloor nt_1 \rfloor -1}
                 \sum_{\ell=1}^k
                 a^{k-\ell} \big(\ee^{\ii (k-\ell +1)\theta_1} - 1 - \ii (k-\ell +1)\theta_1 \big)\\
      &\phantom{=\;} + \sum_{j=0}^{\lfloor nt_1 \rfloor-1} a^j
                  \big(\ee^{\ii (j + 1)\theta_1 } - 1 - \ii (j + 1)\theta_1 \big)
 \end{align*}
 \begin{align*}
      &=  \sum_{k=1}^{\lfloor nt_1 \rfloor-1} a^k (\ee^{\ii k \theta_1 } - 1 - \ii k\theta_1)
           + \frac{a^{\lfloor nt_1 \rfloor}}{1-a}
            \big(\ee^{\ii \lfloor nt_1 \rfloor \theta_1 } - 1 - \ii \lfloor nt_1 \rfloor \theta_1 \big) \\
      &\phantom{=\;}  + (1-a)\sum_{k=1}^{\lfloor nt_1 \rfloor -1}
                 \sum_{j=0}^{k-1}
                 a^j \big(\ee^{\ii (j+1)\theta_1} - 1 - \ii (j+1)\theta_1 \big)
          + \sum_{j=1}^{\lfloor nt_1 \rfloor} a^{j-1}
                  \big(\ee^{\ii j \theta_1 } - 1 - \ii j \theta_1 \big)\\
      &=  \sum_{k=1}^{\lfloor nt_1 \rfloor-1} a^k (\ee^{\ii k \theta_1 } - 1 - \ii k\theta_1)
           + \frac{a^{\lfloor nt_1 \rfloor}}{1-a}
            \big(\ee^{\ii \lfloor nt_1 \rfloor \theta_1 } - 1 - \ii \lfloor nt_1 \rfloor \theta_1 \big) \\
      &\phantom{=\;}  + (1-a)\sum_{k=1}^{\lfloor nt_1 \rfloor -1}
                 \sum_{j=1}^k
                 a^{j-1} \big(\ee^{\ii j\theta_1} - 1 - \ii j\theta_1 \big)
      + \sum_{j=1}^{\lfloor nt_1 \rfloor} a^{j-1}
                  \big(\ee^{\ii j \theta_1 } - 1 - \ii j \theta_1 \big).
 \end{align*}
Next, by induction with respect to \ $p\in\NN$, \ we prove that
 \begin{align}\label{help_formula_eq1}
 \begin{split}
    &\frac{1}{1-a} \Bigg[  (\ee^{\ii \theta_1} - 1 - \ii \theta_1)p
                            + \sum_{1\leq k_1<k_2\leq p}
                              a^{k_2-k_1}
                            \bigl(\ee^{\ii \theta_1} - 1\bigr)^2
                             \ee^{\ii (k_2-k_1-1)\theta_1}
                    \Bigg]\\
    &\qquad =  \sum_{k=1}^{p-1} a^k (\ee^{\ii k \theta_1 } - 1 - \ii k\theta_1)
           + \frac{a^p}{1-a}
            \big(\ee^{\ii p \theta_1 } - 1 - \ii p \theta_1 \big) \\
      &\phantom{\qquad =\;}  + (1-a)\sum_{k=1}^{p -1}
                 \sum_{j=1}^k
                 a^{j-1} \big(\ee^{\ii j\theta_1} - 1 - \ii j\theta_1 \big)
         + \sum_{j=1}^p a^{j-1}
                  \big(\ee^{\ii j \theta_1 } - 1 - \ii j \theta_1 \big)
 \end{split}
 \end{align}
 for all \ $a\in(0,1)$, \ which yields that \ $\frac{1}{1-a}K_n(a) = \widetilde K_n(a)$, \ $n\in\NN$, \ $a\in(0,1)$, \
 in case of \ $m=1$, \ as desired.
For \ $p=1$, \ \eqref{help_formula_eq1} takes the form
 \[
   \frac{1}{1-a} (\ee^{\ii\theta_1} - 1 - \ii\theta_1)
      = \frac{a}{1-a} (\ee^{\ii\theta_1} - 1 - \ii\theta_1)
        + (\ee^{\ii\theta_1} - 1 - \ii\theta_1),
 \]
 which readily holds.
Let us suppose that \eqref{help_formula_eq1} holds for \ $1,\ldots,p$, \ where \ $p\in\NN$.
\ Then, using the induction hypothesis, the left-hand side of \eqref{help_formula_eq1} with \ $p$ \ replaced by \ $p+1$ \ takes the form
 \begin{align*}
   &\frac{1}{1-a} \Bigg[  (\ee^{\ii \theta_1} - 1 - \ii \theta_1)(p+1)
                            + \sum_{1\leq k_1<k_2\leq p+1}
                              a^{k_2-k_1}
                            \bigl(\ee^{\ii \theta_1} - 1\bigr)^2
                             \ee^{\ii (k_2-k_1-1)\theta_1}
                    \Bigg]\\
  & = \frac{1}{1-a} \Bigg[  (\ee^{\ii \theta_1} - 1 - \ii \theta_1)p
                            + \sum_{1\leq k_1<k_2\leq p}
                              a^{k_2-k_1}
                            \bigl(\ee^{\ii \theta_1} - 1\bigr)^2
                             \ee^{\ii (k_2-k_1-1)\theta_1}
                    \Bigg]
      + \frac{\ee^{\ii\theta_1} - 1 - \ii\theta_1}{1-a}\\
  &\phantom{=\;}
      + \frac{1}{1-a} \sum_{k_1=1}^p
         a^{p+1-k_1} (\ee^{\ii \theta_1} - 1)^2
         \ee^{\ii(p + 1 - k_1 -1)\theta_1}=
\end{align*}
\begin{align*}
 &= \sum_{k=1}^{p-1} a^k (\ee^{\ii k \theta_1 } - 1 - \ii k\theta_1)
           + \frac{a^p}{1-a}
            \big(\ee^{\ii p \theta_1 } - 1 - \ii p \theta_1 \big) \\
 &\phantom{=\;}  + (1-a)\sum_{k=1}^{p -1}
                 \sum_{j=1}^k
                 a^{j-1} \big(\ee^{\ii j\theta_1} - 1 - \ii j\theta_1 \big)
         + \sum_{j=1}^p a^{j-1}
                  \big(\ee^{\ii j \theta_1 } - 1 - \ii j \theta_1 \big)\\
 &\phantom{=\;} + \frac{\ee^{\ii\theta_1} - 1 - \ii\theta_1}{1-a}
      + \frac{1}{1-a} \sum_{k_1=1}^p
         a^{p+1-k_1} (\ee^{\ii \theta_1} - 1)^2
         \ee^{\ii(p - k_1)\theta_1}.
 \end{align*}
The right-hand side of \eqref{help_formula_eq1} with \ $p$ \ replaced by \ $p+1$ \ takes the form
 \begin{align*}
   & \sum_{k=1}^p a^k (\ee^{\ii k \theta_1 } - 1 - \ii k\theta_1)
           + \frac{a^{p+1}}{1-a}
            \big(\ee^{\ii (p+1) \theta_1 } - 1 - \ii (p+1) \theta_1 \big) \\
  & + (1-a)\sum_{k=1}^p
                 \sum_{j=1}^k
                 a^{j-1} \big(\ee^{\ii j\theta_1} - 1 - \ii j\theta_1 \big)
         + \sum_{j=1}^{p+1} a^{j-1}
                  \big(\ee^{\ii j \theta_1 } - 1 - \ii j \theta_1 \big) ,
 \end{align*}
 so to prove that \eqref{help_formula_eq1} holds with \ $p$ \ replaced by \ $p+1$, \ it is enough to check that
 \begin{align*}
   &a^p(\ee^{\ii p \theta_1 } - 1 - \ii p\theta_1)
     + \frac{a^{p+1}}{1-a} \big(\ee^{\ii (p+1) \theta_1 } - 1 - \ii (p+1) \theta_1 \big)
     - \frac{a^p}{1-a}(\ee^{\ii p \theta_1 } - 1 - \ii p\theta_1) \\
   & + (1-a) \sum_{j=1}^p a^{j-1} \big(\ee^{\ii j\theta_1} - 1 - \ii j\theta_1 \big)
     + a^p\big(\ee^{\ii (p+1) \theta_1 } - 1 - \ii (p+1) \theta_1 \big) \\
   & - \frac{\ee^{\ii\theta_1} - 1 - \ii\theta_1}{1-a}
     - \frac{1}{1-a} \sum_{k_1=1}^p
         a^{p+1-k_1} (\ee^{\ii \theta_1} - 1)^2
         \ee^{\ii(p - k_1)\theta_1}\\
  & =0,\qquad a\in(0,1), \;\; \theta_1\in\RR,
 \end{align*}
 or equivalently
 \begin{align*}
  &-\frac{a^{p+1}}{1-a}  (\ee^{\ii p \theta_1} - 1 - \ii p\theta_1)
  + \frac{a^p}{1-a} \big(\ee^{\ii (p+1) \theta_1 } - 1 - \ii (p+1) \theta_1 \big)
  + (1-a) \sum_{j=1}^p a^{j-1} \big(\ee^{\ii j\theta_1} - 1 - \ii j\theta_1 \big)\\
  &- \frac{\ee^{\ii\theta_1} - 1 - \ii\theta_1}{1-a}
  - \frac{a}{1-a} \frac{(a\ee^{\ii\theta_1})^p - 1}{a\ee^{\ii\theta_1} - 1} (\ee^{\ii \theta_1} - 1)^2
  =0, \qquad a\in(0,1), \;\; \theta_1\in\RR.
 \end{align*}
After multiplying both sides by \ $(1-a)(a\ee^{\ii\theta_1} - 1)$, \
 to prove that \eqref{help_formula_eq1} holds it remains to verify that
\begin{align}\label{help_formula_eq2}
\begin{split}
  &-(\ee^{\ii p \theta_1} - 1 - \ii p\theta_1) a^{p+1}(a\ee^{\ii\theta_1} - 1)
   + \big(\ee^{\ii (p+1) \theta_1 } - 1 - \ii (p+1) \theta_1 \big) a^p (a\ee^{\ii\theta_1} - 1)\\
  &+ (1-a)^2 (a\ee^{\ii\theta_1} - 1)  \sum_{j=1}^p a^{j-1} \big(\ee^{\ii j\theta_1} - 1 - \ii j\theta_1 \big)
   - (a\ee^{\ii\theta_1} - 1)(\ee^{\ii\theta_1} - 1 - \ii\theta_1) \\
  &- a(a^p\ee^{\ii p\theta_1} - 1)(\ee^{\ii 2\theta_1} - 2\ee^{\ii\theta_1} +1 )
   = 0, \qquad a\in(0,1), \;\; \theta_1\in\RR.
 \end{split}
\end{align}
 Since
 \[
    \sum_{j=1}^p a^{j-1} \big(\ee^{\ii j\theta_1} - 1 - \ii j\theta_1 \big)
      = \ee^{\ii\theta_1} \frac{a^p \ee^{\ii p\theta_1} -1}{a\ee^{\ii\theta_1} - 1}
         - \frac{a^p-1}{a-1} - \ii\theta_1 \sum_{j=1}^p ja^{j-1},
 \]
 we get that \eqref{help_formula_eq2} is equivalent to
\begin{align}\label{help_formula_eq3}
\begin{split}
  &-(\ee^{\ii p \theta_1} - 1 - \ii p\theta_1) a^{p+1}(a\ee^{\ii\theta_1} - 1)
   + \big(\ee^{\ii (p+1) \theta_1 } - 1 - \ii (p+1) \theta_1 \big) a^p (a\ee^{\ii\theta_1} - 1)\\
  &+ (1-a)^2 \ee^{\ii\theta_1}(a^p \ee^{\ii p\theta_1} -1)
   + (1-a)(a\ee^{\ii\theta_1} - 1)(a^p-1)
   - \ii\theta_1 (1-a)^2 (a\ee^{\ii\theta_1} - 1) \sum_{j=1}^p ja^{j-1}\\
  & - (a\ee^{\ii\theta_1} - 1)(\ee^{\ii\theta_1} - 1 - \ii\theta_1)
    - a(a^p\ee^{\ii p\theta_1} - 1)(\ee^{\ii 2\theta_1} - 2\ee^{\ii\theta_1} +1 )
   = 0, \qquad a\in(0,1), \;\; \theta_1\in\RR.
 \end{split}
\end{align}
The validity of \eqref{help_formula_eq3} can be checked by calculating the coefficients of
 \ $\ee^{\ii(p+2)\theta_1}$, \ $\ee^{\ii(p+1)\theta_1}$, \ $\ee^{\ii p \theta_1}$, \ $\theta_1\ee^{\ii\theta_1}$,
 \ $\ee^{\ii 2\theta_1}$, \ $\ee^{\ii \theta_1}$, \ $\theta_1$, \ and the constant term (not depending on \ $\theta_1$),
 and verifying that these are all \ $0$.
\ We provide the details for \ $\ee^{\ii(p+1)\theta_1}$, \ $\theta_1\ee^{\ii\theta_1}$, \ and \ $\theta_1$.
\ The coefficient of \ $\ee^{\ii(p+1)\theta_1}$ \  on the left-hand side of \eqref{help_formula_eq3} is
 \[
   -a^{p+2} - a^p + (1-a)^2 a^p + 2a^{p+1} = 0,
 \]
 the coefficient of \ $\theta_1\ee^{\ii\theta_1}$ \  at the left-hand side of \eqref{help_formula_eq3} is
 \begin{align*}
  &\ii p a^{p+2} -\ii (p+1)a^{p+1} - \ii(a-2a^2 + a^3) \sum_{j=1}^p j a^{j-1} + \ii a\\
  &  = \ii p a^{p+2} -\ii (p+1)a^{p+1} + \ii a
    - \ii \left( \sum_{j=1}^p j a^j - 2 \sum_{j=1}^p j a^{j+1} + \sum_{j=1}^p j a^{j+2}\right) \\
  & = \ii p a^{p+2} -\ii (p+1)a^{p+1} + \ii a
    - \ii \left( \sum_{j=1}^p j a^j - 2 \sum_{j=2}^{p+1} (j-1) a^j + \sum_{j=3}^{p+2} (j-2) a^j\right)
   =0,
 \end{align*}
 and the coefficient of \ $\theta_1$ \  at the left-hand side of \eqref{help_formula_eq3} is
 \begin{align*}
 & -\ii p a^{p+1} + \ii(p+1) a^p + \ii(1-2a+a^2) \sum_{j=1}^p j a^{j-1} -\ii\\
 & =-\ii p a^{p+1} + \ii(p+1) a^p -\ii +\ii\left( \sum_{j=1}^p j a^{j-1} - 2 \sum_{j=2}^{p+1} (j-1) a^{j-1} + \sum_{j=3}^{p+2} (j-2) a^{j-1}\right)
   =0.
 \end{align*}
 \proofend

\noindent{\bf Proof of Theorem \ref{simultaneous_aggregation_random_2}.}
To prove this limit theorem we have to show that for any sequence
 \ $(N_n)_{n\in\NN}$ \ of positive integers with \ $N_n^{\frac{-\beta}{1+\beta}} n^{-1} \to \infty$, \ we have
 \[
   n^{-1} N_n^{-\frac{1}{2(1+\beta)}} \,
   S^{(N_n,n)}
   \distrf (V_{2(1+\beta)} t)_{t\in\RR_+} \qquad
   \text{as \ $n \to \infty$.}
 \]
For this, by continuous mapping theorem, it is enough to verify that for any
 \ $m \in \NN$ \ and \ $t_0, t_1, \ldots, t_m \in \RR_+$ \ with
 \ $0 =: t_0 < t_1 < \ldots < t_m$, \ we have
 \begin{align*}
  &n^{-1} N_n^{-\frac{1}{2(1+\beta)}}
   \sum_{j=1}^{N_n}
    \Biggl(\sum_{k=1}^{\lfloor nt_1\rfloor}
            \Bigl(X_k^{(j)}-\frac{\lambda}{1-\alpha^{(j)}}\Bigr),
           \sum_{k=\lfloor nt_1\rfloor+1}^{\lfloor nt_2\rfloor}
            \Bigl(X_k^{(j)}-\frac{\lambda}{1-\alpha^{(j)}}\Bigr), \dots,
           \sum_{k=\lfloor nt_{m-1}\rfloor+1}^{\lfloor nt_m\rfloor}
            \Bigl(X_k^{(j)}-\frac{\lambda}{1-\alpha^{(j)}}\Bigr)\Biggr) \\
  &\distr V_{2(1+\beta)}(t_1, t_2 - t_1, \dots, t_m - t_{m-1}) \qquad
   \text{as \ $n \to \infty$.}
 \end{align*}
So, by continuity theorem, we have to check that for any \ $m \in \NN$,
 \ $t_0, t_1, \ldots, t_m \in \RR_+$ \ with \ $0 = t_0 < t_1 < \ldots < t_m$ \ and
 \ $\theta_1, \dots, \theta_m \in \RR$ \ the convergence
 \begin{align*}
   &\EE\Biggl(\exp\Biggl\{\ii
                          \sum_{\ell=1}^m
                           \theta_\ell n^{-1} N_n^{-\frac{1}{2(1+\beta)}}
                           \sum_{j=1}^{N_n}
                            \sum_{k=\lfloor nt_{\ell-1}\rfloor+1}^{\lfloor nt_\ell\rfloor}
                             \Bigl(X_k^{(j)}
                                   - \frac{\lambda}
                                          {1-\alpha^{(j)}}\Bigr)\Biggr\}\Biggr) \\
   &\qquad
    =\EE\Biggl(\exp\Biggl\{\ii n^{-1} N_n^{-\frac{1}{2(1+\beta)}}
                           \sum_{j=1}^{N_n}
                            \sum_{\ell=1}^m
                             \theta_\ell
                             \sum_{k=\lfloor nt_{\ell-1}\rfloor+1}^{\lfloor nt_\ell\rfloor}
                              \Bigl(X_k^{(j)}
                                    - \frac{\lambda}
                                           {1-\alpha^{(j)}}\Bigr)\Biggr\}\Biggr)\\
   &\qquad
    =\Biggl[\EE\Biggl(\exp\Biggl\{\ii n^{-1} N_n^{-\frac{1}{2(1+\beta)}}
                                  \sum_{\ell=1}^m
                                   \theta_\ell
                                   \sum_{k=\lfloor nt_{\ell-1}\rfloor+1}^{\lfloor nt_\ell\rfloor}
                                    \Bigl(X_k
                                          - \frac{\lambda}
     {1-\alpha}\Bigr)\Biggr\}\Biggr)\Biggr]^{N_n}\\
   &\qquad\to \EE\Biggl(\ee^{\ii
                       \sum_{\ell=1}^m
                        \theta_\ell (t_\ell - t_{\ell-1}) V_{2(1+\beta)}}\Biggr)
    = \ee^{-K_\beta |\sum_{\ell=1}^m \theta_\ell (t_\ell - t_{\ell-1})|^{2(1+\beta)}}
    \qquad \text{as \ $n \to \infty$}
 \end{align*}
 holds.
Note that it suffices to show
 \begin{align*}
  \Theta_n
  &:=N_n
     \Biggl[1 - \EE\Biggl(\exp\Biggl\{\ii n^{-1} N_n^{-\frac{1}{2(1+\beta)}}
                                    \sum_{\ell=1}^m
                                     \theta_\ell
                                      \sum_{k=\lfloor nt_{\ell-1}\rfloor+1}^{\lfloor nt_\ell\rfloor}
                                       \Bigl(X_k
                                             - \frac{\lambda}
                                                    {1-\alpha}\Bigr)\Biggr\}\Biggr)\Biggr] \\
  &\to K_\beta
       \Biggl|\sum_{\ell=1}^m \theta_\ell (t_\ell - t_{\ell-1})\Biggr|^{2(1+\beta)}
   \qquad \text{as \ $n \to \infty$,}
 \end{align*}
 since it implies that
 \ $(1 - \Theta_n/N_n)^{N_n}
    \to \ee^{-K_\beta |\sum_{\ell=1}^m \theta_\ell (t_\ell - t_{\ell-1})|^{2(1+\beta)}}$
 \ as \ $n \to \infty$, \ as desired.
By applying \eqref{help1_alter} (or \eqref{egyuttes_char_aggregalt1}) to the left hand side, we get
 \begin{equation*}
  \begin{split}
   \Theta_n
   &= N_n \EE\Biggl[1 -  \ee^{-\ii n^{-1} N_n^{-\frac{1}{2(1+\beta)}}\frac{\lambda}{1-\alpha}
            \sum_{\ell=1}^m \theta_\ell (\lfloor{nt_{\ell}}\rfloor - \lfloor nt_{\ell-1}\rfloor)} \\
   &\phantom{=\;}
     \times F_{0,\dots,\lfloor nt_m\rfloor-1}
      \Bigl( \underbrace{ \ee^{\ii n^{-1} N_n^{-\frac{1}{2(1+\beta)}}\theta_1},\ldots, \ee^{\ii n^{-1} N_n^{-\frac{1}{2(1+\beta)}}\theta_1} }_{\lfloor nt_1\rfloor \ \text{items}},
      \dots,
   \underbrace{ \ee^{\ii n^{-1} N_n^{-\frac{1}{2(1+\beta)}}\theta_m}, \ldots, \ee^{\ii n^{-1} N_n^{-\frac{1}{2(1+\beta)}}\theta_m}  }_{\lfloor nt_m\rfloor - \lfloor nt_{m-1}\rfloor \ \text{items}}  \,\Big|\, \alpha\Bigr)
         \Biggr]\\
   &= N_n \EE\left[1 - \ee^{\frac{\lambda}{1-\alpha}A_n(\alpha)}\right]
    = N_n \int_0^1
         \left(1 - \ee^{\frac{\lambda}{1-a}A_n(a)}\right) \psi(a) (1-a)^\beta
         \, \dd a
  \end{split}
 \end{equation*}
 with
 \begin{align*}
  &A_n(a)
   := \sum_{\ell=1}^m
                \Bigl(\ee^{\ii n^{-1} N_n^{-\frac{1}{2(1+\beta)}}\theta_\ell} - 1 - \ii n^{-1} N_n^{-\frac{1}{2(1+\beta)}} \theta_\ell \Bigr)
                (\lfloor nt_\ell\rfloor - \lfloor nt_{\ell-1}\rfloor) \\
   &+\sum_{1\leq\ell_1<\ell_2\leq m}
      \sum_{k_1=\lfloor nt_{\ell_1-1}\rfloor+1}^{\lfloor nt_{\ell_1}\rfloor}
       \sum_{k_2=\lfloor nt_{\ell_2-1}\rfloor+1}^{\lfloor nt_{\ell_2}\rfloor}
        a^{k_2-k_1}
        \bigl(\ee^{\ii n^{-1} N_n^{-\frac{1}{2(1+\beta)}}\theta_{\ell_1}} - 1\bigr)
        \bigl(\ee^{\ii n^{-1} N_n^{-\frac{1}{2(1+\beta)}}\theta_{\ell_2}} - 1\bigr) \\
   &\phantom{+\sum_{1\leq\ell_1<\ell_2\leq m}
               \sum_{k_1=\lfloor nt_{\ell_1-1}\rfloor+1}^{\lfloor nt_{\ell_1}\rfloor}
                \sum_{k_2=\lfloor nt_{\ell_2-1}\rfloor+1}^{\lfloor nt_{\ell_2}\rfloor}}
        \times
        \ee^{\ii n^{-1} N_n^{-\frac{1}{2(1+\beta)}}
             \bigl((\lfloor nt_{\ell_1}\rfloor-k_1)\theta_{\ell_1}
                   +\sum_{\ell=\ell_1+1}^{\ell_2-1}
                     \theta_\ell (\lfloor nt_\ell\rfloor - \lfloor nt_{\ell-1}\rfloor)
                   +(k_2 -1 - \lfloor nt_{\ell_2-1}\rfloor)\theta_{\ell_2}\bigr)} \\
   &+\sum_{\ell=1}^m
      \sum_{\lfloor nt_{\ell-1}\rfloor+1\leq k_1<k_2\leq\lfloor nt_\ell\rfloor}
       a^{k_2-k_1}
       \bigl(\ee^{\ii n^{-1} N_n^{-\frac{1}{2(1+\beta)}}\theta_\ell} - 1\bigr)^2
        \ee^{\ii n^{-1} N_n^{-\frac{1}{2(1+\beta)}} (k_2-k_1-1)\theta_\ell}
 \end{align*}
 for \ $a \in [0, 1]$.
\ The aim of the following discussion is to apply Lemma \ref{ordo} with \ $z_n(a):=A_n(a)$, \ $n\in\NN$, \ $a\in(0,1)$, \
 \ $\vare_n :=  N_n^{\frac{\beta}{1+\beta}}$, \ $n \in \NN$, \ and
 \[
   I:= \psi_1^{-1}K_\beta \left|\sum_{\ell=1}^m \theta_\ell (t_\ell - t_{\ell-1})\right|^{2(1+\beta)}.
 \]
Since \ $\beta\in(-1,0)$, \ we have \ $\vare_n\in(0,1)$ \ for \ $n\geq n_0$, \ where \ $n_0$ \ is sufficiently large, and \ $\lim_{n\to\infty}\vare_n=0$.
\ First we check \eqref{help_ordo1}.
Using \eqref{exp:3}, for any \ $a \in (0, 1)$ \ we get
\begin{align*}
   |A_n(a)|
   &\leq \sum_{\ell=1}^m
             n^{-2} N_n^{-\frac{1}{1+\beta}} \frac{\theta_\ell^2}{2}
             (\lfloor nt_\ell\rfloor - \lfloor nt_{\ell-1}\rfloor) \\
   &\quad
    + \sum_{1\leq\ell_1<\ell_2\leq m}
         n^{-2} N_n^{-\frac{1}{1+\beta}} |\theta_{\ell_1}| |\theta_{\ell_2}|
         (\lfloor nt_{\ell_1}\rfloor - \lfloor nt_{\ell_1-1}\rfloor)
         (\lfloor nt_{\ell_2}\rfloor - \lfloor nt_{\ell_2-1}\rfloor) \\
   &\quad
    + \sum_{\ell=1}^m
         n^{-2} N_n^{-\frac{1}{1+\beta}} \frac{\theta_\ell^2}{2}
         (\lfloor nt_\ell\rfloor - \lfloor nt_{\ell-1}\rfloor)
         (\lfloor nt_\ell\rfloor - \lfloor nt_{\ell-1}\rfloor - 1) \\
  &= \frac{1}{2} n^{-2} N_n^{-\frac{1}{1+\beta}}
      \biggl(\sum_{\ell=1}^m
              |\theta_\ell| (\lfloor nt_\ell\rfloor - \lfloor nt_{\ell-1}\rfloor)\biggr)^2
    \leq \frac{1}{2} N_n^{-\frac{1}{1+\beta}}
         \biggl(\sum_{\ell=1}^m |\theta_\ell| (t_\ell - t_{\ell-1} + 1)\biggr)^2 ,
 \end{align*}
 since \ $\frac{1}{n}(\lfloor nt_\ell\rfloor - \lfloor nt_{\ell-1}\rfloor) \leq \frac{1}{n}(nt_\ell - nt_{\ell-1} + 1) = t_\ell - t_{\ell-1} + \frac{1}{n}
  \leq t_\ell - t_{\ell-1} +1$.
\ Consequently, since \ $\vare_n^{-1} N_n = N_n^{\frac{1}{1+\beta}}$, \ we have
 \[
   \sup_{n\geq n_0} \vare_n^{-1} N_n \sup_{a\in(0,1-\vare_n)} |A_n(a)|
   \leq \frac{1}{2} \biggl(\sum_{\ell=1}^m |\theta_\ell| (t_\ell - t_{\ell-1} + 1)\biggr)^2
   < \infty ,
 \]
 i.e., \eqref{help_ordo1} is satisfied.
Therefore, by Lemma \ref{ordo}, substituting \ $a = 1 - z^{-1} N_n^{-\frac{1}{1+\beta}}$ \ with \ $z>0$,
 \ the statement of the theorem will follow from
 \begin{gather}\label{limsup_an}
  \begin{aligned}
   &\limsup_{n\to\infty}
     N_n \int^1_{1-N_n^{\frac{\beta}{1+\beta}}}
        \Bigl|1 - \ee^{\frac{\lambda}{1-a}A_n(a)}\Bigr| (1-a)^\beta \, \dd a \\
   &= \limsup_{n\to\infty}
       \int^\infty_{N_n^{-1}}
        \Bigl|1 - \ee^{\lambda zN_n^{\frac{1}{1+\beta}}
                       A_n\bigl(1 - z^{-1}N_n^{-\frac{1}{1+\beta}}\bigr)}\Bigr|
        z^{-(2+\beta)} \, \dd z
    < \infty
  \end{aligned}
 \end{gather}
 and
 \begin{gather}\label{lim_an}
  \begin{aligned}
   &\lim_{n\to\infty}
     \left|N_n \int_{1-N_n^{\frac{\beta}{1+\beta}}}^1
                \Bigl(1 - \ee^{\frac{\lambda}{1-a}A_n(a)}\Bigr) (1-a)^\beta \, \dd a
           - I\right| \\
   &= \lim_{n\to\infty}
        \left|\int^\infty_{N_n^{-1}}
               \Bigl(1 - \ee^{\lambda zN_n^{\frac{1}{1+\beta}}
                              A_n\bigl(1 - z^{-1}N_n^{-\frac{1}{1+\beta}}\bigr)}\Bigr)
               z^{-(2+\beta)} \, \dd z
             - I\right|
    = 0
  \end{aligned}
 \end{gather}
 with
 \begin{align*}
  I & = \psi_1^{-1} K_\beta \left|\sum_{\ell=1}^m \theta_\ell (t_\ell - t_{\ell-1})\right|^{2(1+\beta)}
      = \left(\frac{\lambda}{2}
             \left|\sum_{\ell=1}^m \theta_\ell (t_\ell - t_{\ell-1})\right|^2\right)^{1+\beta}
        \int_0^\infty (1-\ee^{-z}) z^{-(2+\beta)} \, \dd z \\
    & =  \int_0^\infty
         \biggl(1-\ee^{-\frac{\lambda z}{2}
                        \bigl(\sum_{\ell=1}^m \theta_\ell (t_\ell - t_{\ell-1})\bigr)^2}\biggr)
          z^{-(2+\beta)}
          \, \dd z ,
 \end{align*}
 where the first equality is justified by Lemma 2.2.1 in Zolotarev \cite{Zol} (be
 careful for the misprint in \cite{Zol}: a negative sign is superfluous) or by
 Li \cite[formula (1.28)]{Li}.

Next we check \eqref{limsup_an} and \eqref{lim_an}.
By Taylor expansion,
 \begin{gather*}
  \ee^{\ii n^{-1} N_n^{-\frac{1}{2(1+\beta)}} \theta_\ell} - 1
  = \ii n^{-1} N_n^{-\frac{1}{2(1+\beta)}} \theta_\ell + n^{-2} N_n^{-\frac{1}{1+\beta}} \OO(1)
  = n^{-1} N_n^{-\frac{1}{2(1+\beta)}} \OO(1), \\
  \ee^{\ii n^{-1} N_n^{-\frac{1}{2(1+\beta)}} \theta_\ell} - 1
  - \ii n^{-1} N_n^{-\frac{1}{2(1+\beta)}}\theta_{\ell}
  = - n^{-2} N_n^{-\frac{1}{1+\beta}} \frac{\theta_{\ell}^2}{2}
    + n^{-3} N_n^{-\frac{3}{2(1+\beta)}} \OO(1)
  =  n^{-2} N_n^{-\frac{1}{1+\beta}} \OO(1)
 \end{gather*}
 for all \ $\ell \in\{1, \ldots, m\}$, \ resulting
 \begin{equation}\label{an}
  \lambda z N_n^{\frac{1}{1+\beta}} A_n\biggl(1 - \frac{1}{zN_n^{\frac{1}{1+\beta}}}\biggr)
  = - \frac{\lambda z
            \bigl(\sum_{\ell=1}^m
                   \theta_\ell (\lfloor nt_\ell\rfloor - \lfloor nt_{\ell-1}\rfloor)\bigr)^2}
           {2n^2}
    + \frac{z\OO(1)}{N_n^{\frac{1}{2(1+\beta)}}}
    + \frac{n\OO(1)}{N_n^{\frac{1}{1+\beta}}}
 \end{equation}
 for \ $z > N_n^{-1}$.
\ Indeed, for \ $z > N_n^{-1}$, \ we have \ $z> N_n^{-\frac{1}{1+\beta}}$ \ yielding \ $1 - z^{-1}N_n^{-\frac{1}{1+\beta}}\in(0,1)$, \ and
 \begin{align*}
   &A_n\left(1-\frac{1}{zN_n^{\frac{1}{1+\beta}}}\right) \\
   &= \sum_{\ell=1}^m
       \bigl(\ee^{\ii n^{-1} N_n^{-\frac{1}{2(1+\beta)}}\theta_\ell} - 1
             - \ii n^{-1} N_n^{-\frac{1}{2(1+\beta)}}\theta_\ell\bigr)
       (\lfloor nt_\ell\rfloor - \lfloor nt_{\ell-1}\rfloor) \\
   &\quad
      +\sum_{1\leq\ell_1<\ell_2\leq m}
        \sum_{k_1=\lfloor nt_{\ell_1-1}\rfloor+1}^{\lfloor nt_{\ell_1}\rfloor}
         \sum_{k_2=\lfloor nt_{\ell_2-1}\rfloor+1}^{\lfloor nt_{\ell_2}\rfloor}
          \left(1 - \frac{1}{zN_n^{\frac{1}{1+\beta}}}\right)^{k_2-k_1}
          \bigl(\ee^{\ii n^{-1} N_n^{-\frac{1}{2(1+\beta)}}\theta_{\ell_1}} - 1\bigr) \\
   &\phantom{\quad+}\times
        \ee^{\ii n^{-1} N_n^{-\frac{1}{2(1+\beta)}}
             \bigl((\lfloor nt_{\ell_1}\rfloor-k_1)\theta_{\ell_1}
                   +\sum_{\ell=\ell_1+1}^{\ell_2-1}
                     \theta_\ell (\lfloor nt_\ell\rfloor - \lfloor nt_{\ell-1}\rfloor)
                   +(k_2  - 1- \lfloor nt_{\ell_2-1}\rfloor)\theta_{\ell_2}\bigr)}
        \bigl(\ee^{\ii n^{-1} N_n^{-\frac{1}{2(1+\beta)}}\theta_{\ell_2}} - 1\bigr)
 \end{align*}
 \begin{align*}
   &\quad+\sum_{\ell=1}^m
      \sum_{\lfloor nt_{\ell-1}\rfloor+1\leq k_1<k_2\leq\lfloor nt_\ell\rfloor}
       \left(1 - \frac{1}{zN_n^{\frac{1}{1+\beta}}}\right)^{k_2-k_1}
       \bigl(\ee^{\ii n^{-1} N_n^{-\frac{1}{2(1+\beta)}}\theta_\ell} - 1\bigr)^2
        \ee^{\ii n^{-1} N_n^{-\frac{1}{2(1+\beta)}} (k_2-k_1-1)\theta_\ell}\\
  & = \sum_{\ell=1}^m
       \biggl(- \frac{\theta_\ell^2}{2n^2N_n^{\frac{1}{1+\beta}}}
              + \frac{\OO(1)}{n^3N_n^{\frac{3}{2(1+\beta)}}}\biggr)
       (\lfloor nt_\ell\rfloor - \lfloor nt_{\ell-1}\rfloor) \\
   &\quad
      +\sum_{1\leq\ell_1<\ell_2\leq m}
        \left(1 + \frac{n\OO(1)}{zN_n^{\frac{1}{1+\beta}}}\right)
        \biggl(\frac{\ii\theta_{\ell_1}}{nN_n^{\frac{1}{2(1+\beta)}}}
               + \frac{\OO(1)}{n^2N_n^{\frac{1}{1+\beta}}}\biggr)
        \biggl(1 + \frac{\OO(1)}{N_n^{\frac{1}{2(1+\beta)}}}\biggr)
        \biggl(\frac{\ii\theta_{\ell_2}}{nN_n^{\frac{1}{2(1+\beta)}}}
               + \frac{\OO(1)}{n^2N_n^{\frac{1}{1+\beta}}}\biggr) \\
   &\phantom{\quad+\sum_{1\leq\ell_1<\ell_2\leq m}} \times
        (\lfloor nt_{\ell_1}\rfloor - \lfloor nt_{\ell_1-1}\rfloor)
        (\lfloor nt_{\ell_2}\rfloor - \lfloor nt_{\ell_2-1}\rfloor) \\
   &\quad
      +\frac{1}{2}\sum_{\ell=1}^m
        \left(1 + \frac{n\OO(1)}{zN_n^{\frac{1}{1+\beta}}}\right)
        \biggl(\frac{\ii\theta_\ell}{nN_n^{\frac{1}{2(1+\beta)}}}
               + \frac{\OO(1)}{n^2N_n^{\frac{1}{1+\beta}}}\biggr)^2
        \biggl(1 + \frac{\OO(1)}{N_n^{\frac{1}{2(1+\beta)}}}\biggr) \\
   &\phantom{=+\sum_{\ell=1}^m}\times
        (\lfloor nt_\ell\rfloor - \lfloor nt_{\ell-1}\rfloor)
        (\lfloor nt_\ell\rfloor - \lfloor nt_{\ell-1} - 1\rfloor)\\
  &= - \frac{\sum_{\ell=1}^m \theta_\ell^2 (\lfloor nt_\ell\rfloor - \lfloor nt_{\ell-1}\rfloor)}
             {2n^2N_n^{\frac{1}{1+\beta}}}
      + \frac{\OO(1)}{n^2N_n^{\frac{3}{2(1+\beta)}}} \\
   &\quad
      - \frac{\sum_{1\leq\ell_1<\ell_2\leq m}
               \theta_{\ell_1} \theta_{\ell_2}
               (\lfloor nt_{\ell_1}\rfloor - \lfloor nt_{\ell_1-1}\rfloor)
               (\lfloor nt_{\ell_2}\rfloor - \lfloor nt_{\ell_2-1}\rfloor)}
             {n^2N_n^{\frac{1}{1+\beta}}}
      + \frac{\OO(1)}{N_n^{\frac{3}{2(1+\beta)}}}
      + \frac{n\OO(1)}{zN_n^{\frac{2}{1+\beta}}} \\
   &\quad
      - \frac{\sum_{\ell=1}^m
               \theta_\ell^2 (\lfloor nt_\ell\rfloor - \lfloor nt_{\ell-1}\rfloor)
               (\lfloor nt_\ell\rfloor - \lfloor nt_{\ell-1} - 1\rfloor)}
             {2n^2N_n^{\frac{1}{1+\beta}}}
      + \frac{\OO(1)}{N_n^{\frac{3}{2(1+\beta)}}}
      + \frac{n\OO(1)}{zN_n^{\frac{2}{1+\beta}}} \\
   &= - \frac{\bigl(\sum_{\ell=1}^m
                     \theta_\ell (\lfloor nt_\ell\rfloor - \lfloor nt_{\ell-1}\rfloor)\bigr)^2}
             {2n^2N_n^{\frac{1}{1+\beta}}}
      + \frac{\OO(1)}{N_n^{\frac{3}{2(1+\beta)}}}
      + \frac{n\OO(1)}{zN_n^{\frac{2}{1+\beta}}} ,
 \end{align*}
 where we used the following facts:
 \begin{itemize}
 \item
 \begin{align}\label{help_ordo3}
 \begin{split}
  &\ee^{\ii n^{-1} N_n^{-\frac{1}{2(1+\beta)}}
             \bigl((\lfloor nt_{\ell_1}\rfloor-k_1)\theta_{\ell_1}
                   + \sum_{\ell=\ell_1+1}^{\ell_2-1}
                     \theta_\ell (\lfloor nt_\ell\rfloor - \lfloor nt_{\ell-1}\rfloor)
                   +(k_2 - 1 - \lfloor nt_{\ell_2-1}\rfloor)\theta_{\ell_2}\bigr)} \\
  &\qquad = \ee^{\ii N_n^{-\frac{1}{2(1+\beta)}} \OO(1)} = 1+ N_n^{-\frac{1}{2(1+\beta)}} \OO(1)
  \end{split}
 \end{align}
 \item
  \begin{align}\label{help_ordo4}
  \begin{split}
     \ee^{\ii n^{-1} N_n^{-\frac{1}{2(1+\beta)}} (k_2-k_1+1)\theta_\ell}
      =  \ee^{\ii N_n^{-\frac{1}{2(1+\beta)}} \OO(1)}
      = 1+ N_n^{-\frac{1}{2(1+\beta)}} \OO(1),
   \end{split}
 \end{align}
   due to \ $\lfloor nt_{\ell-1}\rfloor+1 \leq k_1< k_2 \leq \lfloor nt_\ell\rfloor$,
 \item
   \[
     \left(1 - \frac{1}{zN_n^{\frac{1}{1+\beta}}}\right)^{k_2-k_1}
      = 1 + \frac{n\OO(1)}{zN_n^{\frac{1}{1+\beta}}} ,
    \]
    following from an application of Bernoulli's inequality:
   \[
     \Bigg|\left(1 - \frac{1}{zN_n^{\frac{1}{1+\beta}}}\right)^{k_2-k_1} - 1\Bigg|
     \leq \frac{k_2-k_1}{zN_n^{\frac{1}{1+\beta}}}
     \leq \frac{\lfloor nt_m\rfloor}{zN_n^{\frac{1}{1+\beta}}} .
    \]
\end{itemize}
By \eqref{an}, for large enough \ $n$ \ and for any \ $z \in [1, \infty)$, \ we have
 \begin{align*}
  &\lambda z N_n^{\frac{1}{1+\beta}}
   \Re A_n\bigl(1 - z^{-1} N_n^{-\frac{1}{1+\beta}}\bigr) \\
  &= - \frac{\lambda z
             (\sum_{\ell=1}^m
               \theta_\ell (\lfloor nt_\ell\rfloor - \lfloor nt_{\ell-1}\rfloor))^2}
            {2n^2} \,
       \left(1 - \frac{\Re\OO(1)}{N_n^{\frac{1}{2(1+\beta)}}}\right)
     + \frac{n\Re\OO(1)}{N_n^{\frac{1}{1+\beta}}} \\
  &\leq - \frac{\lambda z (\sum_{\ell=1}^m \theta_\ell (t_\ell - t_{\ell-1}))^2}{4}
        + \frac{n|\OO(1)|}{N_n^{\frac{1}{1+\beta}}}
  \leq - \frac{\lambda (\sum_{\ell=1}^m \theta_\ell (t_\ell - t_{\ell-1}))^2}{4}
            + \frac{n|\OO(1)|}{N_n^{\frac{1}{1+\beta}}}
   \leq 0 ,
 \end{align*}
 since \ $N_n^{\frac{1}{2(1+\beta)}} \to\infty$ \ as \ $n\to\infty$, \ and
 \ $n N_n^{-\frac{1}{1+\beta}} \leq n N_n^{\frac{\beta}{1+\beta}} \to 0$ \ as
 \ $n \to \infty$, \ hence we obtain for large enough \ $n$,
 \begin{align}\label{help10A}
  \begin{split}
   &\int_1^\infty
     \left|1 - \ee^{\lambda z N_n^{\frac{1}{1+\beta}}
                    A_n(1 - z^{-1} N_n^{-\frac{1}{1+\beta}})}\right|
     z^{-(\beta+2)} \, \dd z \\
   &\leq \int_1^\infty
          \left(1 + \ee^{\lambda z N_n^{\frac{1}{1+\beta}}
                         \Re A_n(1 - z^{-1} N_n^{-\frac{1}{1+\beta}})}\right)
          z^{-(\beta+2)} \, \dd z
    \leq 2 \int_1^\infty z^{-(\beta+2)} \, \dd z
    < \infty .
  \end{split}
 \end{align}
Again by \eqref{an}, for large enough \ $n$ \ and for any \ $z \in \bigl(N_n^{-1}, 1\bigr]$,
 \ we have
 \begin{align*}
  &\Bigl|\lambda z N_n^{\frac{1}{1+\beta}}
         A_n\Bigl(1 - z^{-1} N_n^{-\frac{1}{1+\beta}}\Bigr)\Bigr|
   \leq \frac{\lambda z
              (\sum_{\ell=1}^m
                \theta_\ell (\lfloor nt_\ell\rfloor - \lfloor nt_{\ell-1}\rfloor))^2}
             {2n^2}
        + \frac{z|\OO(1)|}{N_n^{\frac{1}{2(1+\beta)}}}
        + \frac{n|\OO(1)|}{N_n^{\frac{1}{1+\beta}}} \\
  &\leq z \biggl(\frac{\lambda(\sum_{\ell=1}^m \vert \theta_\ell \vert (t_\ell - t_{\ell-1} +1))^2}{2}
                 + \frac{|\OO(1)|}{N_n^{\frac{1}{2(1+\beta)}}}
                 + |\OO(1)| \biggr)
   \leq  z |\OO(1)|
   \leq |\OO(1)| ,
 \end{align*}
 where we used that \ $z \in \bigl(N_n^{-1}, 1\bigr]$ \ and \ $n N_n^{\frac{\beta}{1+\beta}}\to 0$ \
  as \ $n\to\infty$ \ imply that
 \begin{align}\label{help_cond_needed1}
  \frac{1}{z} \frac{n|\OO(1)|}{N_n^{\frac{1}{1+\beta}}}
    \leq  N_n \frac{n|\OO(1)|}{N_n^{\frac{1}{1+\beta}}}
    = \frac{n|\OO(1)|}{N_n^{-\frac{\beta}{1+\beta}}}
    = |\OO(1)|.
 \end{align}
Hence, using \eqref{exp:1}, we obtain for large enough \ $n$
 \begin{align*}
   &\int_{N_n^{-1}}^1
     \left|1 - \ee^{\lambda z N_n^{\frac{1}{1+\beta}}
                    A_n\bigl(1 - z^{-1} N_n^{-\frac{1}{1+\beta}}\bigr)}\right|
     z^{-(2+\beta)} \, \dd z \\
   &\leq \int_{N_n^{-1}}^1
          \left|\lambda z N_n^{\frac{1}{1+\beta}}
                A_n\bigl(1 - z^{-1} N_n^{-\frac{1}{1+\beta}}\bigr)\right|
          \ee^{\left|\lambda z N_n^{\frac{1}{1+\beta}}
                     A_n\bigl(1 - z^{-1} N_n^{-\frac{1}{1+\beta}}\bigr)\right|}
          z^{-(2+\beta)} \, \dd z \\
   &\leq |\OO(1)| \ee^{|\OO(1)|}
         \int_0^1 z^{-(1+\beta)} \, \dd z
     < \infty ,
 \end{align*}
 which, together with \eqref{help10A}, imply \eqref{limsup_an}.

Now we turn to prove \eqref{lim_an}.
By \eqref{exp:4}, we have
 \begin{equation*}
  \begin{split}
   &\left|\int_0^{N_n^{-1}}
           {\Bigl(1 - \ee^{-\frac{\lambda z}{2}
                          (\sum_{\ell=1}^m \theta_\ell (t_\ell - t_{\ell-1}))^2}\Bigr)}
           z^{-(2+\beta)} \, \dd z\right|
    \leq \int_0^{N_n^{-1}}
          {\frac{\lambda z(\sum_{\ell=1}^m \theta_\ell (t_\ell - t_{\ell-1}))^2}{2}}
          z^{-(2+\beta)} \, \dd z \\
   &= \frac{\lambda(\sum_{\ell=1}^m \theta_\ell (t_\ell - t_{\ell-1}))^2}{2}
            \int_0^{N_n^{-1}} z^{-(1+\beta)}
             \, \dd z
          = \frac{\lambda(\sum_{\ell=1}^m \theta_\ell (t_\ell - t_{\ell-1}))^2}{2}
            \frac{N_n^\beta}{(-\beta)}
    \to 0
  \end{split}
 \end{equation*}
 as \ $n \to \infty$, \ hence \eqref{lim_an} reduces to check that
 \ $\lim_{n\to\infty} I_n = 0$, \ where
 \begin{align*}
  I_n
  := \int^\infty_{N_n^{-1}}
      \Bigl[\ee^{\lambda z N_n^{\frac{1}{1+\beta}}
                 A_n(1 - z^{-1} N_n^{-\frac{1}{1+\beta}})}
            - \ee^{-\frac{\lambda z}{2}
                    (\sum_{\ell=1}^m \theta_\ell (t_\ell - t_{\ell-1}))^2}\Bigr]
      z^{-(2+\beta)} \, \dd z .
 \end{align*}
Applying again \eqref{an}, we obtain
 \begin{align*}
   |I_n|
    & \leq \int^\infty_{N_n^{-1}}
         \ee^{-\frac{\lambda z}{2n^2} (\sum_{\ell=1}^m \theta_\ell (\lfloor nt_\ell\rfloor - \lfloor nt_{\ell-1}\rfloor) )^2   }
          \Bigl|\ee^{z N_n^{-\frac{1}{2(1+\beta)}} \OO(1)
                     + n N_n^{-\frac{1}{1+\beta}} \OO(1)}
                - 1\Bigr|
          z^{-(2+\beta)} \, \dd z \\
    & \phantom{\leq} + \int^\infty_{N_n^{-1}}
        \Bigl| \ee^{-\frac{\lambda z}{2n^2} (\sum_{\ell=1}^m \theta_\ell (\lfloor nt_\ell\rfloor - \lfloor nt_{\ell-1}\rfloor) )^2 }
                - \ee^{-\frac{\lambda z}{2} (\sum_{\ell=1}^m \theta_\ell (t_\ell - t_{\ell-1}) )^2 } \Bigr| z^{-(2+\beta)} \,\dd z \\
    & =: I_n^{(1)} + I_n^{(2)}.
 \end{align*}
Here, for \ $z \in (N_n^{-1}, \infty)$, \ we have
 \[
   \bigl|z N_n^{-\frac{1}{2(1+\beta)}} \OO(1) + n N_n^{-\frac{1}{1+\beta}} \OO(1)\bigr|
   \leq z \bigl(N_n^{-\frac{1}{2(1+\beta)}} + n N_n^{\frac{\beta}{1+\beta}}\bigr) |\OO(1)| ,
 \]
 and hence, by \eqref{exp:1}, we get
 \begin{align*}
  &\Bigl|\ee^{z N_n^{-\frac{1}{2(1+\beta)}}\OO(1)+nN_n^{-\frac{1}{1+\beta}}\OO(1)}-1\Bigr| \\
  &\leq \bigl|z N_n^{-\frac{1}{2(1+\beta)}} \OO(1) + nN_n^{-\frac{1}{1+\beta}}  \OO(1)\bigr|\,
        \ee^{\bigl|zN_n^{-\frac{1}{2(1+\beta)}}\OO(1)+nN_n^{-\frac{1}{1+\beta}}\OO(1)\bigr|}\\
  &\leq z \bigl(N_n^{-\frac{1}{2(1+\beta)}} + n N_n^{\frac{\beta}{1+\beta}}\bigr) |\OO(1)| \,
     \ee^{z \bigl(N_n^{-\frac{1}{2(1+\beta)}} + n N_n^{\frac{\beta}{1+\beta}}\bigr) |\OO(1)|} .
 \end{align*}
Consequently, for large enough \ $n$,
 \begin{align*}
  I^{(1)}_n
  &\leq \bigl(N_n^{-\frac{1}{2(1+\beta)}} + n N_n^{\frac{\beta}{1+\beta}}\bigr) |\OO(1)|
        \int^\infty_{N_n^{-1}}
         \ee^{-\frac{\lambda z}{2n^2}(\sum_{\ell=1}^m \theta_\ell (\lfloor nt_\ell\rfloor - \lfloor nt_{\ell-1}\rfloor) )^2
              +z \bigl(N_n^{-\frac{1}{2(1+\beta)}} + n N_n^{\frac{\beta}{1+\beta}}\bigr) |\OO(1)|}
         z^{-(1+\beta)} \, \dd z \\
  &\leq \bigl(N_n^{-\frac{1}{2(1+\beta)}} + n N_n^{\frac{\beta}{1+\beta}}\bigr) |\OO(1)|
        \int_0^\infty
         \ee^{-\frac{\lambda z}{4}(\sum_{\ell=1}^m \theta_\ell (t_\ell - t_{\ell-1}))^2}
         z^{-(1+\beta)} \, \dd z ,
 \end{align*}
 that gets arbitrarily close to zero as \ $n$ \ approaches infinity, since the integral is
 finite due to the fact that
 \[
   \frac{1}{\Gamma(-\beta)}
   \left(\frac{\lambda}{4}
         \left(\sum_{\ell=1}^m \theta_\ell (t_\ell - t_{\ell-1})\right)^2\right)^{-\beta}
   \ee^{-\lambda z {(\sum_{\ell=1}^m \theta_\ell (t_\ell - t_{\ell-1}))^2/4}} \, z^{-(1+\beta)},
   \qquad z > 0,
 \]
 is the density function of a Gamma distributed random variable with
 parameters \ $-\beta$ \ and
 \ $\lambda(\sum_{\ell=1}^m \theta_\ell (t_\ell - t_{\ell-1}))^2/4$.
\ Further, by \eqref{exp:4},
 \begin{align*}
  I^{(2)}_n
    & = \int^\infty_{N_n^{-1}}
         \ee^{-\frac{\lambda z}{2} ( \sum_{\ell=1}^m \theta_\ell (t_\ell - t_{\ell-1}) )^2 }
          \Bigl| \ee^{-\frac{\lambda z}{2n^2} (\sum_{\ell=1}^m \theta_\ell (\lfloor nt_\ell\rfloor - \lfloor nt_{\ell-1}\rfloor) )^2
                       + \frac{\lambda z}{2} ( \sum_{\ell=1}^m \theta_\ell (t_\ell - t_{\ell-1}) )^2 }
                   - 1 \Bigr| z^{-(2+\beta)} \,\dd z\\
    & \leq \frac{\lambda}{2}
         \left(
              \frac{1}{n^2} \left(\sum_{\ell=1}^m \theta_\ell (\lfloor nt_\ell\rfloor - \lfloor nt_{\ell-1}\rfloor) \right)^2
               - \left( \sum_{\ell=1}^m \theta_\ell (t_\ell - t_{\ell-1}) \right)^2
            \right)\\
    &\phantom{\leq}\times\int^\infty_{N_n^{-1}}
            \ee^{-\frac{\lambda z}{2} ( \sum_{\ell=1}^m \theta_\ell (t_\ell - t_{\ell-1}) )^2 }
             z^{-(1+\beta)} \,\dd z \\
    & \to 0 \qquad \text{as \ $n\to\infty$.}
 \end{align*}
This yields \eqref{lim_an} completing the proof.
\proofend

\medskip

\noindent{\bf Proof of Theorem \ref{simultaneous_aggregation_random_0}.}
To prove this limit theorem we have to show that for any sequence
 \ $(N_n)_{n\in\NN}$ \ of positive integers with \ $(\log N_n)^{2} n^{-1} \to \infty$, \ we have
 \[
   n^{-1} {(N_n\log N_n)^{-\frac{1}{2}}} \,
   S^{(N_n,n)}
   \distrf (W_{\lambda \psi_1} t)_{t\in\RR_+} \qquad
   \text{as \ $n \to \infty$.}
 \]
For this, by continuous mapping theorem, it is enough to verify that for any
 \ $m \in \NN$ \ and \ $t_0, t_1, \ldots, t_m \in \RR_+$ \ with
 \ $0 =: t_0 < t_1 < \ldots < t_m$, \ we have
 \begin{align*}
  &n^{-1} {(N_n\log N_n)^{-\frac{1}{2}}}
   \sum_{j=1}^{N_n}
    \Biggl(\sum_{k=1}^{\lfloor nt_1\rfloor}
            \Bigl(X_k^{(j)}-\frac{\lambda}{1-\alpha^{(j)}}\Bigr),
           \sum_{k=\lfloor nt_1\rfloor+1}^{\lfloor nt_2\rfloor}
            \Bigl(X_k^{(j)}-\frac{\lambda}{1-\alpha^{(j)}}\Bigr),
           \dots\\
   &\phantom{n^{-1} {(N_n\log N_n)^{-\frac{1}{2}}} \sum_{j=1}^{N_n} \Bigl( }
         \ldots,
         \sum_{k=\lfloor nt_{m-1}\rfloor+1}^{\lfloor nt_m\rfloor}
            \Bigl(X_k^{(j)}-\frac{\lambda}{1-\alpha^{(j)}}\Bigr)\Biggr) \\
  &\distr W_{\lambda \psi_1}(t_1, t_2 - t_1, \dots, t_m - t_{m-1}) \qquad
   \text{as \ $n \to \infty$.}
 \end{align*}
So, by continuity theorem, we have to check that for any \ $m \in \NN$,
 \ $t_0, t_1, \ldots, t_m \in \RR_+$ \ with \ $0 = t_0 < t_1 < \ldots < t_m$ \ and
 \ $\theta_1, \dots, \theta_m \in \RR$ \ the convergence
 \begin{equation*}
  \begin{split}
   &\EE\Biggl(\exp\Biggl\{\ii
                          \sum_{\ell=1}^m
                           \theta_\ell n^{-1} {(N_n\log N_n)^{-\frac{1}{2}}}
                           \sum_{j=1}^{N_n}
                            \sum_{k=\lfloor nt_{\ell-1}\rfloor+1}^{\lfloor nt_\ell\rfloor}
                             \Bigl(X_k^{(j)}
                                   - \frac{\lambda}
                                          {1-\alpha^{(j)}}\Bigr)\Biggr\}\Biggr) \\
   &\qquad
    =\EE\Biggl(\exp\Biggl\{\ii n^{-1} {(N_n\log N_n)^{-\frac{1}{2}}}
                           \sum_{j=1}^{N_n}
                            \sum_{\ell=1}^m
                             \theta_\ell
                             \sum_{k=\lfloor nt_{\ell-1}\rfloor+1}^{\lfloor nt_\ell\rfloor}
                              \Bigl(X_k^{(j)}
                                    - \frac{\lambda}
                                           {1-\alpha^{(j)}}\Bigr)\Biggr\}\Biggr)\\
   &\qquad
    =\Biggl[\EE\Biggl(\exp\Biggl\{\ii n^{-1} {(N_n\log N_n)^{-\frac{1}{2}}}
                                  \sum_{\ell=1}^m
                                   \theta_\ell
                                   \sum_{k=\lfloor nt_{\ell-1}\rfloor+1}^{\lfloor nt_\ell\rfloor}
                                    \Bigl(X_k
                                          - \frac{\lambda}
                                                 {1-\alpha}\Bigr)\Biggr\}\Biggr)\Biggr]^{N_n} \\
   &\qquad
    \to \EE\Biggl(\ee^{\ii
                       \sum_{\ell=1}^m
                        \theta_\ell (t_\ell - t_{\ell-1}) W_{\lambda \psi_1}}\Biggr)
    = \ee^{-\frac{\lambda \psi_1(\sum_{\ell=1}^m \theta_\ell (t_\ell - t_{\ell-1}))^{2}}{2}}
    \qquad \text{as \ $n \to \infty$}
  \end{split}
 \end{equation*}
 holds.
Note that it suffices to show
 \begin{align*}
  \Theta_n
  &:=N_n
     \Biggl[1 - \EE\Biggl(\exp\Biggl\{\ii n^{-1} {(N_n\log N_n)^{-\frac{1}{2}}}
                                    \sum_{\ell=1}^m
                                     \theta_\ell
                                      \sum_{k=\lfloor nt_{\ell-1}\rfloor+1}^{\lfloor nt_\ell\rfloor}
                                       \Bigl(X_k
                                             - \frac{\lambda}
                                                    {1-\alpha}\Bigr)\Biggr\}\Biggr)\Biggr] \\
  &\to \frac{\lambda \psi_1(\sum_{\ell=1}^m \theta_\ell (t_\ell - t_{\ell-1}))^{2}}{2}
   \qquad \text{as \ $n \to \infty$,}
 \end{align*}
 since it implies that
 \ $(1 - \Theta_n/N_n)^{N_n}
    \to \ee^{-\frac{\lambda \psi_1(\sum_{\ell=1}^m \theta_\ell (t_\ell - t_{\ell-1}))^{2}}{2}}$
 \ as \ $n \to \infty$, \ as desired.
By applying \eqref{help1_alter} (or \eqref{egyuttes_char_aggregalt1}) to the left hand side, we get
 \begin{equation*}
  \begin{split}
   \Theta_n
   &= N_n \EE\Biggl[1 -  \ee^{-\ii n^{-1} {(N_n\log N_n)^{-\frac{1}{2}}}\frac{\lambda}{1-\alpha}
            \sum_{\ell=1}^m \theta_\ell (\lfloor nt_{\ell}\rfloor - \lfloor nt_{\ell-1} \rfloor)} \\
   &\phantom{=\;}
     \times F_{0,\dots,\lfloor nt_m\rfloor-1}
      \Bigl( \underbrace{ \ee^{\ii n^{-1} {(N_n\log N_n)^{-\frac{1}{2}}}\theta_1},\ldots, \ee^{\ii n^{-1} {(N_n\log N_n)^{-\frac{1}{2}}}\theta_1} }_{\lfloor nt_1\rfloor \ \text{items}},
      \dots\\
    &\phantom{=\;\times F_{0,\dots,\lfloor nt_m\rfloor-1} \Bigl( }\ldots,
   \underbrace{ \ee^{\ii n^{-1} {(N_n\log N_n)^{-\frac{1}{2}}}\theta_m}, \ldots, \ee^{\ii n^{-1} {(N_n\log N_n)^{-\frac{1}{2}}}\theta_m} }_{\lfloor nt_m\rfloor - \lfloor nt_{m-1}\rfloor \ \text{items}}  \,\Big|\, \alpha\Bigr)
         \Biggr]\\
   &= N_n \EE\left[1 - \ee^{\frac{\lambda}{1-\alpha}B_n(\alpha)}\right]
    = N_n \int_0^1
         \left(1 - \ee^{\frac{\lambda}{1-a}B_n(a)}\right) \psi(a)  \, \dd a
  \end{split}
 \end{equation*}
 with
 \begin{align*}
  &B_n(a)
   := \sum_{\ell=1}^m
                \Bigl(\ee^{\ii n^{-1} {(N_n\log N_n)^{-\frac{1}{2}}}\theta_\ell} - 1 - \ii n^{-1} {(N_n\log N_n)^{-\frac{1}{2}}}\theta_\ell\Bigr)
                (\lfloor nt_\ell\rfloor - \lfloor nt_{\ell-1}\rfloor)  \\
   &+\sum_{1\leq\ell_1<\ell_2\leq m}
      \sum_{k_1=\lfloor nt_{\ell_1-1}\rfloor+1}^{\lfloor nt_{\ell_1}\rfloor}
       \sum_{k_2=\lfloor nt_{\ell_2-1}\rfloor+1}^{\lfloor nt_{\ell_2}\rfloor}
        a^{k_2-k_1}
        \bigl(\ee^{\ii n^{-1} {(N_n\log N_n)^{-\frac{1}{2}}}\theta_{\ell_1}} - 1\bigr)
        \bigl(\ee^{\ii n^{-1} {(N_n\log N_n)^{-\frac{1}{2}}}\theta_{\ell_2}} - 1\bigr) \\
   &\phantom{+
               \sum_{k_1=\lfloor nt_{\ell_1-1}\rfloor+1}^{\lfloor nt_{\ell_1}\rfloor}
                \sum_{k_2=\lfloor nt_{\ell_2-1}\rfloor+1}^{\lfloor nt_{\ell_2}\rfloor}}
        \times
        \ee^{\ii n^{-1} {(N_n\log N_n)^{-\frac{1}{2}}}
             \bigl((\lfloor nt_{\ell_1}\rfloor-k_1)\theta_{\ell_1}
                   +\sum_{\ell=\ell_1+1}^{\ell_2-1}
                     \theta_\ell (\lfloor nt_\ell\rfloor - \lfloor nt_{\ell-1}\rfloor)
                   +(k_2 - 1 - \lfloor nt_{\ell_2-1}\rfloor)\theta_{\ell_2}\bigr)} \\
   &+\sum_{\ell=1}^m
      \sum_{\lfloor nt_{\ell-1}\rfloor+1\leq k_1<k_2\leq\lfloor nt_\ell\rfloor}
       a^{k_2-k_1}
      \bigl(\ee^{\ii n^{-1} {(N_n\log N_n)^{-\frac{1}{2}}}\theta_\ell} - 1\bigr)^2
        \ee^{\ii n^{-1} {(N_n\log N_n)^{-\frac{1}{2}}} (k_2-k_1-1)\theta_\ell}
 \end{align*}
 for \ $a \in [0, 1]$.
\ The aim of the following discussion is to apply Lemma \ref{ordo} with \ $z_n(a):=B_n(a)$, \ $n\in\NN$, \ $a\in(0,1)$, \
 \ $\vare_n := {(\log N_n)}^{-1}$, \ $n \in \NN$, \ and \ $I:=\frac{ \lambda}{2}  (\sum_{\ell=1}^m \theta_\ell (t_\ell - t_{\ell-1}))^{2}$.
\ Note that \ $\vare_n\in(0,1)$ \ for \ $n\geq n_0$, \ where \ $n_0$ \ is sufficiently large, and \ $\lim_{n\to\infty}\vare_n=0$.
\ First we check \eqref{help_ordo1}.
Using \eqref{exp:3}, for any \ $a \in (0, 1)$ \ we get
 \begin{align*}
   |B_n(a)|
   &\leq \sum_{\ell=1}^m
             n^{-2} {(N_n\log N_n)^{-1}} \frac{\theta_\ell^2}{2}
             (\lfloor nt_\ell\rfloor - \lfloor nt_{\ell-1}\rfloor) \\
   &\quad
    + \sum_{1\leq\ell_1<\ell_2\leq m}
         n^{-2} {(N_n\log N_n)^{-1}} |\theta_{\ell_1}| |\theta_{\ell_2}|
         (\lfloor nt_{\ell_1}\rfloor - \lfloor nt_{\ell_1-1}\rfloor)
         (\lfloor nt_{\ell_2}\rfloor - \lfloor nt_{\ell_2-1}\rfloor) \\
   &\quad
    + \sum_{\ell=1}^m
         n^{-2} {(N_n\log N_n)^{-1}} \frac{\theta_\ell^2}{2}
         (\lfloor nt_\ell\rfloor - \lfloor nt_{\ell-1}\rfloor)
         (\lfloor nt_\ell\rfloor - \lfloor nt_{\ell-1}\rfloor - 1)\\
  & = \frac{1}{2} n^{-2} {(N_n\log N_n)^{-1}}
      \biggl(\sum_{\ell=1}^m
              |\theta_\ell| (\lfloor nt_\ell\rfloor - \lfloor nt_{\ell-1}\rfloor)\biggr)^2 \\
  &\leq \frac{1}{2} {(N_n\log N_n)^{-1}}
         \biggl(\sum_{\ell=1}^m |\theta_\ell| (t_\ell - t_{\ell-1} + 1)\biggr)^2 ,
 \end{align*}
 since \ $\frac{1}{n}(\lfloor nt_\ell\rfloor - \lfloor nt_{\ell-1}\rfloor) \leq \frac{1}{n}(nt_\ell - nt_{\ell-1} + 1) = t_\ell - t_{\ell-1} + \frac{1}{n}
  \leq t_\ell - t_{\ell-1} +1$.
\ Consequently, since \ $\vare_n= (\log N_n)^{-1}$, \ we have
 \[
   \sup_{n\geq n_0} \vare_n^{-1} N_n \sup_{a\in(0,1-\vare_n)} |B_n(a)|
   \leq \frac{1}{2} \biggl(\sum_{\ell=1}^m |\theta_\ell| (t_\ell - t_{\ell-1} + 1)\biggr)^2
   < \infty ,
 \]
 i.e., \eqref{help_ordo1} is satisfied.
Therefore, by Lemma \ref{ordo}, substituting \ $a = 1 - z N_n^{-1}$ \ with \ $z>0$,
 \ the statement of the theorem will follow from
 \begin{gather}\label{limsup_bn}
  \begin{aligned}
   &\limsup_{n\to\infty}
     N_n \int^1_{1-(\log N_n)^{-1}}
        \Bigl|1 - \ee^{\frac{\lambda}{1-a}B_n(a)}\Bigr|\, \dd a \\
   &= \limsup_{n\to\infty}
       \int_0^{\frac{N_n} {\log N_n}}
        \Bigl|1 - \ee^{\lambda \frac{N_n}{z}
                       B_n\bigl(1 - z{N_n}^{-1}\bigr)}\Bigr| \, \dd z
    < \infty
  \end{aligned}
 \end{gather}
 and
 \begin{gather}\label{lim_bn}
  \begin{aligned}
   &\lim_{n\to\infty}
     \left|N_n \int_{1-(\log N_n)^{-1}}^1
                \Bigl(1 - \ee^{\frac{\lambda}{1-a}B_n(a)}\Bigr)\, \dd a
           - I\right| \\
   &= \lim_{n\to\infty}
        \left|\int_0^{\frac{N_n} {\log N_n}}
               \Bigl(1 - \ee^{\lambda \frac{N_n}{z}
                       B_n\bigl(1 - z{N_n}^{-1}\bigr)}\Bigr)
              \, \dd z
             - I\right|
    = 0
  \end{aligned}
 \end{gather}
 with \ $I = \frac{\lambda}{2}  (\sum_{\ell=1}^m \theta_\ell (t_\ell - t_{\ell-1}))^{2}$.

Next we check \eqref{limsup_bn} and \eqref{lim_bn}.
By Taylor expansion,
 \begin{gather*}
\begin{split}
&\ee^{\ii n^{-1} {(N_n\log N_n)^{-\frac{1}{2}} \theta_\ell}} - 1
  = \ii n^{-1} {(N_n\log N_n)^{-\frac{1}{2}}} \theta_\ell + n^{-2} {(N_n\log N_n)^{-1}} \OO(1)
  = n^{-1} {(N_n\log N_n)^{-\frac{1}{2}}} \OO(1), \\
&\ee^{\ii n^{-1} {(N_n\log N_n)^{-\frac{1}{2}} \theta_\ell}} - 1
  - \ii n^{-1} {(N_n\log N_n)^{-\frac{1}{2}}}\theta_{\ell}
  = - n^{-2} {(N_n\log N_n)^{-1}} \frac{\theta_{\ell}^2}{2}
    + n^{-3} {(N_n\log N_n)^{-\frac{3}{2}}} \OO(1)\\
  & \phantom{\ee^{\ii n^{-1} {(N_n\log N_n)^{-\frac{1}{2}} \theta_\ell}} - 1
  - \ii n^{-1} {(N_n\log N_n)^{-\frac{1}{2}}}\theta_{\ell}}
 =  n^{-2} {(N_n\log N_n)^{-1}} \OO(1)
\end{split}
 \end{gather*}
 for all \ $\ell \in\{1, \ldots, m\}$, \ resulting
 \begin{equation}\label{bn}
  \lambda \frac{N_n}{z}B_n\biggl(1 - \frac{z}{N_n}\biggr)
  = - \frac{\lambda
            \bigl(\sum_{\ell=1}^m
                   \theta_\ell (\lfloor nt_\ell\rfloor - \lfloor nt_{\ell-1}\rfloor)\bigr)^2}
           {2 z n^2 \log N_n}
    + \frac{\OO(1)}{z N_n^{\frac{1}{2}}(\log N_n)^{\frac{3}{2}}}
    + \frac{n\OO(1)}{N_n \log N_n}
 \end{equation}
 for \ $z <  N_n$.
\ Indeed, \ $z < N_n$ \  yields that \ $1 - z/N_n\in(0,1)$, \ and
 \begin{align*}
   &B_n\biggl(1 - \frac{z}{N_n}\biggr) \\
   &= \sum_{\ell=1}^m
       \bigl(\ee^{\ii n^{-1} {(N_n\log N_n)^{-\frac{1}{2}}}\theta_\ell} - 1
             - \ii n^{-1} {(N_n\log N_n)^{-\frac{1}{2}}}\theta_\ell\bigr)
       (\lfloor nt_\ell\rfloor - \lfloor nt_{\ell-1}\rfloor) \\
   &\quad
      +\sum_{1\leq\ell_1<\ell_2\leq m}
        \sum_{k_1=\lfloor nt_{\ell_1-1}\rfloor+1}^{\lfloor nt_{\ell_1}\rfloor}
         \sum_{k_2=\lfloor nt_{\ell_2-1}\rfloor+1}^{\lfloor nt_{\ell_2}\rfloor}
          \left(1 - \frac{z}{N_n}\right)^{k_2-k_1}
          \bigl(\ee^{\ii n^{-1}(N_n\log N_n)^{-\frac{1}{2}}\theta_{\ell_1}} - 1\bigr) \\
   &\phantom{+}\times
        \ee^{\ii n^{-1} {(N_n\log N_n)^{-\frac{1}{2}}}
             \bigl((\lfloor nt_{\ell_1}\rfloor-k_1)\theta_{\ell_1}
                   + \sum_{\ell=\ell_1+1}^{\ell_2-1}
                     \theta_\ell (\lfloor nt_\ell\rfloor - \lfloor nt_{\ell-1}\rfloor)
                   +(k_2- 1 - \lfloor nt_{\ell_2-1}\rfloor)\theta_{\ell_2}\bigr)}
        \bigl(\ee^{\ii n^{-1} {(N_n\log N_n)^{-\frac{1}{2}}}\theta_{\ell_2}} - 1\bigr) \\
   &+\sum_{\ell=1}^m
      \sum_{\lfloor nt_{\ell-1}\rfloor+1\leq k_1<k_2\leq\lfloor nt_\ell\rfloor}
       \left(1 - \frac{z}{N_n}\right)^{k_2-k_1}
       \bigl(\ee^{\ii n^{-1} {(N_n\log N_n)^{-\frac{1}{2}}}\theta_\ell} - 1\bigr)^2
        \ee^{\ii n^{-1} {(N_n\log N_n)^{-\frac{1}{2}}} (k_2-k_1-1)\theta_\ell} \\
   &= \sum_{\ell=1}^m
       \biggl(- \frac{\theta_\ell^2}{2n^{2} {N_n\log N_n}}
              + \frac{\OO(1)}{n^{3} {(N_n\log N_n)^{\frac{3}{2}}}}\biggr)
       (\lfloor nt_\ell\rfloor - \lfloor nt_{\ell-1}\rfloor) \\
   &\quad
      +\sum_{1\leq\ell_1<\ell_2\leq m}
        \left(1 + \frac{n\, z\OO(1)}{N_n}\right)
        \biggl(\frac{\ii\theta_{\ell_1}}{n(N_n\log N_n)^{\frac{1}{2}}}
               + \frac{\OO(1)}{n^{2} {N_n\log N_n}}\biggr)
        \biggl(1 + \frac{\OO(1)}{(N_n\log N_n)^{\frac{1}{2}}}\biggr)
         \\
   &\phantom{\quad+\sum_{1\leq\ell_1<\ell_2\leq m}} \times
	\biggl(\frac{\ii\theta_{\ell_2}}{n(N_n\log N_n)^{\frac{1}{2}}}
               + \frac{\OO(1)}{n^2N_n\log N_n}\biggr)
        (\lfloor nt_{\ell_1}\rfloor - \lfloor nt_{\ell_1-1}\rfloor)
        (\lfloor nt_{\ell_2}\rfloor - \lfloor nt_{\ell_2-1}\rfloor) \\
   &\quad
      + \frac{1}{2}\sum_{\ell=1}^m
        \left(1 + \frac{n \, z\OO(1)}{N_n}\right)
        \biggl(\frac{\ii\theta_\ell}{n(N_n\log N_n)^{\frac{1}{2}}}
               + \frac{\OO(1)}{n^2N_n\log N_n}\biggr)^2
        \biggl(1 + \frac{\OO(1)}{(N_n\log N_n)^{\frac{1}{2}}}\biggr) \\
   &\phantom{=+\sum_{\ell=1}^m}\times
        (\lfloor nt_\ell\rfloor - \lfloor nt_{\ell-1}\rfloor)
        (\lfloor nt_\ell\rfloor - \lfloor nt_{\ell-1} - 1\rfloor)\\
   &= - \frac{\sum_{\ell=1}^m \theta_\ell^2 (\lfloor nt_\ell\rfloor - \lfloor nt_{\ell-1}\rfloor)}
             {2n^2N_n\log N_n}
      + \frac{\OO(1)}{n^2(N_n\log N_n)^{\frac{3}{2}}} \\
   &\quad
      - \frac{\sum_{1\leq\ell_1<\ell_2\leq m}
               \theta_{\ell_1} \theta_{\ell_2}
               (\lfloor nt_{\ell_1}\rfloor - \lfloor nt_{\ell_1-1}\rfloor)
               (\lfloor nt_{\ell_2}\rfloor - \lfloor nt_{\ell_2-1}\rfloor)}
             {n^2N_n\log N_n}
      + \frac{\OO(1)}{(N_n\log N_n)^{\frac{3}{2}}}
      + \frac{n \, z\OO(1)}{N_n^2\log N_n} \\
   &\quad
      - \frac{\sum_{\ell=1}^m
               \theta_\ell^2 (\lfloor nt_\ell\rfloor - \lfloor nt_{\ell-1}\rfloor)
               (\lfloor nt_\ell\rfloor - \lfloor nt_{\ell-1} - 1\rfloor)}
             {2n^2N_n\log N_n}
      + \frac{\OO(1)}{(N_n\log N_n)^{\frac{3}{2}}}
      + \frac{n \, z\OO(1)}{N_n^2\log N_n}
 \end{align*}
 \begin{align*}
   &= - \frac{\bigl(\sum_{\ell=1}^m
                     \theta_\ell (\lfloor nt_\ell\rfloor - \lfloor nt_{\ell-1}\rfloor)\bigr)^2}
             {2n^2N_n\log N_n}
      + \frac{\OO(1)}{(N_n\log N_n)^{\frac{3}{2}}}
      + \frac{n \, z\OO(1)}{N_n^2\log N_n} ,
 \end{align*}
where we used the corresponding versions of \eqref{help_ordo3} and \eqref{help_ordo4} replacing \ $N_n^{-\frac{1}{2(1+\beta)}}$ \
 by \ $(N_n\log N_n)^{-\frac{1}{2}}$ \ and that
   \[
     \left(1 - \frac{z}{N_n}\right)^{k_2-k_1}
      = 1 + \frac{n \, z\OO(1)}{N_n}
    \]
following from Bernoulli's inequality.
By \eqref{bn}, for large enough \ $n$ \ and for any \ $z\in(0,N_n)$, \ we have
 \begin{align*}
  &\lambda \frac{N_n}{z}\Re B_n\biggl(1 - \frac{z}{N_n}\biggr)
  = - \frac{\lambda
            \bigl(\sum_{\ell=1}^m
                   \theta_\ell (\lfloor nt_\ell\rfloor - \lfloor nt_{\ell-1}\rfloor)\bigr)^2}
           {2zn^2 \log N_n} \left(1 - \frac{\Re\OO(1)}{(N_n \log N_n)^{\frac{1}{2}}}\right)
    + \frac{n\Re\OO(1)}{N_n \log N_n} \\
  &\leq - \frac{\lambda \big(\sum_{\ell=1}^m \theta_\ell (\lfloor nt_\ell\rfloor - \lfloor nt_{\ell-1}\rfloor)\big)^2}{4 z n^2 \log N_n}
        + \frac{n|\OO(1)|}{N_n\log N_n},
 \end{align*}
hence we obtain that
 \begin{align}\label{help10B}
  \begin{split}
   &\int_0^{(\log N_n)^{-1}}
     \left|1 - \ee^{\lambda \frac{N_n}{z}B_n\left(1 - \frac{z}{N_n}\right)}\right|
   \, \dd z
   \leq \int_0^{(\log N_n)^{-1}}
          \left(1 + \ee^{\lambda \frac{N_n}{z}\Re B_n\left(1 - \frac{z}{N_n}\right)}\right)
 \, \dd z \\
& \leq {(\log N_n)^{-1}}
\left(
1+\exp{\left\{ - \frac{\lambda \big(\sum_{\ell=1}^m \theta_\ell ( \lfloor nt_\ell\rfloor - \lfloor nt_{\ell-1}\rfloor )  \big)^2}{4 z n^2 \log N_n}
        + \frac{n|\OO(1)|}{N_n\log N_n} \right\}}
\right)\to 0
  \end{split}
 \end{align}
as \ $n\to \infty$, \ since
 \[
   \lim_{n\to\infty} \frac{\big(\sum_{\ell=1}^m \theta_\ell (  \lfloor nt_\ell\rfloor - \lfloor nt_{\ell-1}\rfloor   )  \big)^2}
                           {n^2 \log N_n}=0,
 \]
 and, due to the assumption \ $(\log N_n)^2 n^{-1}\to \infty$ \ as \ $n\to \infty$, \ we have
 \[
 \frac{n}{N_n \log N_n} = \frac{n}{(\log N_n)^2}\frac{\log N_n}{N_n}  \to 0 \qquad \text{as \ $n\to\infty$.}
 \]
Note that for every \ $z\in ((\log N_n)^{-1}, N_n (\log N_n)^{-1})$ \ we have
 \begin{align}\label{help10AB}
  \begin{split}
& \left| \lambda \frac{N_n}{z}B_n\left(1 - \frac{z}{N_n}\right) \right|
  \leq
  \frac{\lambda\bigl(\sum_{\ell=1}^m
                     \vert\theta_\ell\vert (t_\ell - t_{\ell-1} +1 )\bigr)^2}
             {2 z \log N_n}
      + \frac{|\OO(1)|}{z N_n^{\frac{1}{2}}(\log N_n)^{\frac{3}{2}}}
      + \frac{n \, |\OO(1)|}{N_n\log N_n}\\
& \leq
\frac{ \lambda\bigl(\sum_{\ell=1}^m
                      \vert \theta_\ell\vert  (t_\ell - t_{\ell-1}+1 )\bigr)^2}{2}
      + \frac{|\OO(1)|}{ N_n^{\frac{1}{2}}(\log N_n)^{\frac{1}{2}}}
      + \frac{n \, |\OO(1)|}{N_n\log N_n}= \vert\OO(1)\vert,
  \end{split}
 \end{align}
since \ $n (N_n\log N_n)^{-1}\to 0$ \ as \ $n\to \infty$, \ as we have seen before.

Hence, using \eqref{exp:1}, we obtain for large enough \ $n$
 \begin{align*}
   &\int_{(\log N_n)^{-1}}^{N_n(\log N_n)^{-1}}
     \left|1 - \ee^{\lambda \frac{N_n}{z}B_n\left(1 - \frac{z}{N_n}\right)}\right|
      \, \dd z \leq
 \end{align*}
 \begin{align*}
   &\leq \int_{(\log N_n)^{-1}}^{N_n(\log N_n)^{-1}}
          \left|\lambda \frac{N_n}{z}B_n\left(1 - \frac{z}{N_n}\right)\right|
          \ee^{\left|\lambda \frac{N_n}{z}B_n\left(1 - \frac{z}{N_n}\right)\right|}
           \, \dd z \\
   &\leq\ee^{|\OO(1)|}
         \int_{(\log N_n)^{-1}}^{N_n(\log N_n)^{-1}}
\left[
 \frac{\lambda \bigl(\sum_{\ell=1}^m
                     \vert \theta_\ell\vert (t_\ell - t_{\ell-1} + 1)\bigr)^2}
             {2 z \log N_n}
      + \frac{|\OO(1)|}{z N_n^{\frac{1}{2}}(\log N_n)^{\frac{3}{2}}}
      + \frac{n \, |\OO(1)|}{N_n\log N_n}
\right] \, \dd z
     < \infty ,
 \end{align*}
 since for every \ $N_n \in \NN$, \ we have
 \begin{align}\label{help_int1}
 \frac{1}{\log N_n} \int_{(\log N_n)^{-1}}^{N_n(\log N_n)^{-1}} \frac{1}{z}  \, \dd z=1,
 \end{align}
 and
 \begin{align}\label{help_int2}
   \frac{n}{N_n \log N_n}
      \int_{(\log N_n)^{-1}}^{N_n(\log N_n)^{-1}} 1 \, \dd z
    = \frac{n(N_n -1)}{N_n (\log N_n)^2}
    = \frac{n}{(\log N_n)^2} \left(1 - \frac{1}{N_n}\right)
    \to 0
 \end{align}
 as \ $n\to\infty$ \ due to the assumption \ $n^{-1}(\log N_n)^2\to\infty$ \ as \ $n\to\infty$.
 Together with \eqref{help10B}, this implies \eqref{limsup_bn}.

Now we turn to prove \eqref{lim_bn}.
By \eqref{help10B}, the convergence \eqref{lim_bn} reduces to showing that
 \[
   \left\vert
       \int_{(\log N_n)^{-1}}^{N_n(\log N_n)^{-1}} \Big(1 - \ee^{\lambda \frac{N_n}{z}B_n\left(1 - \frac{z}{N_n}\right)} \Big)\, \dd z
            - \frac{\lambda(\sum_{\ell=1}^m \theta_\ell (t_\ell - t_{\ell-1}))^2}{2}
      \right\vert \to 0
 \]
 as \ $n\to\infty$.
\ Using \eqref{help_int1}, it is enough to check that
 \[
   \left|\int_{(\log N_n)^{-1}}^{N_n(\log N_n)^{-1}}
           \left(
\ee^{\lambda \frac{N_n}{z}B_n\left(1 - \frac{z}{N_n}\right)} -1 + \frac{\lambda(\sum_{\ell=1}^m \theta_\ell (t_\ell - t_{\ell-1}))^2}{2 z \log N_n}
\right) \, \dd z\right| \to 0
 \]
 as \ $n\to\infty$.
\ By applying \eqref{exp:2}, \eqref{bn} and \eqref{help10AB}, we have
 \begin{align*}
 &\left|\int_{(\log N_n)^{-1}}^{N_n(\log N_n)^{-1}}
   \left(
    \ee^{\lambda \frac{N_n}{z}B_n\left(1 - \frac{z}{N_n}\right)} -1 + \frac{\lambda(\sum_{\ell=1}^m \theta_\ell (t_\ell - t_{\ell-1}))^2}{2 z \log N_n}
   \right) \, \dd z\right|\\
 & \leq \int_{(\log N_n)^{-1}}^{N_n(\log N_n)^{-1}}
 \Bigg( \frac{1}{2} \left|
   \lambda \frac{N_n}{z}B_n\left(1 - \frac{z}{N_n}\right) \right|^2
  \ee^{\left| \lambda \frac{N_n}{z}B_n\left(1 - \frac{z}{N_n}\right) \right|} \\
 & \phantom{\leq \int_{(\log N_n)^{-1}}^{N_n(\log N_n)^{-1}} \Bigg(}
  +\left| \lambda \frac{N_n}{z}B_n\left(1 - \frac{z}{N_n}\right)+\frac{\lambda(\sum_{\ell=1}^m \theta_\ell (t_\ell - t_{\ell-1}))^2}{2 z \log N_n}
   \right|\Bigg) \, \dd z\\
 &\leq\int_{(\log N_n)^{-1}}^{N_n(\log N_n)^{-1}}
 \Bigg( \frac{1}{2} \left(
\frac{\lambda(\sum_{\ell=1}^m \vert \theta_\ell\vert (t_\ell - t_{\ell-1} +1))^2}{2 z \log N_n}
+\frac{|\OO(1)|}{zN_n^{\frac{1}{2}}(\log N_n)^\frac{3}{2}}+\frac{n|\OO(1)|}{N_n \log N_n}
\right)^2
\ee^{|\OO(1)|}  \\
 & \phantom{\leq\int_{(\log N_n)^{-1}}^{N_n(\log N_n)^{-1}} \Bigg(}
  + \frac{|\OO(1)|}{z N_n^{\frac{1}{2}}(\log N_n)^\frac{3}{2}}+\frac{n|\OO(1)|}{N_n \log N_n}\\
 & \phantom{\leq\int_{(\log N_n)^{-1}}^{N_n(\log N_n)^{-1}} \Bigg(}
  + \frac{\lambda }{2z\log N_n}\left\vert  \Big(\sum_{\ell=1}^m \theta_\ell (t_\ell - t_{\ell-1}) \Big)^2
                                             -  n^{-2}\Big(\sum_{\ell=1}^m \theta_\ell (\lfloor nt_\ell\rfloor - \lfloor nt_{\ell-1}\rfloor) \Big)^2  \right\vert
  \Bigg) \, \dd z
 \end{align*}
 \begin{align*}
 &\leq\int_{(\log N_n)^{-1}}^{N_n(\log N_n)^{-1}}
 \Bigg( \frac{3}{2} \left(
\frac{\lambda^2(\sum_{\ell=1}^m {\vert \theta_\ell\vert} (t_\ell - t_{\ell-1}+1))^4}{4 z^2 (\log N_n)^2}
+\frac{|\OO(1)|}{z^2N_n(\log N_n)^3}+\frac{n^2|\OO(1)|}{N_n^2 (\log N_n)^2}
\right)
\ee^{|\OO(1)|} \\
 & \phantom{\leq\int_{(\log N_n)^{-1}}^{N_n(\log N_n)^{-1}} \Bigg(}
  + \frac{|\OO(1)|}{z N_n^{\frac{1}{2}}(\log N_n)^\frac{3}{2}}+\frac{n|\OO(1)|}{N_n \log N_n}\\
 & \phantom{\leq\int_{(\log N_n)^{-1}}^{N_n(\log N_n)^{-1}} \Bigg(}
   + \frac{\lambda }{2z\log N_n}\left\vert  \Big(\sum_{\ell=1}^m \theta_\ell (t_\ell - t_{\ell-1}) \Big)^2
                                             -  n^{-2}\Big(\sum_{\ell=1}^m \theta_\ell (\lfloor nt_\ell\rfloor - \lfloor nt_{\ell-1}\rfloor) \Big)^2  \right\vert
  \Bigg) \, \dd z.
 \end{align*}
Indeed, the last but one inequality follows from
 \begin{align*}
 &\left|
    \lambda \frac{N_n}{z}B_n\left(1 - \frac{z}{N_n}\right)+\frac{\lambda(\sum_{\ell=1}^m \theta_\ell (t_\ell - t_{\ell-1}))^2}{2 z \log N_n}
  \right|\\
  &\leq \left|
     \lambda \frac{N_n}{z}B_n\left(1 - \frac{z}{N_n}\right)
     +\frac{\lambda\Big(\sum_{\ell=1}^m \theta_\ell (\lfloor nt_\ell\rfloor - \lfloor nt_{\ell-1}\rfloor) \Big)^2}{2 z n^2 \log N_n}
     \right| \\
  &\phantom{\leq} + \left|
        -\frac{\lambda\Big(\sum_{\ell=1}^m \theta_\ell (\lfloor nt_\ell\rfloor - \lfloor nt_{\ell-1}\rfloor) \Big)^2}{2 z n^2 \log N_n}
        +\frac{\lambda(\sum_{\ell=1}^m \theta_\ell (t_\ell - t_{\ell-1}))^2}{2 z \log N_n}
    \right| \\
 & \leq
     \frac{|\OO(1)|}{z N_n^{\frac{1}{2}}(\log N_n)^\frac{3}{2}}+\frac{n|\OO(1)|}{N_n \log N_n} \\
 &\phantom{\leq}
     + \frac{\lambda }{2z\log N_n}\left\vert  \Big(\sum_{\ell=1}^m \theta_\ell (t_\ell - t_{\ell-1}) \Big)^2
                                             -  n^{-2}\Big(\sum_{\ell=1}^m \theta_\ell (\lfloor nt_\ell\rfloor - \lfloor nt_{\ell-1}\rfloor) \Big)^2  \right\vert,
 \end{align*}
 and the last inequality from \ $(a+b+c)^2\leq 3(a^2 + b^2 + c^2)$, \ $a,b,c\in\RR$.
\ Consequently,
 \begin{align*}
 \left|\int_{(\log N_n)^{-1}}^{N_n(\log N_n)^{-1}}
           \left(
  \ee^{\lambda \frac{N_n}{z}B_n\left(1 - \frac{z}{N_n}\right)} -1 + \frac{\lambda(\sum_{\ell=1}^m \theta_\ell (t_\ell - t_{\ell-1}))^2}{2 z \log N_n}
 \right) \, \dd z\right|\to 0
 \qquad \text{as \ $n\to\infty$.}
 \end{align*}
Indeed,
 \begin{align*}
  \frac{1}{(\log N_n)^2}\int_{(\log N_n)^{-1}}^{N_n(\log N_n)^{-1}} \frac{1}{z^2}\,\dd z
    = \frac{1}{(\log N_n)^2} \left( \log N_n - \frac{\log N_n}{N_n} \right)
    = \frac{1}{\log N_n} - \frac{1}{N_n \log N_n}
    \to 0
 \end{align*}
 as \ $n\to\infty$, \ and hence
 \begin{align*}
  \frac{1}{N_n(\log N_n)^3}\int_{(\log N_n)^{-1}}^{N_n(\log N_n)^{-1}} \frac{1}{z^2}\,\dd z
       \to 0 \qquad \text{as \ $n\to\infty$.}
 \end{align*}
Further, using the assumption \ $(\log N_n)^2n^{-1}\to\infty$ \ as \ $n\to\infty$, \ we have
 \begin{align*}
  \frac{n^2}{N_n^2(\log N_n)^2}\int_{(\log N_n)^{-1}}^{N_n(\log N_n)^{-1}} 1\,\dd z
    & = \frac{n^2}{N_n^2(\log N_n)^2} (\log N_n)^{-1}(N_n -1)
      = \frac{n^2}{N_n(\log N_n)^3}\left(1-\frac{1}{N_n}\right)\\
    & = \left(\frac{n}{(\log N_n)^2}\right)^2\frac{\log N_n}{N_n}\left(1-\frac{1}{N_n}\right)
      \to 0 \qquad \text{as \ $n\to\infty$.}
 \end{align*}
Moreover, \eqref{help_int1} yields that
 \begin{align*}
  \frac{1}{N_n^{\frac{1}{2}}(\log N_n)^{\frac{3}{2}}} \int_{(\log N_n)^{-1}}^{N_n(\log N_n)^{-1}} \frac{1}{z}  \, \dd z
     \to 0\qquad \text{as \ $n\to\infty$,}
 \end{align*}
 and
 \begin{align*}
 &\left(\frac{1}{\log N_n}\int_{(\log N_n)^{-1}}^{N_n(\log N_n)^{-1}} \frac{1}{z}  \, \dd z \right)
   \left\vert  \Big(\sum_{\ell=1}^m \theta_\ell (t_\ell - t_{\ell-1}) \Big)^2
                                             -  n^{-2}\Big(\sum_{\ell=1}^m \theta_\ell (\lfloor nt_\ell\rfloor - \lfloor nt_{\ell-1}\rfloor) \Big)^2  \right\vert \\
 &= \left\vert  \Big(\sum_{\ell=1}^m \theta_\ell (t_\ell - t_{\ell-1}) \Big)^2
                                             -  n^{-2}\Big(\sum_{\ell=1}^m \theta_\ell (\lfloor nt_\ell\rfloor - \lfloor nt_{\ell-1}\rfloor) \Big)^2  \right\vert
 \to 0  \qquad \text{as \ $n\to\infty$.}
 \end{align*}
This together with \eqref{help_int2} yield \eqref{lim_bn}, completing the proof.
\proofend

\vspace*{5mm}

\appendix

\vspace*{5mm}

\noindent{\bf\Large Appendices}

\section{Generator function of finite-dimensional distributions of
         stationary INAR(1) processes with Poisson immigrations}
\label{prel}

Consider a strictly stationary (usual) INAR(1) process \ $(X_k)_{k\in\ZZ_+}$ \ with thinning parameter \ $a\in(0,1)$ \ and
 with Poisson immigration distribution having parameter \ $\lambda\in(0,\infty)$.
\ Namely,
 \begin{align}\label{help_stac_INAR1_1}
  \begin{split}
  & \PP(\xi_{1,1} = 1) = a = 1 - \PP(\xi_{1,1} = 0) , \\
  & \PP(\vare_1 = \ell)
     = \frac{\lambda^\ell}{\ell!} \ee^{-\lambda} , \qquad
       \ell \in \ZZ_+ ,\\
  & \PP(X_0=k) =  \frac{((1-a)^{-1} \lambda)^k }{k!}\ee^{-(1-a)^{-1} \lambda}, \qquad k\in\ZZ_+.
  \end{split}
 \end{align}
As it was recalled in Section \ref{Sec_Int_Main}, \ $(X_k)_{k\in\ZZ_+}$ \ is indeed a strictly stationary INAR(1) process.

\begin{Pro}\label{Pro_egyuttes_mom_gen}
Under the assumption \eqref{help_stac_INAR1_1}, the joint generator function of \ $(X_0, X_1, \ldots, X_k)$, \ $k \in \ZZ_+$,
 \ takes the form
 \begin{equation}\label{help1_alter}
  \begin{aligned}
   F_{0,\dots, k}(z_0, \dots, z_k)
   &:= \EE(z_0^{X_0} z_1^{X_1} \cdots z_k^{X_k}) \\
   &= \exp\left\{\frac{\lambda}{1-a}
                 \sum_{0\leq i\leq j\leq k}
                  a^{j-i} (z_i - 1) z_{i+1} \cdots z_{j-1} (z_j - 1)\right\}
  \end{aligned}
 \end{equation}
 for all \ $k \in \NN$ \ and \ $z_0, \ldots, z_k \in \CC$, \ where, for \ $i = j$,
 \ the term in the sum above is \ $z_i - 1$.
\end{Pro}

For the proof of Proposition \ref{Pro_egyuttes_mom_gen}, see the proof of Proposition 2.1
 in Barczy et al.\ \cite{BarNedPap_Arxiv} (see also Barczy et al.\ \cite[Proposition 2.1]{BarNedPap}).

\section{Infinite series representation of strictly stationary INAR(1) processes}
\label{App_Stac_series}

\begin{Lem}\label{lem:reprX}
Under the assumption \eqref{help_stac_INAR1_1}, we have
 \begin{equation}\label{repr}
  (X_k)_{k\in\ZZ}
  \distre \biggl(\vare_k + \sum_{i=1}^\infty a_k^{(k-i)} \circ \cdots \circ a_{k-i+1}^{(k-i)} \circ \vare_{k-i}\biggr)_{k\in\ZZ} ,
 \end{equation}
 where \ $\{\vare_k : k \in \ZZ\}$ \ are independent random variables with the same distribution as \ $\vare_1$ \ (given in assumption \eqref{help_stac_INAR1_1}),
 and \ $a_k^{(\ell)}$, \ $k, \ell \in \ZZ$, \ are given by
 \[
   a_k^{(\ell)} \circ i
    := \begin{cases}
       \sum_{j=1}^i \xi_{k,j}^{(\ell)} , & \text{if \ $i \in \NN$,} \\
       0 , & \text{if \ $i = 0$,}
      \end{cases}
 \]
 where \ $\xi_{k,j}^{(\ell)}$, \ $j \in \NN$, \ $k, \ell \in \ZZ$, \ have the same distribution as
 \ $\xi_{1,1}$ \ (given in assumption \eqref{help_stac_INAR1_1}), and \ $\{\vare_k : k \in \ZZ\}$ \ and \ $a_k^{(\ell)}$, \ $k, \ell \in \ZZ$,
 \ are independent in the sense that the families $\{\vare_k : k \in \ZZ\}$ and $\{\xi_{k,j}^{(\ell)} : j \in \NN\}$, $k, \ell \in \ZZ$,
 occurring in $a_k^{(\ell)}$, $k, \ell \in \ZZ$,
 are independent families of independent random variables, and the series in the representation \eqref{repr} converge with probability one.
\end{Lem}

Lemma \ref{lem:reprX} is a special case of Lemma E.2 in Barczy et al.\ \cite{BarBasKevPapPla}, where one can find a proof as well.

\section{Approximations of the exponential function and some of its integrals}
\label{App1}


In this appendix we collect some useful approximations of the exponential function and
 some of its integrals.

We will frequently use the following the well-known inequalities:
 \begin{gather}\label{exp:4}
  1 - \ee^{-x} \leq x, \qquad x \in \RR, \\  \label{exp:3}
  |\ee^{\ii u} - 1| \leq |u|, \qquad |\ee^{\ii u} - 1 - \ii u| \leq u^2/2 ,
  \qquad u \in \RR .
 \end{gather}

The next lemma is about how the inequalities in \eqref{exp:3} change if we replace \ $u\in\RR$ \ by an arbitrary complex number
(for a proof, see, e.g., the proof of Lemma B.1 in Barczy et al.\ \cite{BarNedPap_Arxiv}).

\begin{Lem}\label{exp}
For any \ $z \in \CC$ \ it holds that
 \begin{gather}\label{exp:1}
  |\ee^z - 1| \leq |z| \ee^{|z|} , \\
  \label{exp:2}
  |\ee^z - 1 - z| \leq \frac{|z|^2}{2}\ee^{|z|} .
 \end{gather}
\end{Lem}


The next lemma is a variant of Lemma B.2 in Barczy et al.\ \cite{BarNedPap} (developed for proving limit theorems for iterated aggregation
 of randomized INAR(1) processes), and we use it in the proofs of Theorems \ref{simultaneous_aggregation_random_2} and \ref{simultaneous_aggregation_random_0}.

\begin{Lem}\label{ordo}
Suppose that \ $(0, 1) \ni x \mapsto \psi(x) (1 - x)^\beta$ \ is a probability density,
 where \ $\psi$ \ is a function on \ $(0, 1)$ \ having a limit
 \ $\lim_{x\uparrow 1} \psi(x) = \psi_1 \in (0, \infty)$ \ (and then necessarily \ $\beta\in(-1,\infty)$).
For all \ $a \in (0, 1)$, \ let \ $(z_n(a))_{n\in\NN}$ \ be a sequence of complex
 numbers, let \ $n_0\in\NN$,
 \ $(\vare_n)_{n\geq n_0}$ \ be a sequence in \ $(0, 1)$ \ with
 \ $\lim_{n\to\infty} \vare_n = 0$, \ and let \ $(N_n)_{n\in\NN}$ \ be a sequence of positive
 integers such that
 \begin{gather}\label{help_ordo1}
  \sup_{ n\geq n_0} \vare_n^{-1} N_n \sup_{a\in(0,1-\vare_n)} |z_n(a)| < \infty , \\ \nonumber
  \limsup_{n\to\infty}
   N_n \int_{1-\vare_n}^1
        \left|1 - \ee^{\frac{\lambda}{1-a}z_n(a)}\right| (1-a)^\beta \, \dd a
   < \infty , \\ \nonumber
  \lim_{n\to\infty}
   \left|N_n \int_{1-\vare_n}^1
              \left(1 - \ee^{\frac{\lambda}{1-a}z_n(a)}\right) (1-a)^\beta \, \dd a
         - I\right|
  = 0
 \end{gather}
 with some \ $I \in \CC$.
\ Then
 \[
   \lim_{n\to\infty}
    N_n \int_0^1
         \left(1 - \ee^{\frac{\lambda}{1-a}z_n(a)}\right) \psi(a) (1-a)^\beta
         \, \dd a
   = \psi_1 I .
 \]
\end{Lem}

\noindent{\bf Proof.}
For all \ $a \in (0, 1)$ \ and for sufficiently large \ $n \in \NN$, \ we have
 \ $1 - \vare_n > a$, \ hence, by
 \eqref{help_ordo1},
 \begin{equation}\label{help_ordo2}
   N_n |z_n(a)| \leq \vare_n \vare_n^{-1} N_n \sup_{b\in(0,1-\vare_n)} |z_n(b)|
   \to 0 \qquad \text{as \ $n \to \infty$,}
 \end{equation}
 thus we conclude \ $\lim_{n\to\infty} N_n |z_n(a)| = 0$.
\ By applying \eqref{exp:1} and using \eqref{help_ordo2}, for any \ $n \in \NN$ \ and
 \ $a \in (0, 1)$, \ we get
 \begin{equation}\label{ordo1}
  \Bigl|N_n \left(1 - \ee^{\frac{\lambda}{1-a} z_n(a)}\right)\Bigr|
  \leq N_n\Bigl|\frac{\lambda}{1-a}z_n(a)\Bigr|
       \ee^{\left|\frac{\lambda}{1-a}z_n(a)\right|}
  \to 0 \qquad \text{as \ $n \to \infty$.}
 \end{equation}
If \ $n \geq n_0$ \ and \ $a \in (0, 1 - \vare_n)$, \ then \ $\frac{1}{1-a}<\vare_n^{-1}$ \ and
 \[
   \left|N_n\left(1 - \ee^{\frac{\lambda}{1-a}z_n(a)}\right)\right|
   \leq \lambda
        \left( \sup_{n\geq n_0} \vare_n^{-1} N_n \sup_{a\in (0,1-\vare_n)} |z_n(a)|\right)
        \ee^{\lambda \sup_{n\geq n_0} \vare_n^{-1}\sup_{a\in (0,1-\vare_n)} |z_n(a)|} =: C ,
 \]
 where \ $C \in \RR_+$ \ (due to \eqref{help_ordo1}).
\ Since \ $\int_0^1 \psi(a) (1-a)^{\beta} \, \dd a = 1$, \ we have
 \begin{align*}
  &\left| N_n\int_0^{1-\vare_n}
              \left(1 - \ee^{\frac{\lambda}{1-a}z_n(a)}\right) \psi(a)(1-a)^\beta
              \, \dd a\right|\\
  &\qquad = \left| \int_0^1
                 N_n\left(1 - \ee^{\frac{\lambda}{1-a}z_n(a)}\right) \bone_{(0,1-\vare_n)}(a) \psi(a)(1-a)^\beta
              \, \dd a\right|  \\
  &\qquad \leq \int_0^1 C \psi(a) (1-a)^\beta \, \dd a
   < \infty
 \end{align*}
 for \ $n\geq n_0$.
\ Therefore, \ $(0, 1) \ni a \mapsto C \psi(a)(1-a)^\beta$ \ serves as a dominating
 integrable function.
Thus, by the dominated convergence theorem, the pointwise convergence in \eqref{ordo1} results
 \begin{align}\label{help_ordo5}
   \lim_{n\to\infty}
    N_n \int_0^{1-\vare_n}
         \left(1 - \ee^{\frac{\lambda}{1-a}z_n(a)}\right) \psi(a) (1 - a)^\beta
         \, \dd a
   = 0 .
 \end{align}
Moreover, for all \ $ n \geq n_0$, \ we have
 \begin{align}\label{help_ordo6}
  \begin{split}
  &\left|N_n \int_0^1
              \left(1 - \ee^{\frac{\lambda}{1-a}z_n(a)}\right) \psi(a) (1-a)^\beta
              \, \dd a
        - \psi_1 I\right| \\
  &\qquad\qquad
   \leq \left|N_n \int_0^{1-\vare_n}
                   \left(1 - \ee^{\frac{\lambda}{1-a}z_n(a)}\right) \psi(a) (1-a)^\beta
                   \, \dd a\right| \\
  &\qquad\qquad\quad
        + \left|N_n \int_{1-\vare_n}^1
                     \left(1 - \ee^{\frac{\lambda}{1-a}z_n(a)}\right)
                     (\psi(a) - \psi_1) (1-a)^\beta
                     \, \dd a\right| \\
  &\qquad\qquad\quad
        + \psi_1 \left|N_n \int_{1-\vare_n}^1
                            \left(1 - \ee^{\frac{\lambda}{1-a}z_n(a)}\right)
                            (1-a)^\beta
                            \, \dd a
                       - I\right| ,
 \end{split}
 \end{align}
 where
 \begin{align*}
  &\left|N_n \int_{1-\vare_n}^1
              \left(1 - \ee^{\frac{\lambda}{1-a}z_n(a)}\right) (\psi(a) - \psi_1)
              (1-a)^\beta
              \, \dd a\right| \\
  &\leq \biggl(\sup_{a\in[1-\vare_n,1)} |\psi(a) - \psi_1|\biggr)
        N_n \int_{1-\vare_n}^1
             \left|1 - \ee^{\frac{\lambda}{1-a}z_n(a)}\right| (1-a)^\beta
             \, \dd a ,
 \end{align*}
 with \ $\sup_{a\in[1-\vare_n,1)} |\psi(a) - \psi_1| \to 0$ \ as \ $n \to \infty$, \ by
 the assumption \ $\lim_{x\uparrow 1}\psi(x) = \psi_1$.
\ Taking \ $\limsup_{n\to\infty}$ \ of both sides of \eqref{help_ordo6},
by \eqref{help_ordo5} and the assumptions of the lemma, we obtain the statement.
\proofend

\section*{Acknowledgements}
We would like to thank the referees for their comments that helped us improve the paper.


\bibliographystyle{plain}
\bibliography{aggr5}

\def\polhk#1{\setbox0=\hbox{#1}{\ooalign{\hidewidth
  \lower1.5ex\hbox{`}\hidewidth\crcr\unhbox0}}}
\begin{thebibliography}{10}

\bibitem{AloAlz}
M.~A. Al-Osh and A.~A. Alzaid.
\newblock First-order integer-valued autoregressive ({INAR}({$1$})) process.
\newblock {\em J. Time Ser.\ Anal.}, 8(3):261--275, 1987.

\bibitem{BarBasKevPapPla}
M.~{Barczy}, B.~{Basrak}, P.~{Kevei}, G.~{Pap}, and H.~{Planini{\'c}}.
\newblock Statistical inference of subcritical strongly stationary
  {G}alton--{W}atson processes with regularly varying immigration.
\newblock 2019.
\newblock To appear in Stochastic Process. Appl. See also ar{X}iv:1910.01420.

\bibitem{BarNedPap_Arxiv}
M.~Barczy, F.~Ned\'enyi, and G.~Pap.
\newblock Iterated scaling limits for aggregation of randomized {INAR}(1)
  processes with idiosyncratic {P}oisson innovations.
\newblock 2015.
\newblock ar{X}iv: 1509.05149. This is an extended version of [3].

\bibitem{BarNedPap}
M.~Barczy, F.~Ned\'enyi, and G.~Pap.
\newblock Iterated scaling limits for aggregation of randomized {INAR}(1)
  processes with idiosyncratic {P}oisson innovations.
\newblock {\em J. Math. Anal. Appl.}, 451(1):524--543, 2017.

\bibitem{DomKaj}
C.~Dombry and I.~Kaj.
\newblock The on-off network traffic model under intermediate scaling.
\newblock {\em Queueing Syst.}, 69(1):29--44, 2011.

\bibitem{Gra}
C.~W.~J. Granger.
\newblock Long memory relationships and the aggregation of dynamic models.
\newblock {\em J. Econometrics}, 14(2):227--238, 1980.

\bibitem{Jir}
M.~Jirak.
\newblock Limit theorems for aggregated linear processes.
\newblock {\em Adv.\ in Appl.\ Probab.}, 45(2):520--544, 2013.

\bibitem{LeiPhiPilSur}
R.~Leipus, A.~Philippe, V.~Pilipauskait\.{e}, and D.~Surgailis.
\newblock Sample covariances of random-coefficient {${\rm AR}(1)$} panel model.
\newblock {\em Electron. J. Stat.}, 13(2):4527--4572, 2019.

\bibitem{Li}
Z.~Li.
\newblock {\em Measure-valued branching {M}arkov processes}.
\newblock Probability and its {A}pplications ({N}ew {Y}ork). Springer,
  Heidelberg, 2011.

\bibitem{McK}
E.~McKenzie.
\newblock Some simple models for discrete variate time series.
\newblock {\em JAWRA Journal of the American Water Resources Association},
  21(4):645--650, 1985.

\bibitem{MikResRooSte}
T.~Mikosch, S.~Resnick, H.~Rootz{\'e}n, and A.~Stegeman.
\newblock Is network traffic approximated by stable {L}\'evy motion or
  fractional {B}rownian motion?
\newblock {\em Ann.\ Appl.\ Probab.}, 12(1):23--68, 2002.

\bibitem{NedPap}
F.~Ned{\'e}nyi and G.~Pap.
\newblock Iterated scaling limits for aggregation of random coefficient {AR}(1)
  and {INAR}(1) processes.
\newblock {\em Statist. Probab. Lett.}, 118:16--23, 2016.

\bibitem{PilSkoSur}
V.~Pilipauskait{\.e}, V.~Skorniakov, and D.~Surgailis.
\newblock Joint temporal and contemporaneous aggregation of random-coefficient
  {${\rm AR}(1)$} processes with infinite variance.
\newblock {\em Adv.\ in Appl.\ Probab.}, 52(1):237--265, 2020.

\bibitem{PilSur}
V.~Pilipauskait{\.e} and D.~Surgailis.
\newblock Joint temporal and contemporaneous aggregation of random-coefficient
  {AR}(1) processes.
\newblock {\em Stochastic Process.\ Appl.}, 124(2):1011--1035, 2014.

\bibitem{PupSur1}
D.~Puplinskait{\.e} and D.~Surgailis.
\newblock Aggregation of random-coefficient {AR}(1) process with infinite
  variance and common innovations.
\newblock {\em Lith.\ Math.\ J.}, 49(4):446--463, 2009.

\bibitem{PupSur2}
D.~Puplinskait{\.e} and D.~Surgailis.
\newblock Aggregation of a random-coefficient {${\rm AR}(1)$} process with
  infinite variance and idiosyncratic innovations.
\newblock {\em Adv.\ in Appl.\ Probab.}, 42(2):509--527, 2010.

\bibitem{Rob}
P.~M. Robinson.
\newblock Statistical inference for a random coefficient autoregressive model.
\newblock {\em Scand.\ J. Statist.}, 5(3):163--168, 1978.

\bibitem{SteHar}
F.~W. Steutel and K.~van Harn.
\newblock Discrete analogues of self-decomposability and stability.
\newblock {\em Ann.\ Probab.}, 7(5):893--899, 1979.

\bibitem{Str}
D.~W. Stroock.
\newblock {\em Probability {T}heory. An {A}nalytic {V}iew, 2nd edition}.
\newblock Cambridge University Press, Cambridge, 2011.

\bibitem{TaqWilShe}
M.~S. Taqqu, W.~Willinger, and R.~Sherman.
\newblock Proof of a fundamental result in self-similar traffic modeling.
\newblock {\em ACM SIGCOMM Computer Communication Review}, 27(2):5--23, 1997.

\bibitem{The}
H.~Theil.
\newblock {\em Linear Aggregation of Economic Relations. Contributions to
  Economic Analysis, VII}.
\newblock Amsterdam, North-Holland Publishing Company, 1954.

\bibitem{TurScoBer}
K.~F. Turkman, M.~G. Scotto, and P.~de~Zea~Bermudez.
\newblock {\em Non-linear time series. Extreme events and integer value
  problems}.
\newblock Springer, Cham, 2014.

\bibitem{Wei}
Ch.~H. Wei{\ss}.
\newblock Thinning operations for modeling time series of counts---a survey.
\newblock {\em AStA Adv.\ Stat.\ Anal.}, 92(3):319--341, 2008.

\bibitem{Zel}
A.~Zellner.
\newblock On the aggregation problem: a new approach to a troublesome problem.
  {I}n: {E}conomic models, estimation and risk programming: {E}ssays in honor
  of {G}erhard {T}intner. {E}dited by {K}. {A}. {F}ox, {J}. {K}. {S}engupta and
  {G}. {V}. {L}. {N}arasimham.
\newblock {\em Lecture Notes in Operations Research and Mathematical
  Economics}, 15, Springer-Verlag, Berlin-New York, 1969.

\bibitem{Zol}
V.~M. Zolotarev.
\newblock {\em One-{D}imensional {S}table {D}istributions}.
\newblock American Mathematical Society, Providence, RI, 1986.

\end{thebibliography}

\end{document}